# Grassmannian stochastic analysis and the stochastic quantization of Euclidean Fermions

Sergio Albeverio, Luigi Borasi,
Francesco C. De Vecchi, Massimiliano Gubinelli

Hausdorff Center for Mathematics &
Institute for Applied Mathematics
University of Bonn, Germany

*Email:* albeverio@iam.uni-bonn.de, borasi@iam.uni-bonn.de
francesco.devecchi@uni-bonn.de, gubinelli@iam.uni-bonn.de

**Abstract**

We introduce a stochastic analysis of Grassmann random variables suitable for the stochastic quantization of Euclidean fermionic quantum field theories. Analysis on Grassmann algebras is developed here from the point of view of quantum probability: a Grassmann random variable is an homomorphism of an abstract Grassmann algebra into a quantum probability space, i.e. a $C^*$-algebra endowed with a suitable state. We define the notion of Gaussian processes, Brownian motion and stochastic (partial) differential equations taking values in Grassmann algebras. We use them to study the long time behavior of finite and infinite dimensional Langevin Grassmann stochastic differential equations driven by Gaussian space-time white noise and to describe their invariant measures. As an application we give a proof of the stochastic quantization and of the removal of the space cut-off for the Euclidean Yukawa model.

# Table of contents

# 1 Introduction

Euclidean quantum field theories (QFTs) [167, 69, 141] are tools to construct and analyze mathematical models of relativistic quantum fields, that is the quantum theory of elementary particle in interaction satisfying the basic requirements of special relativity, i.e. Poincaré covariance and locality [85, 170, 98, 32, 82, 15, 14, 171, 63]. In the case of Bose–Einstein particles, the Euclidean theory is given by a probability measure on Schwartz distributions over the Euclidean space $\mathbb{R}^d$. In other words, bosonic Euclidean fields $\varphi$ are random distributions [173, 127, 27]. By a well–known result of Osterwalder and Schrader [69] (see also [187]), a basic set of properties, essentially Euclidean covariance, reflection positivity and regularity, is sufficient to reconstruct a well-behaved relativistic QFT from this probabilistic data.

Stochastic quantization (SQ) [95, 48] is an approach to the construction of the correlation functions of an Euclidean QFT introduced by Parisi and Wu [133]. The basic idea is to consider an additional variable (usually a fictitious time) and interpret the Euclidean fields $\varphi$ as the stationary solution to a stochastic partial differential equations (SPDEs) involving this additional variable and an external Gaussian source of noise $\xi$.

This strategy has been employed to rigorously study bosonic Euclidean theories starting with the work by Jona-Lasinio and Mitter [96, 97, 30] on the stochastic quantization of the $\Phi^4$ model for a scalar particle in two (Euclidean) dimensions and with quartic interaction. More recently the work of Hairer on regularity structures [84] opened the way to the study of the three dimensional $\Phi^4$ model, see also [42, 106] for other approaches. While the original implementation of stochastic quantization gives rise to parabolic SPDEs, variants can be constructed involving elliptic [3, 4, 17] or hyperbolic equations [80, 79].

From physicists' point of view, stochastic quantization gives an alternative approach to define and regularize quantum theories (especially theories with gauge invariance) [48, 30]. However, in the past few years, it has been realized that stochastic quantization has also interesting properties from the mathematical point of view. By solving the stochastic evolution it is possible to express the random fields $\varphi$ of a (bosonic) Euclidean theory as well behaved functionals of the external noise and this allows to study non-perturbative features of the former by leveraging the Gaussian structure of the latter. This new perspective led to a series of results on the global space-time control of the stochastic dynamics [123, 76, 6, 122] and to a novel proof of the constructions of non-Gaussian bosonic Euclidean quantum field theories in two and three dimensions, including the $\Phi_3^4$ model [77] (see also [7] for the construction of $\Phi_3^4$ model using rotation invariant cut-offs), the $\exp(\varphi)_2$ model [4], the $\sinh(\varphi)_2$ model [17] (first introduced in [5]). See the introductions to [6, 7] and to [77] for further references on stochastic quantization and its use in constructive QFT.

In the case of Fermi–Dirac particles (e.g. the electron), the Euclidean theory is expressed in terms of a linear functional over a Grassmann algebra [157, 158, 126, 130, 131]. The value of this functional on Grassmann monomials gives the correlation functions (Schwinger functions) of the Euclidean theory, which are represented by Berezin integrals [28, 29] of the type

$$\langle O(\psi,\bar\psi)\rangle = \frac{\int \mathrm{d}\psi \mathrm{d}\bar\psi\, O(\psi,\bar\psi) e^{-S_E(\psi,\bar\psi)}}{\int \mathrm{d}\psi \mathrm{d}\bar\psi\, e^{-S_E(\psi,\bar\psi)}}, \qquad (1)$$

where the fields $\psi, \bar\psi$ are the generators of the Grassmann algebra, $O$ is a functional of $\psi$ and $\bar\psi$ generating the Schwinger functions, and $S_E(\psi,\bar\psi)$ is the Euclidean action function, usually the sum of a quadratic part and a polynomial interaction (possibly involving also bosonic fields, which then will have to be averaged according to an appropriate probability measure) [61, 53, 118].

The main motivation of the present work is to develop a suitable setting for the stochastic quantization of Euclidean fermionic quantum field theories. In consideration of the Grassmannian nature of the fermionic Euclidean fields as expressed by (1), this must involve stochastic partial differential equations taking values in Grassmann algebras, including the study of their long time behavior, invariant measures and regularity properties. The main differences between the bosonic and the fermionic cases is the need of extending the notion of stochastic processes to a non-commutative framework. The approach we follow is to frame this problem in the context of *algebraic* probability (sometimes referred to as *quantum* or *non-commutative* probability) [134, 31, 121, 169, 1] and to construct the relevant Grassmannian objects as non-commutative random variables.

An algebraic probability space $(\mathcal{A}, \omega)$ is given by a $C^*$-algebra $\mathcal{A}$ and a state $\omega$, which is a linear normalized positive functional on $\mathcal{A}$. Inspired by the general approach of Accardi et al. [2] we will define a Grassmann random variable as a homomorphism of a Grassmann algebra $\Lambda$ into $\mathcal{A}$ (cf. Definition 1 below). We will not require the homomorphisms to respect the natural involution of $\mathcal{A}$, since there is no canonical candidate for that in Grassmann algebras.

It is useful to keep in mind how classical commutative random variables fit in the algebraic approach. Consider a random variable $X\colon \Omega \to \mathcal{M}$ on a probability space $(\Omega, \mathcal{F}, \mathbb{P})$, taking values in a manifold $\mathcal{M}$ and let $\mathcal{B} = L^\infty(\mathcal{M})$ be the algebra of all measurable, bounded, complex functions on $\mathcal{M}$. All the relevant probabilistic information about $X$ is encoded in the homomorphism of algebras $f \in \mathcal{B} \mapsto f(X) \in \mathcal{A}$ where $\mathcal{A}$ is the algebra of measurable, bounded, complex random variables on the basic space $\Omega$ endowed with the linear functional given by the expectation $\mathbb{E}$ associated to the probability measure $\mathbb{P}$. In particular, this is the point of view with which SDEs on manifolds are defined [88].

The embedding of $\Lambda$ into $\mathcal{A}$ allows to use the topology of $\mathcal{A}$ to do analysis on Grassmann algebras. Of course the associated analysis will not be canonical from the point of view of the Grassmann algebra itself, but it will turn out to be powerful enough to allow us to obtain a satisfactory theory of stochastic quantization for fermions. Let us mention another analogy which can help the reader to understand our point of view. The standard approach to study Gaussian processes in Hilbert spaces [71, 94, 86] requires to have a non-canonical embedding of the Hilbert space into a larger Banach space $B$ which supports the Gaussian measure. This Banach space is not canonical and there are various possible choices for any given Hilbert space $H$. For example Brownian motion is associated with the Hilbert space $H$ of functions on $\mathbb{R}_+$ with square integrable derivative, but Wiener measure is supported on the Banach space $C^\gamma$ of $\gamma$ Hölder functions for any $\gamma < 1/2$. In the case of random variables taking values in a Grassmann algebra $\Lambda$, the role of this "bigger" Banach space is played by a suitable Clifford sub-algebra of $\mathcal{A}$ which provides us with enough Grassmann generators to realize our homomorphisms and which has a canonical notion of norm. It will be for us just "a convenient place where to hang our (analytic) hat on".

The main finding of our work is that, once one accepts this point of view, the rest of the analysis falls in place quite naturally. Objects like Grassmann white noise, Brownian motion and free fields are relatively simple to define and control. On top of them a theory of stochastic differential equations (SDEs) with values in Grassmann algebras can be initiated and carried to the point to be able to discuss their long time evolution and invariant measures. Moreover many considerations carry over to SDE in infinite dimensions and allow to define and solve, in particular, the non-commutative nonlinear stochastic partial differential equations (SPDEs) appearing in stochastic quantization, again in relatively simple terms. We found that, in some respects, the analysis of Grassmann SPDEs is *simpler* than its commutative counterpart: the discussion of certain parts of the theory, like the existence of global space-time dynamics, are made relatively trivial by the boundedness of Grassmann Gaussian variables.

To test the effectiveness of our approach we give a proof of existence of the infinite volume limit of the massive Euclidean Yukawa model for small coupling and with an ultraviolet regularization (see e.g. [130, 131, 113] and the discussion of the literature below for additional references).

We have left open the substantial problem of removing the small-scale regularization in stochastic quantization. In order to make progress in this direction one will have to understand the renormalization problem for singular Grassmann SPDEs in our setting, along the line of the recent work in (commutative) singular SPDEs [84, 78, 106]. Apart from the more delicate analysis required by the low regularity of the random fields, the removal of these divergences will make appear renormalized Wick products of Euclidean fermion fields. Here we have in mind for example the local term $\bar\psi(x)\cdot\psi(x)$ appearing in the Yukawa$_2$ model (cf. (67) below) which is well known to require renormalization (Wick ordering) to be defined after removing the small-scale regularization. After renormalization, these random variables can presumably only be represented as unbounded operators acting in some underlying Hilbert space. This problem is similar to the one of definition of the stochastic integral in a Clifford algebra [18, 19, 20] which involves unbounded operators belonging to Segal's non-commutative $L^2$ space. We plan to investigate this problem in a future work.

**Structure of the paper.** In Section 2 we lay down the basic definition for Grassmann random variables and we construct the basic example of the Grassmann Gaussian together with corresponding notions of Gaussian white noise and Brownian motions. Section 3 studies finite dimensional Grassmann differential equations with additive Gaussian white noise, addressing the problem of existence, uniqueness, long time solutions and invariant measures. Section 4 is devoted to extend the previous results to SDE in infinite dimensions of the kind relevant for stochastic quantization. In Section 5 we present our main application to the stochastic quantization of the Yukawa model. The model is introduced in Section 5.1 in two Euclidean dimensions and with ultraviolet (UV) and infrared (IR) cut-offs, following the construction in [131]. In Section 5.2 the corresponding stochastic quantization equation (SQE) is formulated and studied and we formulate our main result, Theorem 60 were we formulate the stochastic quantization of the finite volume model and subsequently use it to prove the existence of the infinite volume limit. The proof uses a series of approximations discussed in Section 5.3. In Appendix A we include some side results about the convergence of series expansions for the solutions to finite dimensional Grassmann SDEs, while in Appendix B we collect various technical results about functional spaces of Banach-valued functions. Finally, Appendix C contains the proof of convergence of perturbative series in the infinite volume Yukawa$_2$ model.

**Acknowledgements.** This work is supported by DFG via the grant AL 214/50-1 "Invariant measures for SPDEs and Asymptotics". The authors are funded by the DFG under Germany's Excellence Strategy - GZ 2047/1, project-id 390685813. M. G. would like to thank Hao Shen for comments on an earlier version of the paper.

## 1.1 Relations with previous works

This part of the introduction is dedicated to discuss the connections of our work with previous research.

**Non-commutative analysis and quantum physics.**

The interest in studying time evolutions and stochastic processes with values in non commutative algebras arose in connection with the development of quantum mechanics in the 20s–30s of last century and later in quantum field theory, particularly around the 50s. In these theories one associates physical observables with operators evolving in time and acting on an Hilbert space. Mathematically, this lead to the representation theory of certain operator algebras (CCR-algebras resp. CAR-algebras in the case of bosons resp. fermions). In this line, quantum fields require the study of certain infinite dimensional (non-commutative) algebras, the study of which has been particularly stimulated by developments in the period 1940–1960, where the mathematical description of systems with an infinite number of degrees of freedom, in particular for quantum fields and statistical mechanics has been worked out. Let us mention the first works on von Neumann algebras by Murray and von Neumann (see e.g. [180]) and the work on $C^*$-algebras developed originally and quite independently by I. M. Gelfand and his school (see, e.g., [67]), since the beginnings of the 40s, and by I. E. Segal and his school (see e.g. [160, 15]). The mathematical literature on Banach, $C^*$ and $W^*$ algebras is very rich and there are a number of excellent monographs on them, see e.g., Naimark [125] (see also, e.g. [34]), and Zelazko [185] for normed spaces and algebras, and Kadison and Ringrose [100, 99] (and references therein) for the other type of algebras. For applications in mathematical physics, in particular quantum mechanics, quantum field theory and statistical mechanics, there are monographs by Bratteli and Robinson [37, 38], Emch [55], Baumgärtel [22], as well as Baumgärtel and Wollenberg [23], Holevo [87], Segal [162], Schmüdgen [156]; for relativistic quantum fields, we refer particularly to the books by Araki [14], by Baez, Segal and Zhou [15] and Haag [82]. The articles by Summers [172] and Doplicher [16] survey particularly some recent developments related to application of such algebras in quantum field theory.

**Grassmann algebras and Euclidean fermions.**

The algebraic features of a quantum theory change drastically as we perform the analytic continuation of the correlations functions to imaginary time, which is the main idea underlying the Euclidean method. In particular theory with bosons, which satisfy the canonical commutation relations (CCR) are described, in the Euclidean domain, with a standard commutative probability measures [173]. On the other hand, already in the pioneering work of Schwinger and Nakano [157, 158], it was observed that theories with fermions which satisfy canonical anti-commutation relations (CAR) are described in the Euclidean domain by Grassmann algebras. Grassmann algebras constitute a class of associative non commutative algebras that were introduced by Grassmann in the second half of the 19th century.

Osterwalder and Schrader [130, 131] were the first to rigorously show that Euclidean fermions satisfy the relations defining a Grassmann algebra. They constructed an embedding of this Grassmann into a $C^*$ algebra of bounded operators which allowed them to consider rigorously constructions involving infinite dimensional Grassmann algebras, e.g. an Euclidean Feynman–Kac formula from which eq. (1) can be derived. See the review of Palmer [132] for a relatively recent account of their theory. Their approach is not unique and there are various alternative Euclidean theories for fermions [179, 182, 81, 110, 35] proposed in the literature. We stick to the OS approach being simple and effective: our setting for Grassmann random variables is a reformulation of the OS construction in the language of Accardi et al. [2].

After the seminal contribution of OS, it seems however that the starting point of the majority of the subsequent work with Euclidean fermions has been the definition of the Schwinger functions via eq. (1) as formulated in the formalism introduced by Berezin's in [29]. Berezin integrals have been useful in calculating fermionic systems in the theoretical physics literature, see e.g. Zinn-Justin [186], Izkyson and Zuber [91]. Berezin's formalism has been subsequently developed into an analysis on Grassmann algebras, in the form of super-analysis [29]. For super-algebras and super-analysis see also, e.g., Rogers [146] and Leites [111], on a more heuristic level also DeWitt [51]. For a review see Pestov [137].

Such methods are also quite essential in most of the rigorous mathematical literature on fermionic quantum fields, including renormalization group methods, see, e.g., Mastropietro [118], Feldman, Knörrer and Trubowitz [57] or Salmhofer [149]. Euclidean fermionic quantum field theories involving fermions have surprisingly simpler behavior than the bosonic ones [40] and have been very much studied also from the rigorous point of view. In particular removal of both the UV and IR cutoffs has been obtained in various models either by using phase-cell expansions à la Glimm-Jaffe or via renormalization group methods. Without the ambition to gather an exhaustive list of all the relevant literature, let us cite some of the most important works in this direction for various kind of models.

a) Since the main example treated in this work the Yukawa model in two dimensions, let us start from it. It couples a scalar (or pseudo-scalar) massive boson field with a pair of Fermi fields. Historically it was introduced by Yukawa in 1935 and further studied e.g. in [119]. Its mathematical construction was achieved, in space-time dimension 2, as part of the constructive approach to relativistic quantum fields by different methods. The first results have been obtained by Hamiltonian and Minkowski space methods: after a pioneering work of Lanford III [108] ultraviolet (UV) and infrared (IR) cut-offs, the removal of the UV cut-off was realized in [120] and the complete construction obtained by Glimm and Jaffe [68]. After, with the advent of the Euclidean approach Seiler [163] investigated the removal of the UV cutoff using the Euclidean Matthew–Salam formula [119], then Cooper and Rosen [45] and Magnen and Sénéor [114] proved Lorenz invariance and the Wightman axioms by removing both cutoffs and Renouard [140] proved the Borel summability of the perturbation series for the Schwinger functions. Another reference containing a discussion of the Yukawa model in 2 dimensions by renormalization group (RG) methods in Lesniewski [113] where an effective action is constructed and Borel summability is proven. Bonetto and Mastropietro [33] used RG methods to study the limit of small mass (or large coupling). The Yukawa model has also been discussed in 3 dimensions by Magnen and Sénéor [115].

b) Feldman, Magnen, Rivasseau and Sénéor [61] constructed the more singular case of the Gross-Neveu model in 2 dimensions ($GN_2$) with Euclidean methods, for small coupling and also verifying the Osterwalder and Schrader axioms. Similar results are in Gawędzki and Kupiainen [66]. Both groups used renormalization group methods which seems the only approach possible in case of just-renormalizable theories like the $GN_2$ model. Let us note that Disertori and Rivasseau [53] discuss a continuous version of the RG method for $GN_2$ following a suggestion by Salmhofer. More precisely in that paper the RG equations are solved and Borel summability of the perturbation theory is proved. See also [152] for an epsilon-deformed version of the model.

c) The Thirring model has been constructed by Benfatto, De Falco and Mastropietro in [25]. The Luttinger model (a Thirring model with UV cut-off) is treated by RG methods in [117]. RG methods have also been used to study quantum many body problems, see e.g. De Roeck and Salmfhofer [49]. Feldman, Knörrer and Trubowitz [62, 56, 58, 59] developed an impressive extensive analysis of the RG method for fermionic non-relativistic systems at finite temperature where the singularities are due to the slow decay of the correlation functions.

d) More generally we refer to the books [141, 149, 118] and [57] for RG methods leading to removal of the cut-offs in Euclidean models (see also [26] for further references).

**Stochastic quantization of Euclidean fermions.**

After Parisi and Wu seminal paper, the stochastic quantization of fermions has been discussed heuristically in the physics literature, starting with Kakudo et al. [101], Fukai et al. [65], Damgaard and Tsokos [47] and Xue and Hsien [184], up to the recent paper of Efremov [54]. All these authors stress particularly the relations with "Grassmann valued random variables", solving formally partial differential equations with Grassmann valued Gaussian white noise. The invariant measures of physical interest are here described by averages with respect to Berezin integrals. However all these papers never properly discuss the analytical difficulties of considering (stochastic) differential equations in Grassmann algebras and the related stochastic analysis.

Grassmann variables are in a certain sense *completely non-commutative numbers*: as such they appear also in non-commutative central limit theorems [181] and in the setting of non-commutative processes with independent increments [159, 36]. From the point of view of probability theory, pioneering work in the general theory of non-commutative stochastic processes has been done by Accardi, Streater and Hudson and Parthasarathy and others in the 70s–80s, see e.g. [1, 134, 169, 31]. This includes the study of non-commutative Markov semigroups, non-commutative Brownian motion and the associated stochastic calculus. Let us also mention the work of Gross [72, 73, 74, 75] applying Segal's non-commutative integration theory [161, 128] to Clifford algebras.

Non-commutative stochastic calculus has been developed, with various degrees of completeness in the setting of the bosonic and fermionic Fock space and for Clifford algebras in the works of Barnett, Streater and Wilde [20, 19, 18], Applebaum and Hudson [12, 13], Belavkin [24], Carlen and Kree [41], Gordina [70], Sinha and Goswami [168], Kümmer [105]. See also the monographs by Meyer [121], Holevo [87], and the recent survey article by Cipriani [43] and the references therein. In these non-commutative settings is quite difficult to give general constructions of stochastic integrals which match in completeness with the commutative theory. Many technical difficulties stem from the fact that stochastic integrals are represented by unbounded operators and very delicate issues of domain enters the picture and in general restrict the applicability of these theories to study general SDE with non-linear coefficients.

In spite of all these developments, and to our surprise, none of the above works addresses the specific problems of the stochastic analysis in Grassmann algebras. The only relevant previous contributions in the literature that we managed to track down are the following:

a) de Angelis et al. [11] discuss a probabilistic representation for finite-dimensional Grassmann-valued Markov processes using (commutative) Poisson processes.

b) Rogers' approach [142, 143, 144, 145, 112] essentially consists in looking at all the finite dimensional sub-algebras, use the Berezin integral to compute averages and require certain natural consistency conditions in order to obtain a projective system. A similar line of research was also carried on by Kupsch and Haba [107, 83]. These last works follows the observation of Hudson and Parthasarathy [89] and Le Jan [109] (see also Meyer's book [121]) which developed a unified representation of the Bosonic and Fermionic Brownian motion (and related stochastic calculus) using Fock space techniques. However none of these lines of research seems to have reached a stage where the theory is powerful enough to easily and naturally accommodate SPDEs (or even SDEs) of the kind needed in stochastic quantization.

c) A rigorous approach to define Grassmann stochastic quantization was undertaken in a series of papers by Scherbakov [166, 165, 164] using the locally convex topology obtained by considering all the possible correlation functions of the random variables as the family of semi-norms. Convergence in this setting therefore corresponds to convergence of all the correlation functions. This approach has the disadvantage of not allowing the introduction of norms strong enough to study differential equations via standard analytic approaches. Despite this technical shortcoming, which makes the constructions quite indirect, it should be said that these papers define and successfully solve equations with finitely many degrees of freedom (i.e. with small scale and large scale cut-offs) and then address the convergence when the large scale cut-off is removed (infinite volume limit) using cluster expansion techniques. Other results in this direction were obtained later on by Ignatyuk et al. [90], with applications in statistical mechanics. In particular these results hints to the fact that a complete theory of stochastic analysis of Grassmann random variables *should be possible*.

This review of the literature on the various attempts to Grassmann stochastic analysis shows that the main difficulty is the identification of sufficiently flexible analytical tools for non-linear analysis in infinite dimensional Grassmann algebras.

For a mathematical discussion of the problem of analysis in Grassmann algebras and in particular about the introduction of suitable norms see, e.g. [93, 136, 92, 139, 8]. The basic difficulties of a probabilistic approach to Euclidean fermions has been also discussed by Fröhlich and Osterwalder [64] (Theorem 5.7) just after the foundational paper of Osterwalder and Schrader. Let us remark however that, despite these difficulties, norms on Grassmann Gaussian random variables can be naturally introduced via the corresponding covariance operators and have been used to prove interesting constructive results e.g. on Grassmann-valued effective actions in models, see e.g. [57, 60] (and references therein) and the works of Salmhofer et al. [151, 150] where considerable simplifications of remainders estimates in perturbation expansions have been achieved via these tools.

As we already mentioned, the approach we follow in this paper is inspired by Osterwalder and Schrader [130, 131] who were the first to do rigorous analysis of Euclidean fermions. A key aspect of their construction is that Grassmann Gaussian variables can be realized within an $C^*$ algebraic setting, i.e. as bounded operators acting on an Hilbert space. Salmhofer used similar techniques for interacting non-relativistic fermions at finite temperature, see e.g. Appendix B of the book [149], and the paper [135]. This approach has also been recently used to discuss supersymmetry in the context of Parisi–Sourlas dimensional reduction [104, 103] and Euclidean stochastic quantization [3, 4, 50] and, as we already mentioned, is the basis of the approach we use in the present work.

Our work provides non-trivial examples of non-commutative non-linear SPDEs and of their qualitative analysis. The analysis of partial differential equations in non-commutative algebras is not well developed in general. Minkowski QFTs provide examples of PDEs in non-commutative algebras [15]. Other relevant examples we are aware of are the work of Rosenberg [147] on the (linear) Laplace equation on the non-commutative torus and the subsequent work on linear operators and functional spaces on non-commutative spaces (e.g. [153, 183]), works on the algebraic and geometric features in PDEs (see e.g. [138, 44]), the work of Khrennikov on linear equations and differential operators on superspace [102] and the work of Osipov [129] on solutions of quantum field equations via Wick kernels. However none of these works concern SPDEs and all rather use indirect methods to find solutions. Let us mention however that Dabrowski [46] recently studied an interesting example of a non-commutative SPDE driven by multiplicative free noise.

## 2 Grassmann random variables

We denote by $\mathcal{L}(A, B)$ the space of linear maps between vector spaces $A$ and $B$, if both $A, B$ have topologies we consider all the maps to be continuous. We let $\mathcal{L}(A) = \mathcal{L}(A, A)$. With $\mathrm{Hom}(A, B)$ we denote homomorphisms between algebras $A$ and $B$.

### 2.1 Grassmann probability

*Algebraic probability.* We consider a complex Hilbert space $\mathcal{H}$ (to be fixed later) and denote by $\mathcal{A} = \mathcal{L}(\mathcal{H})$ the $C^*$-algebra of bounded operators with the operator norm. Moreover we assume to have a state $\omega$ on $\mathcal{A}$, i.e. a continuous linear map $\omega \colon \mathcal{A} \to \mathbb{C}$ which is positive, that is $\omega(A^*A) \geqslant 0$ for all $A \in \mathcal{A}$ and normalized by $\omega(\mathbb{I}) = 1$ where $*$ denotes the conjugation of the $C^*$-algebra $\mathcal{A}$ and $\mathbb{I} \in \mathcal{A}$ is the unit of $\mathcal{A}$. The pair $(\mathcal{A}, \omega)$ is a (non-commutative or algebraic) probability space. We do not require the state $\omega$ to be either faithful, or tracial.

*Grassmann algebras.* Let $V$ be a (finite dimensional or infinite dimensional and separable) real Hilbert space. Denote by $\Lambda V$ the Grassmann algebra generated by $V$, i.e. the exterior algebra obtained by quotienting the tensor algebra $T(V) = \oplus_{n \geqslant 0} V^{\otimes n}$ by the two-sided ideal generated by elements of the form $x \otimes x$ for $x \in V$. For $V$ finite dimensional, $\Lambda V$ is also finite dimensional with dimension, as a vector space, equal to $2^{\dim V}$. We denote the product of two elements $f, g \in \Lambda V$ by $f \wedge g$ (or if there is no confusion with other products simply by $fg$).

Let $(v_\alpha)_{\alpha \in I}$ be a fixed basis of $V$, where $I$ is a suitable index set. As a consequence $\Lambda V$ is spanned by elements of the form $v_A := v_{\alpha_1} \wedge \cdots \wedge v_{\alpha_n}$ and $A = (\alpha_1, ..., \alpha_n)$ is a $n$-tuple of elements of $I$ with the convention that $v_\emptyset = 1$ is the unit element in $\Lambda V$. When it is clear from the context we will denote the product in $\Lambda V$ simply by $v_A = v_{\alpha_1}...v_{\alpha_n}$. We recall that $\Lambda V$ is $\mathbb{Z}_2$ graded in the sense that it splits into odd and even parts $\Lambda V = \Lambda_{\mathrm{odd}} V \oplus \Lambda_{\mathrm{even}} V$. On $\Lambda V$ there is a super Hopf algebra structure with coproduct $\Delta: \Lambda V \to \Lambda V \otimes \Lambda V$ where $\Lambda V \otimes \Lambda V$ is the $\mathbb{Z}_2$-graded tensor product algebra which satisfies, for homogeneous elements $g, h \in \Lambda_{\mathrm{odd}} V \sqcup \Lambda_{\mathrm{even}} V$,

$$(f \otimes g)(h \otimes k) = (-1)^{|g||h|} fh \otimes gk, \qquad f, k \in \Lambda V,$$

with $|\cdot|: \Lambda_{\mathrm{odd}} V \sqcup \Lambda_{\mathrm{even}} V \to \{0, 1\}$ the even/odd grading. The coproduct $\Delta$ is the algebra homomorphism such that $\Delta v = 1 \otimes v + v \otimes 1$ for all $v \in V \subseteq \Lambda V$ and the counit $\varepsilon: \Lambda V \to \mathbb{C}$ given by $\varepsilon(v_A) = \mathbb{I}_{A = \emptyset}$.

*Random variables.*

**Definition 1.** *A $V$-Grassmann random variable $\Psi$ is an algebra homomorphism from the Grassmann algebra $\Lambda V$ into $\mathcal{A}$. We denote by $\mathcal{G}(V) = \mathrm{Hom}(\Lambda V, \mathcal{A})$ the set of all such homomorphisms. We call the law of $\Psi \in \mathcal{G}(V)$ the family of its moments $\omega^\Psi(F) := \omega(\Psi(F))$ for all $F \in \Lambda V$, also represented by the linear functional $\omega^\Psi: \Lambda V \to \mathbb{R}$.*

Note that $\Psi \in \mathcal{G}(V)$ cannot be assumed to be a $*$-algebra-homomorphism since $\Lambda V$ has no (natural) $*$-operation. If $F \in \Lambda V$ has representation $F = \sum_A F_A v_A$ we shall employ the following dual notation

$$F(\Psi) := \Psi(F) = \sum_A F_A \Psi^A,$$

where $\Psi^A := \Psi(v_A)$ and $F_A \in \mathbb{R}$. Since $\Psi$ is assumed to be an algebra homomorphism, we have e.g. $\Psi^\alpha \Psi^\beta = -\Psi^\beta \Psi^\alpha$, where $\Psi^\alpha = \Psi(v_\alpha)$ and $(v_\alpha)_\alpha$ is a fixed basis of $V$. Moreover, $\Psi(\Lambda V)$ is a Grassmann sub-algebra of $\mathcal{A}$ and $\Psi(F)$ has the same degree (even or odd) as $F$. As shown in [177, 178], even if some arguments are formulated in a basis dependent way, the definition of $\Lambda V$ and its characterization by anti-commutation relations is independent of the basis.

When the context is clear we will abbreviate $X \in \mathcal{G}(V)$ as $X \in \mathcal{G}$.

**Definition 2.** *Let $X \in \mathcal{G}(V)$ and $Y \in \mathcal{G}(W)$, we say that they are compatible if the linear map $Z: V \oplus W \to \mathcal{A}$ given by $Z(v) = X(v)$ if $v \in V$ and $Z(w) = Y(w)$ if $w \in W$, extends to a homomorphism $Z: \Lambda(V \oplus W) \to \mathcal{A}$. In this case we write $Z \in \mathcal{G}(V \oplus W)$ or briefly $(X, Y) \in \mathcal{G}$. Compatibility can of course be defined for $X_1, ..., X_n \in \mathcal{G}$ in a corresponding way. We shall express that $X_1, ..., X_n \in \mathcal{G}$ are compatible by writing $(X_1, ..., X_n) \in \mathcal{G}$.*

Note that we have a super- algebra isomorphism $\Lambda(V \oplus W) \approx \Lambda(V) \otimes \Lambda(W)$ where tensor product is in the sense of super-algebras. From this isomorphism we get in particular $Z(F) = m_{\mathcal{A}}[(X \otimes Y)(F)]$ for all $F \in \Lambda(V \oplus W)$ where $m_{\mathcal{A}}: \mathcal{A} \otimes \mathcal{A} \to \mathcal{A}$ denotes the multiplication of $\mathcal{A}$.

**Remark 3.** The notion of compatibility which we introduce here is not standard in algebraic probability. It ensures that one can extend multiple Grassmannian random variables to a well-defined "multi-component" Grassmannian random variable. In particular it encodes the fact that different component properly anti-commute. A related concept is that of kinematic independence, see e.g. the review of Accardi [1].

Given $X, Y \in \mathcal{G}(V)$ which are compatible we define $X + Y \in \mathcal{G}(V)$ as

$$(X+Y)(F) = F(X+Y) := m_{\mathcal{A}}^2[(X \otimes Y)\Delta F], \qquad F \in \Lambda V,$$

where $m_{\mathcal{A}}^k: \mathcal{A} \otimes \mathcal{A} \to \mathcal{A}$, $k \in \mathbb{N}$, is defined as

$$m_{\mathcal{A}}^k(a_1 \otimes \cdots \otimes a_k) = a_1 \cdots a_k. \tag{2}$$

By compatibility $Z = m_{\mathcal{A}}^2 \circ (X \otimes Y): \Lambda(V \oplus V) \approx \Lambda V \otimes \Lambda V \to \mathcal{A}$ is an algebra homomorphism and therefore

$$(X+Y)(FG) = Z(\Delta(FG)) = Z(\Delta F \Delta G) = Z(\Delta F)Z(\Delta G) = (X+Y)(F)(X+Y)(G).$$

The notation is justified from the fact that $(X+Y)(v) = X(v) + Y(v)$ for all $v \in V$. Similarly for any $\lambda \in \mathbb{C}$ we can define $\lambda X$ as the only homomorphism such that $(\lambda X)(v) = \lambda X(v)$ for all $v \in V$.

**Definition 4.** *If $(X_1, ..., X_n) \in \mathcal{G}(V_1 \oplus \cdots \oplus V_n)$ are compatible Grassmann variables with values in the probability space $(\mathcal{A}, \omega)$, then we say that $X_1, ..., X_n$ are (tensor) independent (with respect to the state $\omega$) if, for all $F_j \in \Lambda V_j$, we have that [52]*

$$\omega\left(\prod_{j=1}^k X_j(F_j)\right) = \prod_{j=1}^k \omega(X_j(F_j)).$$

For example if $(X, Y) \in \mathcal{G}(V \oplus W)$ are (compatible and) independent Grassmann random variables with values in $(\mathcal{A}, \omega)$, then for all $v \in V \approx V \oplus 0$, $w \in W \approx 0 \oplus W$ we have $vw \approx ((v \oplus 0) \wedge (0 \oplus w)) \in \Lambda(V \oplus W) \approx \Lambda V \otimes \Lambda W$ and

$$\omega((X \otimes Y)(vw)) = \omega(X(v)Y(w)) = \omega(X(v))\omega(Y(w)).$$

**Lemma 5.** *We can always arrange to realize two Grassmann variables on the same probability space in such a way that they are compatible and independent while preserving their respective laws.*

**Proof.** First of all lets observe that as long as we are interested only in the law of a Grassmann random variable $X \in \mathcal{G}(V)$, we can assume that it is defined on $\mathcal{A} = \mathcal{L}(\mathcal{H})$ for some Hilbert space $\mathcal{H}$ with a vector state $\omega(\cdot) = (\Omega, \cdot \Omega)$ and that it comes with an involution $R_X \colon \mathcal{H} \to \mathcal{H}$ such that $R_X X(F) = X((-1)^{|F|}F) R_X$ for all $F \in \Lambda V$. Indeed by the GNS construction we can construct a Hilbert space $\mathcal{H}$ and a vector $\Omega \in \mathcal{H}$ such that there is $\tilde{X} \in \mathrm{Hom}(\Lambda V; \mathcal{L}(\mathcal{H}))$ such that $\omega(X(F)) = \langle \Omega, \tilde{X}(F)\Omega\rangle_{\mathcal{H}}$. By restriction, we can always redefine $\mathcal{H}$ to be $\mathcal{H} = \overline{\tilde{X}(\Lambda V)\Omega}$ while still keeping this relation. Here the bar denotes closure with respect to the Hilbert space topology of $\mathcal{K}$. Then, on $\mathcal{K}$, we can introduce by density an involution $R$ that acts as $R\tilde{X}(F)\Omega = \tilde{X}((-1)^{|F|}F)\Omega$ for homogeneous elements $F$, where recall that $|\cdot|$ is the grading on $\Lambda V$.

Now for $i = 1, 2$, consider the Grassmann variable $X_i \colon \Lambda V_i \to \mathcal{A}_i$ defined on the probability spaces $(\omega_i, \mathcal{A}_i = \mathcal{L}(\mathcal{H}_i))$ with involution $R_i = R_{X_i}$ and vector state $\omega_i = (\Omega_i, \cdot \Omega_i)$. Let $\mathcal{H} = \mathcal{H}_1 \otimes \mathcal{H}_2$, and $\omega \colon \mathcal{L}(\mathcal{H}) \to \mathbb{C}$ be given by $\omega(x) := (\Omega_1 \otimes \Omega_2, x\, \Omega_1 \otimes \Omega_2)_{\mathcal{H}}$, $x \in \mathcal{L}(\mathcal{H})$. Moreover we define the random variables $\tilde{X}_i \colon \Lambda(V_1 \oplus V_2) \to \mathcal{L}(\mathcal{H})$, obtained by extending the relations $\tilde{X}_1(v) = X_1(v) \otimes R_2$ and $\tilde{X}_2(w) = \mathbb{I}_{\mathcal{H}_1} \otimes X_2(w)$ for $v \in V_1, w \in V_2$. Then $\tilde{X}_1(v)\tilde{X}_2(w) = -\tilde{X}_2(w)\tilde{X}_1(v)$ and since $\Lambda(V_1 \oplus V_2) \approx \Lambda V_1 \otimes \Lambda V_2$ we have also that the map $\tilde{X}(v \oplus w) = \tilde{X}_1(v \oplus 0) + \tilde{X}_2(0 \oplus w)$ with $v, w \in V$ can be extended to a homomorphism from $\Lambda(V_1 \oplus V_2)$ to $\mathcal{A}$. Hence the variables $\tilde{X}_i$, $i = 1, 2$, are compatible. Moreover they have the same law as $X_i$ because $\omega(\tilde{X}_i(F_i)) = \omega_i(X_i(F_i))$ for all $F_1 \in \Lambda V_1, F_2 \in \Lambda V_2$, $i = 1, 2$. Finally they are independent: $\omega(\tilde{X}_1(F_1)\tilde{X}_2(F_2)) = \omega_1(X_1(F_1))\omega_2(X_2(F_2))$. □

## 2.2 Topology and calculus on Grassmann variables

Before beginning we want to specify the topology we consider on $\mathcal{G}(V)$. We consider here only the case where $V$ is finite dimensional, since there is a more or less unique natural topology on $\mathcal{G}(V)$. Some topologies of $\mathcal{G}(V)$ when $V$ is infinite dimensional are discussed in Section 4.

When $V$ is finite dimensional all norms on $V$ are equivalent. For this reason we choose the norm induced by the pre-Hilbert space inner product $\langle \cdot, \cdot \rangle$ (that in the case of $V$ finite dimensional is a Hilbert space inner product on $V$) related to the construction of $V$-valued Gaussian random variables $X$ (see Section 2.3). In this case $\mathcal{G}(V)$ has a natural metric topology given by the distance

$$d_{\mathcal{G}(V)}(X, Y) := \|X - Y\|_{\mathcal{G}(V)} = \sup_{v \in V, |v|_V = 1} \|X(v) - Y(v)\|_{\mathcal{A}}, \tag{3}$$

where $\|\cdot\|_{\mathcal{A}}$ is the norm of the $C^*$-algebra $\mathcal{A}$.

**Remark 6.** In principle the definition of the distance (3) and the related norm $\|\cdot\|_{\mathcal{G}(V)}$ does not use in any way the fact that $V$ is finite dimensional. For this reason in the following, when $V$ is a (in general infinite dimensional) pre-Hilbert space we will use the notation $\|X\|_{\mathcal{G}(V)}$ for the quantity $\sup_{v \in V, |v|_V = 1} \|X(v)\|_{\mathcal{A}}$.

When $V$ is finite dimensional $\mathcal{G}(V)$ is a complete metric space with respect to $d_{\mathcal{G}(V)}$, in fact we have:

**Lemma 7.** *The metric $d_{\mathcal{G}(V)}$ makes $\mathcal{G}(V)$ a complete metric space.*

**Proof.** The quantity $d_{\mathcal{G}(V)}(X,Y)$ as defined in (3) satisfies the usual properties of a metric. If $(X_n)_n$ is a Cauchy sequence with respect to $d_{\mathcal{G}(V)}$, there is an element $X \in \mathcal{L}(V, \mathcal{A})$ such that $(X_n)_n$ converges to $X$ in $\mathcal{L}(V, \mathcal{A})$. The only thing that we have to prove is that the linear map $X \in \mathcal{L}(V, \mathcal{A})$ can be extended to an homomorphism from $\Lambda V$ to $\mathcal{A}$. This is equivalent to prove that for any $v, w \in V$ we have $X(v)X(w) = -X(w)X(v)$ where the product is the natural one in $\mathcal{A}$. On the other hand by the continuity of the product of $\mathcal{A}$ with respect to $\|\cdot\|_{\mathcal{A}}$ we have

$$X(v)X(w) = \lim_{n\to+\infty} X_n(v)X_n(w) = -\lim_{n\to+\infty} X_n(w)X_n(v) = -X(w)X(v). \qquad \square$$

**Remark 8.** The distance $d_{\mathcal{G}(V)}$ is not the unique reasonable choice in $\mathcal{G}(V)$ and it depends on the topology chosen on $\mathcal{A}$ (in this case we choose the topology of the uniform converge of operators). An essential choice to preserve the non-linear structure of $\mathcal{G}(V)$, exploited in the proof of Lemma 7, seems to be that the product on $\mathcal{A}$ is continuous with respect to this topology.

In the following we will denote $\|Y\|_{\mathcal{G}(V)}$ for $Y \in \mathcal{G}(V)$ the positive number

$$\|Y\|_{\mathcal{G}(V)} := \|Y\|_{\mathcal{L}(V,\mathcal{A})} = \sup_{v \in V, \|v\|=1} \|Y(v)\|_{\mathcal{A}},$$

i.e. the norm of the restriction of $Y$ to $\mathcal{L}(V, \mathcal{A})$. This is useful to formulate a Taylor formula on $\mathcal{G}(V)$. We will use the simpler notation $m$ for the multiplication $m_{\mathcal{A}}$ in $\mathcal{A}$.

Consider $G \in \Lambda V$ and define the *right derivative* $\partial_R \colon \Lambda V \to \Lambda V \otimes V$ by

$$\partial_R G := (\mathbb{I} \otimes \Pi_V)(\Delta G) \colon \Lambda V \to \Lambda V \otimes V \qquad (4)$$

where $\Pi_V \colon \Lambda V \to V$ is the projection from the tensor algebra $\Lambda V$ onto $V$. Then

$$G(X+Y) - G(X) - m_{\mathcal{A}}^2[(X \otimes Y)(\partial_R G)] = m_{\mathcal{A}}^2[(X \otimes Y)(\mathbb{I} - \mathbb{I} \otimes \varepsilon - \mathbb{I} \otimes \Pi_V)\Delta G]$$

for $G \in \Lambda V$, where we used that $G(X) = m_{\mathcal{A}}^2[(X \otimes Y)(\mathbb{I} \otimes \varepsilon)\Delta G]$. Given that

$$(\mathbb{I} - \mathbb{I} \otimes \varepsilon - \mathbb{I} \otimes \Pi_V)\Delta G \in \Lambda V \otimes \Lambda_{\geqslant 2} V$$

where $\Lambda_{\geqslant 2} V$ denotes the subspace of $\Lambda V$ of elements of degree $\geqslant 2$. We can define recursively the $k+1$-th derivative as

$$\partial_R^{k+1} = \left( \partial_R \otimes \underbrace{\mathbb{I} \otimes \cdots \otimes \mathbb{I}}_{k} \right) \partial_R^k \colon \Lambda V \to \Lambda V \otimes \underbrace{V \otimes \cdots \otimes V}_{k}$$

where $\partial_R^k$ is the $k$-th order derivative. Note that the right derivative $\partial_R \colon \Lambda V \to \Lambda V \otimes V$ satisfies

$$\partial_R(f_1 \cdots f_n) = \sum_{k=1}^{n} (-1)^{n-k}(f_1 \cdots \slashed{f}_k \cdots f_n) \otimes f_k, \qquad f_1, \ldots, f_n \in V.$$

For example, for $v_1, v_2, v_3 \in V$ we have

$$\partial_R(v_1 v_2 v_3) = v_2 v_3 \otimes v_1 - v_1 v_3 \otimes v_2 + v_1 v_2 \otimes v_3,$$

and

$$\partial_R^2(v_1 v_2 v_3) = v_2 \otimes v_3 \otimes v_1 - v_3 \otimes v_2 \otimes v_1 - v_1 \otimes v_3 \otimes v_2 + v_3 \otimes v_1 \otimes v_2 + v_1 \otimes v_2 \otimes v_3 - v_2 \otimes v_1 \otimes v_3.$$

We can also define also a left derivative $\partial_L = (\Pi_V \otimes \mathbb{I}) \circ \Delta \colon \Lambda V \to V \otimes \Lambda V$ with similar properties and

$$\partial_L(f_1 \cdots f_n) = \sum_{k=1}^n (-1)^{k-1} f_k \otimes (f_1 \cdots \not{f_k} \cdots f_n), \qquad f_1, ..., f_n \in V.$$

We will consider $\Lambda V \otimes V$ as a $\Lambda V$-bimodule and in particular we define the bilinear form $\langle \cdot, \cdot \rangle \colon (\Lambda V \otimes V) \otimes V \to \Lambda V$ by

$$\langle f \otimes v, w \rangle = f \langle v, w \rangle, \qquad f \in \Lambda V, v, w \in V. \tag{5}$$

An analogous definition holds for $\langle \cdot, \cdot \rangle \colon V \otimes (V \otimes \Lambda V) \otimes V \to \Lambda V$. Note that there is no ambiguity on whether $\langle \cdot, \cdot \rangle$ denotes one of this two bilinear forms or whether it denotes the scalar product of $V$.

**Remark 9.** The definition of derivative $\partial_R^k$ and the bilinear form (5) do not depend on the fact that $V$ is finite dimensional. For this reason we can define $\partial_R^k$ and the bilinear form (5) also when $V$ is an infinite dimensional pre-Hilbert space.

**Lemma 10.** *Consider $G \in \Lambda V$ and let $X, Y \in \mathcal{G}(V)$ be two compatible Grassmann random variables such that $\|Y\|_{\mathcal{G}(V)} \leqslant 1$, then*

$$G(X+Y) = G(X) + \sum_{k=1}^n \frac{1}{k!} m_{\mathcal{A}}^{k+1}[(X \otimes Y \otimes \cdots \otimes Y)(\partial_R^k G)] + O(\|Y\|_{\mathcal{G}(V)}^{n+1}) \tag{6}$$

**Proof.** We have

$$G(X+Y) = m_{\mathcal{A}}^2[(X \otimes Y)(\Delta G)] = \sum_{k \geqslant 0} m_{\mathcal{A}}^2[(X \otimes Y)(\mathbb{I} \otimes \Pi_{\Lambda^k V})(\Delta G)],$$

where the sum is finite since $F \in \Lambda V$ is a finite polynomial. Equation (6) easily follows as soon as we prove that

$$m_{\mathcal{A}}^{k+1}[(X \otimes Y \otimes \cdots \otimes Y)(\partial_R^k G)] = (k!) m_{\mathcal{A}}^2[(X \otimes Y)(\mathbb{I} \otimes \Pi_{\Lambda^k V})(\Delta G)].$$

We have

$$\begin{aligned}
\partial_R^k G &= (\partial_R \otimes \mathbb{I}^{\otimes k}) \cdots (\partial_R \otimes \mathbb{I}) \partial_R G \\
&= ((\mathbb{I} \otimes \Pi_V) \Delta \otimes \mathbb{I}^{\otimes k}) \cdots ((\mathbb{I} \otimes \Pi_V) \Delta \otimes \mathbb{I})(\mathbb{I} \otimes \Pi_V) \Delta G \\
&= ((\mathbb{I} \otimes \Pi_V) \otimes \mathbb{I}^{\otimes k})(\Delta \otimes \mathbb{I}^{\otimes k}) \cdots ((\mathbb{I} \otimes \Pi_V) \otimes \mathbb{I})(\Delta \otimes \mathbb{I})(\mathbb{I} \otimes \Pi_V) \Delta G \\
&= ((\mathbb{I} \otimes \Pi_V) \otimes \mathbb{I}^{\otimes k}) \cdots ((\mathbb{I} \otimes \Pi_V) \otimes \mathbb{I})(\mathbb{I} \otimes \Pi_V)(\Delta \otimes \mathbb{I}^{\otimes k}) \cdots (\Delta \otimes \mathbb{I}) \Delta G \\
&= (\mathbb{I} \otimes \Pi_V^{\otimes k+1})(\Delta \otimes \mathbb{I}^{\otimes k}) \cdots (\Delta \otimes \mathbb{I}) \Delta G.
\end{aligned}$$

Let $N := \{1, ..., n\}$, $\#N = n$. For $I = \{j_1, ..., j_k\} \subset N$, we put $v_I := v_{j_1} \cdots v_{j_k}$ for $v_j \in V$, $j \in N$. Let $\mathrm{sgn}(N \setminus I) \in \{\pm 1\}$ be such that $v_N = \mathrm{sgn}(N \setminus I) v_{N \setminus I} v_I$. By definition of the coproduct $\Delta$ in the exterior algebra $\Lambda V$, we have

$$(\mathbb{I} \otimes \Pi_V^{\otimes k+1})(\Delta \otimes \mathbb{I}^{\otimes k}) \cdots (\Delta \otimes \mathbb{I})(\Delta)(v_1 \cdots v_n) = \\ = \sum_{j_1, ..., j_k \in N} \mathrm{sgn}(N \setminus I) \, v_{N \setminus \{j_1, ..., j_k\}} \otimes v_{j_1} \otimes \cdots \otimes v_{j_k}.$$

On the other hand

$$\begin{aligned}(\mathbb{I}\otimes\Pi_{\Lambda^k V})\Delta v_1\cdots v_n &= (\mathbb{I}\otimes\Pi_{\Lambda^k V})\sum_{I\subset N}\operatorname{sgn}(N\setminus I)\,v_{N\setminus I}\otimes v_I\\ &= \frac{1}{k!}\sum_{j_1,\ldots,j_k}\operatorname{sgn}(N\setminus I)\,v_{N\setminus\{j_1,\ldots,j_k\}}\otimes v_{\{j_1,\ldots,j_k\}}.\end{aligned}$$

By linearity, we therefore have

$$(\mathbb{I}\otimes m_{\Lambda V}^k)(\partial_R^k G) = k!(\mathbb{I}\otimes\Pi_{\Lambda^k V})(\Delta G),$$

where $m_{\Lambda V}^k\colon V^{\otimes k}\subset(\Lambda V)^{\otimes k}\to\Lambda V$ is defined as in equation (2). Since $X,Y$ are algebra homomorphisms they commute with the product, hence we have

$$\begin{aligned}m_{\mathcal{A}}^{k+1}[(X\otimes Y\otimes\cdots\otimes Y)(\partial_R^k G)] &= m_{\mathcal{A}}^2[(X\otimes Y)(\mathbb{I}\otimes m_{\Lambda V}^k)(\partial_R^k G)]\\ &= k!\,m_{\mathcal{A}}^2[(X\otimes Y)(\mathbb{I}\otimes\Pi_{\Lambda^k V})(\Delta G)].\end{aligned}$$

This concludes the proof. $\square$

We want to give now some precise bounds on the norms $\|G(X)\|_{\mathcal{A}}$ and $\|G(X)-G(Y)\|_{\mathcal{A}}$, where $G\in\Lambda V$ and $X,Y\in\mathcal{G}(V)$ are two compatible Grassmann random variables. First we introduce a suitable norm $\|\cdot\|_\pi$ on $\Lambda V$ that is the norm induced by the projective norm on $\bigoplus_{n=0}^{\dim(V)}\bigotimes^n V$. Note that there exists an injection $i_{\Lambda^n V}\colon\Lambda^n V\to\bigotimes^n V$ given by the unique linear extension of the following relation

$$i_{\Lambda^n V}(v_1\wedge\cdots\wedge v_n)=\frac{1}{n!}\sum_{\sigma\in S_n}(-1)^\sigma v_{\sigma(1)}\otimes v_{\sigma(2)}\otimes\cdots\otimes v_{\sigma(n)}.\qquad v_1,\ldots,v_n\in V.$$

We then define the map $i_{\Lambda V}:=\bigoplus_{n=0}^{\dim(V)}i_{\Lambda^n V}$. On $\bigoplus_{n=0}^{\dim(V)}\bigotimes^n V$ we consider the projective norm $\|\cdot\|_\pi$ defined as follows. If $f\in\bigotimes^n V$ then

$$\|f\|_\pi:=\inf\left\{\sum_{k=1}^p\|f_1^k\|_V\cdots\|f_n^k\|_V,\text{ where }f=\sum_{k=1}^p f_1^k\otimes\cdots\otimes f_n^k\right\}.$$

For a general element $g\in\bigoplus_{n=0}^{\dim(V)}\bigotimes^n V$ we put $\|g\|_\pi:=\sum_{k=0}^{\dim(V)}\|\Pi_{\otimes^n V}(g)\|_\pi$, and for any $w\in\Lambda V$,

$$\|w\|_{\Lambda_\pi V}:=\|i_{\Lambda V}(w)\|_\pi. \tag{7}$$

Since $i_{\Lambda V}$ is an injection, $\|\cdot\|_\pi$ defines a norm on $\Lambda V$. If $G\in\Lambda V$ we define

$$\deg(G)=\max\{n\in\mathbb{N}\text{ such that }\|\Pi_{\Lambda^n V}(G)\|_\pi\neq 0\}.$$

Let $W$ a vector space and introduce the symmetrizer $S\colon TW\to TW$ as the operator from the tensor algebra generated by $W$ in itself that is the unique linear extension of

$$S(w_1\otimes\cdots\otimes w_n)=\frac{1}{n!}\sum_{\sigma\in S_n}w_{\sigma(1)}\otimes\cdots\otimes w_{\sigma(n)},\qquad w_1,\ldots,w_n\in W.$$

It is important to note that for $S((v+w)^{\otimes n})$ the binomial formula holds, i.e. we have

$$S((v+w)^{\otimes n})=\sum_{k=0}^n\binom{n}{k}S(v^{\otimes k}\otimes w^{\otimes(n-k)}). \tag{8}$$

**Lemma 11.** *If $G_n \in \Lambda^n V$ and $X, Y \in \mathcal{G}(V)$ are two compatible Grassmann random variables we have that*

$$m_{\mathcal{A}}^{k+1}[(X \otimes Y \otimes \cdots \otimes Y)(\partial_R^k G_n)] = \frac{n!}{(n-k)!} m_{\mathcal{A}}^n(S(X^{\otimes(n-k)} \otimes Y^{\otimes k})(i_{\Lambda V}(G_n))).$$

**Proof.** For any $\lambda \in \mathbb{R}$, by Lemma 10

$$G_n(X + \lambda Y) = G_n(X) + \sum_{k=1}^n \frac{\lambda^k}{k!} m_{\mathcal{A}}^{k+1}[(X \otimes Y \otimes \cdots \otimes Y)(\partial_R^k G_n)]. \tag{9}$$

On the other hand, by the binomial formula (8)

$$\begin{aligned} G_n(X + \lambda Y) &= m_{\mathcal{A}}^n[S((X + \lambda Y)^{\otimes n}) i_{\Lambda V}(G_n)] \\ &= m_{\mathcal{A}}^n \left[ \sum_{k=0}^n \binom{n}{k} \lambda^k S(X^{\otimes(n-k)} \otimes Y^{\otimes k}) i_{\Lambda V}(G_n) \right]. \end{aligned} \tag{10}$$

By comparing the expressions (9) and (10) as polynomials in $\lambda$ we conclude. $\square$

**Theorem 12.** *Let $X, Y \in \mathcal{G}(V)$ be two compatible Grassmann random variables and let $G \in \Lambda V$ then, for any $n \leqslant \deg(G) - 1$, we have*

$$\|G(X)\|_{\mathcal{A}} \leqslant \|G\|_{\Lambda_\pi V}(1 + \|X\|_{\mathcal{G}(V)})^{\deg(G)}, \tag{11}$$

*and*

$$\left\| G(Y) - G(X) - \sum_{k=1}^n \frac{1}{k!} m_{\mathcal{A}}^{k+1}[(X \otimes (Y-X)^{\otimes k})(\partial_R^k G)] \right\|_{\mathcal{A}} \leqslant$$
$$\leqslant C_{n, \deg(G)} (1 + \max(\|X\|_{\mathcal{G}(V)}, \|Y\|_{\mathcal{G}(V)}))^{\deg(G) - n - 1} \|G\|_{\Lambda_\pi V} \|Y - X\|_{\mathcal{G}(V)}^{n+1} \tag{12}$$

*where $C_{n, \deg(G)} > 0$ is a suitable constant depending only on $n$ and $\deg(G)$. In the case $n = 0$ we can choose $C_{0, \deg(G)} = \deg(G)$.*

**Proof.** First of all we note that if $(G_1, \ldots, G_n) \in \mathcal{L}(V, \mathcal{A})$ and $\tilde{G} \in \bigotimes^n V$ we have

$$\|m_{\mathcal{A}}^n((G_1 \otimes \cdots \otimes G_n)(\tilde{G}))\|_{\mathcal{A}} \leqslant \|\tilde{G}\|_\pi \|G_1\|_{\mathcal{L}(V, \mathcal{A})} \cdots \|G_n\|_{\mathcal{L}(V, \mathcal{A})}. \tag{13}$$

Furthermore, if $X \in \mathcal{G}(V)$ and $G \in \Lambda V$ we have

$$G(X) = \sum_{n=0}^{\deg(G)} m_{\mathcal{A}}^n(X^{\otimes n}(i_{\Lambda^n V}(\Pi_{\Lambda^n V}(G)))).$$

Using then the definition of $\|X\|_{\mathcal{G}(V)}$ and of $\|G\|_{\Lambda_\pi V}$ we get

$$\|G(X)\|_{\mathcal{A}} \leqslant \sum_{n=0}^{\deg(G)} \|X\|_{\mathcal{G}(V)}^n \|i_{\Lambda^n V}(\Pi_{\Lambda^n V}(G))\|_\pi \leqslant (1 + \|X\|_{\mathcal{G}(V)})^n \|G\|_\pi.$$

In general, writing $G_h = \Pi_{\Lambda^h V}(G)$, for $n + 1 \leqslant h \leqslant \deg(G)$, we have

$$m_{\mathcal{A}}^{k+1}[(X \otimes Y^{\otimes k})(\partial_R^k G_n)] = \frac{n!}{(n-k)!} m_{\mathcal{A}}^n(S(X^{\otimes(n-k)} \otimes Y^{\otimes k})(i_{\Lambda V}(G_n))).$$

$$G_h(Y) - G_h(X) - \sum_{k=1}^{n} \frac{1}{k!} m_{\mathcal{A}}^{k+1}[(X \otimes (Y-X)^{\otimes k})(\partial_R^k G_h)] =$$

$$= \sum_{k=n+1}^{h} \frac{1}{k!} m_{\mathcal{A}}^{k+1}[(X \otimes (Y-X)^{\otimes k})(\partial_R^k G_h)]$$

$$= m_{\mathcal{A}}^h \left[ \left( \sum_{k=n+1}^{h} \binom{h}{k} S(X^{\otimes(h-k)} \otimes (Y-X)^{\otimes k}) \right) (i_{\Lambda V}(G_h)) \right]$$

$$= m_{\mathcal{A}}^h \left[ \left( \sum_{k=n+1}^{h} \sum_{\ell=0}^{k-n-1} \binom{h}{k}\binom{k-n-1}{\ell} \cdot \right. \right.$$

$$\left. \left. \cdot S(X^{\otimes(h-k)} \otimes (-X)^{\otimes \ell} \otimes Y^{\otimes(k-n-1-\ell)} \otimes (X-Y)^{\otimes(n+1)}) \right) (i_{\Lambda V}(G_h)) \right] \quad (14)$$

$$= m_{\mathcal{A}}^h \left[ \left( \sum_{k'=0}^{h-n-1} \sum_{\ell=0}^{k'} \binom{h}{k'+n+1}\binom{k'}{\ell} \cdot \right. \right.$$

$$\left. \left. \cdot S(X^{\otimes(h-k'-n-1)} \otimes (-X)^{\otimes \ell} \otimes Y^{\otimes(k'-\ell)} \otimes (X-Y)^{\otimes(n+1)}) \right) (i_{\Lambda V}(G_h)) \right]$$

$$= m_{\mathcal{A}}^h \left[ \left( \sum_{p=0}^{h-n-1} S(X^{\otimes(h-n-1-p)} \otimes Y^{\otimes p} \otimes (X-Y)^{\otimes(n+1)}) \cdot \right. \right.$$

$$\left. \left. \cdot \left( \sum_{\ell=0}^{h-n-1-p} (-1)^{\ell}\binom{h}{p+\ell+n+1}\binom{p+\ell}{\ell} \right) \right) (i_{\Lambda V}(G_h)) \right],$$

where the first equality follows from (9), the second from Lemma 11; the third by the linearity of $m_{\mathcal{A}}^h$, the binomial theorem for the tensor power $(Y-X)^{\otimes(k-n-1)}$, and the fact that inside the symmetrizer $S$ the tensor product commutes.

Using inequality (13) in relation (14) we get

$$\left\| G_h(Y) - G_h(X) - \sum_{k=1}^{n} \frac{1}{k!} m_{\mathcal{A}}^{k+1}[(X \otimes (Y-X)^{\otimes k})(\partial_R^k G_h)] \right\|_{\mathcal{A}} \leqslant$$
$$\leqslant \|i_{\Lambda V}(G_h)\|_{\Lambda_\pi V} C_{n,h} (\max(\|X\|_{\mathcal{G}(V)}, \|Y\|_{\mathcal{G}(V)}))^{h-n-1} \|Y-X\|_{\mathcal{G}(V)}^{n+1},$$

where we can choose

$$C_{n,h} = \sum_{p=0}^{h-n-1} \left| \sum_{\ell=0}^{h-n-1-p} (-1)^{\ell}\binom{h}{p+\ell+n+1}\binom{p+\ell}{\ell} \right|.$$

In the case $n=0$ we can get a better constant. Indeed we have

$$\|G_h(Y) - G_h(X)\|_{\mathcal{A}} = \left\| \sum_{k=0}^{h-1} S(X^{\otimes k} \otimes Y^{\otimes h-1-k} \otimes (Y-X)) \right\|_{\mathcal{G}(V)} \|G_h\|_{\Lambda_\pi V}$$

$$\leqslant h(\max(\|X\|_{\mathcal{G}(V)}, \|Y\|_{\mathcal{G}(V)}))^{h-1} \|X-Y\|. \qquad \square$$

We now introduce the notion of function "depending on the space variable" $v \in V$ and a suitable norm on the space of such functions. If $F \in \mathcal{L}(V, \Lambda V)$ we can define the composition

$$F(X)(v) := X(F(v)), \qquad X \in \mathcal{G}(V), v \in V.$$

as a linear (and continuous) map $F(X): V \to \mathcal{A}$. We define $\|F\|_{\Lambda_\pi V}$ as

$$\|F\|_{\Lambda_\pi V} := \|F\|_{\mathcal{L}(V, \Lambda_\pi V)} = \sup_{v \in V, |v|=1} \|F(v)\|_{\Lambda_\pi V}. \tag{15}$$

where in the r.h.s. we use the norm $\|\cdot\|_{\Lambda_\pi V}$ defined in (7).

**Remark 13.** If $F: V \to \Lambda_{\mathrm{odd}} V$ we have that for any $v_1, v_2 \in V$

$$F(X)(v_1) F(X)(v_2) = -F(X)(v_2) F(X)(v_1).$$

This means that $F(X)$ can be extend to an homomorphism $\mathrm{Hom}(\Lambda V, \mathcal{A})$ (we shall denote this extension by $F(X)$ too). Furthermore, since

$$F(X)(v_1) X(v_2) = -X(v_2) F(X)(v_1),$$

we have that $F(X)$ and $X$ are compatible Grassmann random variables.

We can define also $\partial_R^n F: V \to \Lambda V \otimes (\bigotimes^n V)$ as $\partial_R^n F(v) = \partial_R^n(F(v))$, and $\deg(F)$ in the obvious way.

**Theorem 14.** *Suppose that $F \in \mathcal{L}(V, \Lambda V)$ and let $X, Y \in \mathcal{G}(V)$ be compatible Grassmann random variables, then we have that*

$$\|F(X)\|_{\mathcal{G}(V)} \leqslant \|F\|_{\Lambda_\pi V} (1 + \|X\|_{\mathcal{G}(V)})^{\deg(F)}$$

$$\left\| F(Y) - F(X) - \sum_{k=1}^n \frac{1}{k!} m_{\mathcal{A}}^{k+1} [(X \otimes (Y-X)^{\otimes k})(\partial_R^k F)] \right\|_{\mathcal{G}(V)} \leqslant$$
$$\leqslant C_{n, \deg(F)} (1 + \max(\|X\|_{\mathcal{G}(V)}, \|Y\|_{\mathcal{G}(V)}))^{\deg(F)-n-1} \|F\|_{\Lambda_\pi V} \|Y - X\|_{\mathcal{G}(V)}^{n+1}$$

**Proof.** The proof is a simple application of Theorem 12 to $F(X)(v)$ and $F(Y)(v)$ for any fixed $v \in V$ and of the definition of the norm (15). $\square$

### 2.3 Grassmann Gaussian variables

Let now $V$ be a real pre-Hilbert space with scalar product $\langle \cdot, \cdot \rangle$ and with an antisymmetric bounded operator $C: V \to V$. By Remark 6 and Remark 9, we can extend the definitions of $\|\cdot\|_{\mathcal{G}(V)}$, $\partial_R^k$ and the bilinear form (5) from the finite dimensional case to the generic (infinite) pre-Hilbert space $V$.

**Definition 15.** *A (V-)Grassmann (centered) Gaussian variable with correlation $C$ is a random variable $X \in \mathcal{G}(V)$ such that*

$$\omega(X(G) X(f)) = \omega(X(\langle \partial_R G, Cf \rangle)), \qquad G \in \Lambda V, f \in V. \tag{16}$$

*We also require that $\|X\|_{\mathcal{G}(V)} < \infty$, i.e that the map $X: V \to \mathcal{A}$ must be continuous with respect the topology induced on $V$ by the pre-Hilbert product structure and the (norm) topology of $\mathcal{A}$.*

If we define $\partial_R G(X) = X(\partial_R G)$, with the understanding that $X(g \otimes w) = X(g) \otimes w$, $g \in \Lambda V$ and $w \in V$, then we have the more familiar expression (similar to the bosonic counterpart)

$$\omega(G(X) X(f)) = \omega(\langle \partial_R G(X), Cf \rangle).$$

Note that the integration by parts formula determines all the moments of the Gaussians. In particular (16) implies that $\omega(X(f_1)\cdots X(f_n)) = 0$ if $n$ is odd and if $n = 2k$ is even

$$\omega(X(f_1)\cdots X(f_{2k})) = \sum_{\sigma} (-1)^{\sigma} \prod_{i=1}^{k} \langle f_{\sigma(2i-1)}, Cf_{\sigma(2i)}\rangle_V \qquad (17)$$

where the sum runs over all the pairings

$$\sigma = \{(i_1, j_1), ..., (i_k, j_k)\}, \quad i_1 < i_2 < \cdots < i_k, i_{\tilde{r}}, j_{r'} \neq j_r$$

(where $\tilde{r}, r', r = 1, ..., k$ and $r' \neq r$) of $\{1, ..., 2k\}$ and $(-1)^{\sigma}$ is an appropriate sign. The sign $(-1)^{\sigma}$ is obtained in the following way: we can identify the pairing $\sigma = \{(i_1, j_1), ..., (i_k, j_k)\}$ with the ordered list $(i_1, j_1, ..., i_k, j_k)$ (where $i_1 < \cdots < i_k$). Then $(-1)^{\sigma} := (-1)^{P_\sigma}$ denotes the sign of the permutation $P_\sigma \in S_{2k}$ (of $\{1, ..., 2k\}$) such that, for any $r = 1, ..., k$, $P_\sigma(2r-1) = i_r$ and $P_\sigma(2r) = j_r$.

Equation (17) is often called Wick's rule. The right hand side of (17) can be written as a Pfaffian:

$$\omega(X(f_1)\cdots X(f_n)) = \Pf_{1 \leqslant i,j \leqslant n} \langle f_i, Cf_j\rangle_V, \qquad (18)$$

where the Pfaffian is defined, for an antisymmetric $n \times n$ matrix $M$, to be zero if $n = 2k+1$ and if $n = 2k$ as the polynomial in the entries of $M$ which satisfy the relation

$$\left[\Pf_{1 \leqslant i,j \leqslant n} M_{ij}\right]^2 = \det_{1 \leqslant i,j \leqslant n} M_{ij}.$$

The existence of a Gaussian variable, in the above sense, implies the inequality

$$\det_{1 \leqslant i,j \leqslant n} \langle f_i, Cf_j\rangle_V = [\omega(X(f_1)\cdots X(f_n))]^2 \leqslant \|X(f_1)\|^2 \cdots \|X(f_n)\|^2 \leqslant \|X\|^{2n} \|f_1\|^2 \cdots \|f_n\|^2,$$

well known in the mathematical physics literature relative to fermionic expansions e.g. see [118].

**Remark 16.** The averages of Gaussian variables depend only on the quadratic form $(f, g) \mapsto \langle f, Cg\rangle_V$, however analysis on the Grassmann algebra relies on the scalar product $\langle \cdot, \cdot \rangle_V$ also via the requirement that $\|X(f)\|_{\mathcal{A}} \leqslant \|f\|_V$. In particular the realization of the Grassmann algebra as a family of bounded operators is not canonical. This is the reason we need to require the above continuity of the map $X: V \to \mathcal{A}$.

In order to construct Grassmann Gaussian variables, we need "a place to hang the hat on", this place will be a (canonical commuting relations) CAR algebra endowed with its vacuum state. To allow for arbitrary correlations $C$ we can use the approach of OS [131] which consists in doubling the generators of the CAR algebra with respect to the generators of the Grassmann algebra.

**Lemma 17.** *For every antisymmetric and bounded $C: V \to V$ there exists a $(V\text{-})$Grassmann Gaussian variable $X$ with correlation $C$ (on a suitable probability space $(\mathcal{A}, \omega)$).*

**Remark 18.** It is important to note that Lemma (17) is peculiar of the Grassmannian setting and it has no equivalent in the commutative case in the following sense: if $V$ is and pre-Hilbert space and $S$ is a bounded positive, self-adjoint operator there is no map $X^{\text{com}}: V \to \mathcal{A}$ such that, for any $v_1, ..., v_{2k} \in V$,

$$\omega(X^{\text{com}}(v_1)\cdots X^{\text{com}}(v_{2k})) = \sum_{\mathcal{M} \in \{\text{perfect matchings of } \{1, ..., 2k\}\}} \prod_{(i,j) \in \mathcal{M}} \langle Sv_i, v_j\rangle$$

for some $C^*$-algebra $\mathcal{A}$ with a positive state $\omega$. Indeed in the commutative case the Gaussian variables are unbounded and so we must realize them in a algebra $\tilde{\mathcal{A}}$ of unbounded operators (for example we can consider $\tilde{\mathcal{A}} = \bigcap_{1 \leqslant p < +\infty} L^p((\Omega, \mathcal{F}, \mathbb{P}), \mathbb{C})$, for a suitable probability space $(\Omega, \mathcal{F}, \mathbb{P})$ and $\omega(\cdot) = \mathbb{E}[\cdot]$).

**Proof of Lemma 17.** Let us consider $\Lambda V$ itself as a real pre-Hilbert space with respect to the scalar product $\langle \cdot, \cdot \rangle_{\Lambda V}$ on $\Lambda V$ given by

$$\langle v_1 \wedge \cdots \wedge v_n, w_1 \wedge \cdots \wedge w_m \rangle_{\Lambda V} = \delta_{nm} \det_{1 \leqslant j,k \leqslant n} \langle v_j, w_k \rangle,$$

for $v_j, w_k \in V$, $j = 1, ..., n$, $k = 1, ..., m$, where $\langle \cdot, \cdot \rangle$ denotes as usual the scalar product on $V$. Note that $\langle \cdot, \cdot \rangle_{\Lambda V}$ on $\Lambda^n V$ is simply the restriction of the Hilbert scalar product on the tensor product $V^{\otimes n}$. Let $\mathcal{H}$ be the completion of $\Lambda V$ with respect to $\langle \cdot, \cdot \rangle_{\Lambda V}$ and denote by $\Omega$ the element in $\mathcal{H}$ which corresponds to $1 \in \Lambda V$, often referred to as the *vacuum* vector. Denote by $\lambda: \Lambda V \to \mathrm{End}(\Lambda V)$ the left action of $\Lambda V$ on itself given by $\lambda(H) G = H \wedge G$, where $H, G \in \Lambda V$. We show that this action extends to a representation of $\Lambda V$ on $\mathcal{H}$. Indeed $\lambda(v)$, $v \in V$ corresponds to a creation operator. Let us denote by $\lambda(v)^{\mathrm{T}}$ the adjoint of the operator $\lambda(v)$ with respect to the scalar product $\langle \cdot, \cdot \rangle_{\Lambda V}$. A standard computation shows that for $v \in V$ we have

$$\lambda(v)^{\mathrm{T}} w_1 \wedge \cdots \wedge w_n = \sum_{\ell=1}^n (-1)^{\ell-1} \langle v, w_\ell \rangle w_1 \wedge \cdots \wedge \cancel{w_\ell} \wedge \cdots \wedge w_n.$$

Namely $\lambda^T(v) x$ can be expressed in terms of the left derivative as $\lambda^T(v) x = \langle v, \partial_L x \rangle$ (recall the definition given in equation (5)) for all $x \in \Lambda V$. In particular we have

$$\{\lambda(v), \lambda(w)\} = \{\lambda(v)^{\mathrm{T}}, \lambda(w)^{\mathrm{T}}\} = 0, \qquad \{\lambda(v)^{\mathrm{T}}, \lambda(w)\} = \langle v, w \rangle, \qquad v, w \in V, \qquad (19)$$

where as usual $\{\cdot, \cdot\}$ denotes the anticommutator. As a consequence we obtain, for any $x \in \Lambda V$

$$\langle \lambda(v) x, \lambda(v) x \rangle_{\Lambda V} + \langle \lambda(v)^{\mathrm{T}} x, \lambda(v)^{\mathrm{T}} x \rangle_{\Lambda V} = \langle x, \{\lambda(v)^{\mathrm{T}}, \lambda(v)\} x \rangle_{\Lambda V} = \langle v, v \rangle \langle x, x \rangle_{\Lambda V}.$$

In particular $\lambda(v), \lambda(v)^{\mathrm{T}}, v \in V$, extend to bounded operators on $\mathcal{H}$ (still denoted by the same symbol). We now define $X: \Lambda V \to \mathcal{L}(\mathcal{H})$ to be the algebra homomorphism given on $V$ by

$$X(v) := \lambda(Cv) + \lambda(v)^{\mathrm{T}}, \qquad v \in V.$$

Then $X$ is extended in a natural way to $\Lambda V$ (still denoted by the same symbol). Note that $X(x)$, $x \in \Lambda V$, is indeed a bounded operator since it is the sum of products of bounded operators. Finally we define $\omega: \mathcal{L}(\mathcal{H}) \to \mathbb{R}$ to be the state on $\mathcal{A} = \mathcal{L}(\mathcal{H})$ defined by $\Omega$, that is

$$\omega(T) := \langle \Omega, T \Omega \rangle_{\Lambda V}, \qquad T \in \mathcal{L}(\mathcal{H}).$$

We claim that $X$ is a Gaussian Grassmann random variable on $(\mathcal{L}(\mathcal{H}), \omega)$. We first show that it is Grassmann:

$$\begin{aligned} \{X(v), X(w)\} &= \{\lambda(Cv), \lambda(w)^{\mathrm{T}}\} + \{\lambda(v)^{\mathrm{T}}, \lambda(Cw)\} \\ &= \langle Cv, w \rangle + \langle v, Cw \rangle = \langle v, C^{\mathrm{T}} w \rangle + \langle v, Cw \rangle = 0, \end{aligned}$$

where $v, w \in \Lambda$, because by assumption $C^{\mathrm{T}} = -C$. Now we have also that

$$\begin{aligned} \omega(X(v) X(w)) &= \langle \Omega, \lambda(v)^{\mathrm{T}} \lambda(Cw) \Omega \rangle_{\Lambda V} \\ &= \langle \lambda(v) \Omega, \lambda(Cw) \Omega \rangle_{\Lambda V} = \langle v, Cw \rangle_{\Lambda V}, \end{aligned}$$

where the first equality follows from the fact that $\lambda(v)^T \Omega = 0$. More generally using that $X(v)\Omega = \lambda(Cv)\Omega$ and the commutation relation (19) we have that, for all $f_1, ..., f_{2k} \in V$

$$\omega(X(f_1)\cdots X(f_{2k})) = \omega(X(f_1)\cdots X(f_{2k-1})\lambda(Cf_{2k}))$$

$$= -\omega(X(f_1)\cdots \lambda(Cf_{2k})X(f_{2k-1})) + \omega(X(f_1)\cdots X(f_{2k-2}))\langle f_{2k-1}, Cf_{2k}\rangle$$

$$= \cdots = \sum_{\ell=1}^{2k-1} (-1)^{2k-1-\ell}\omega\big(X(f_1)\cdots X(\slashed{f_\ell})\cdots X(f_{2k-1})\big)\langle f_\ell, Cf_{2k}\rangle$$

Therefore

$$\omega(X(G)X(f)) = \omega(X(\langle \partial_R G, Cf\rangle)), \qquad G \in \Lambda V, f \in V.$$

Note that, similarly

$$\omega(X(f)X(G)) = \omega(X(\langle f, C\partial_L G\rangle)). \tag{20}$$

$\square$

*Complex Gaussians.* Later on we will need also complex Grassmann Gaussian variables, i.e. Gaussian variables taking values in $\Lambda V$ with a complex pre-Hilbert space $V$ whose Hermitian scalar product we denote by $(\cdot,\cdot)_V$ and we assume anti-linear in the left variable. As in the commutative setting, their definition poses no particular problems, however the interplay of the algebraic and analytic structure is here reflected on the fact that we need to fix a real structure $\varkappa$ on $V$.

**Definition 19.** *A real structure $\varkappa: V \to V$ compatible with $(\cdot,\cdot)_V$ is an anti-unitary involution, i.e a map such that, for all $\alpha, \beta \in \mathbb{C}$, $v, w \in V$,*

$$\varkappa(\alpha v + \beta w) = \bar{\alpha}\varkappa v + \bar{\beta}\varkappa w, \quad (\varkappa v, \varkappa w)_V = \overline{(v,w)}_V, \quad \varkappa\varkappa = 1.$$

Since we will deal with Hilbert spaces we will always assume a real structure to be compatible with the scalar product and we will not explicitly say that it is compatible. For a real structure $\varkappa$ the following properties hold:

- Let $\text{Re}_\varkappa v := (v + \varkappa v)/2$, $\text{Im}_\varkappa v := (v - \varkappa v)/(2i)$, then $v = \text{Re}_\varkappa v + i\text{Im}_\varkappa v$, $\varkappa\text{Re}_\varkappa = \text{Re}_\varkappa$, $\varkappa\text{Im}_\varkappa = -\text{Im}_\varkappa$, $\text{Re}_\varkappa = \text{Im}_\varkappa i$, $\text{Im}_\varkappa = -\text{Re}_\varkappa i$, so

$$V = \text{Re}_\varkappa V + i\text{Im}_\varkappa V = \text{Re}_\varkappa V + i\text{Re}_\varkappa iV,$$

which is a direct sum in the sense of vector spaces but it is not orthogonal with respect to $(\cdot,\cdot)_V$.

- The form

$$\langle\!\langle v, w \rangle\!\rangle_\varkappa := (\varkappa v, w)_V$$

is bilinear, symmetric and non-degenerate.

- For any linear operator $A: V \to V$ we define its $\varkappa$-transpose

$$A^\varkappa := \varkappa A^* \varkappa$$

which satisfies $\langle\!\langle v, Aw \rangle\!\rangle_\varkappa = \langle\!\langle A^\varkappa v, w \rangle\!\rangle_\varkappa$. Here we denote $A^*$ the usual Hermitian adjoint of the linear operator $A: V \to V$.

- If $\kappa$ is another real structure on $V$, then letting $U := \kappa\varkappa$ we have that $U$ is a linear operator which is invertible and $(Uv, Uw)_V = (v,w)_V$, that is $U$ is unitary.

Note that if $V$ is a space of complex functions, then the map $c: V \to V$ given by taking the complex conjugate is a real structure.

**Definition 20.** *Let V a complex pre-Hilbert space, $\varkappa$ a real structure over it and $C\colon V \to V$ a $\varkappa$-antisymmetric (i.e. $C^\varkappa = -C$) bounded linear operator. A $(V,\varkappa)$-Grassmann (centered) Gaussian variable with correlation $C$ is a random variable $X \in \mathcal{G}(V)$ such that*

$$\omega(X(G)X(h)) = \omega(X(\langle\!\langle \partial_R G, Ch \rangle\!\rangle_\varkappa)), \qquad G \in \Lambda V, h \in V,$$

*where $\Lambda V$ is the Grassmann algebra over $\mathbb{C}$ generated by V. We require that $\|X\|_{\mathcal{G}(V)} < \infty$, i.e that the map $X\colon V \to \mathcal{A}$ is continuous.*

**Lemma 21.** *For every $\varkappa$-antisymmetric and bounded $C\colon V \to V$ there exists a $(V,\varkappa)$-Grassmann centered Gaussian variable $X \in \mathcal{G}(V)$ with correlation $C$ (on a suitable probability space $(\mathcal{A},\omega)$).*

**Proof.** The construction of a complex Grassmann Gaussian proceeds as in Lemma 17 by considering the (complex) antisymmetric Fock space $\mathcal{H}$ associated to $V$ with vacuum vector $\Omega \in \mathcal{H}$ and creation operators $(a(v))_{v \in V}$ which are linear and satisfying the canonical anticommutation relations (CAR) $\{a(w)^*, a(v)\} = (w,v)_V$ for $v,w \in V$. Let

$$X(v) = a(Cv)^* + a(\varkappa v), \qquad v \in V,$$

and consider the state $\omega(A) := \langle \Omega, A\Omega \rangle_\mathcal{H}$ for any $A \in \mathcal{A} = \mathcal{L}(\mathcal{H})$. The verification that the bounded operators $(X(v))_{v \in V}$ forms a Grassmann algebra and that

$$\omega(X(v)X(w)) = \langle\!\langle v, Cw \rangle\!\rangle_\varkappa, \qquad v,w \in V,$$

as required, is left to the reader. $\square$

*White noise.* A relevant example of a random variable taking values in an infinite dimensional Grassmann algebra is $d$-dimensional (Gaussian) white noise, defined as follows.

**Definition 22.** *A V-valued d-dimensional (Gaussian) white noise with correlation $C\colon V \to V$ is the (centered) Grassmann Gaussian variable $\Xi \in \mathcal{G}(L^2(\mathbb{R}^d) \otimes V)$ with correlation $\tilde{C}$ given by*

$$(\tilde{C}f)(x) = Cf(x)$$

*for all $x \in \mathbb{R}^d$ and $f \in L^2(\mathbb{R}^d) \otimes V \approx L^2(\mathbb{R}^d; V)$.*

Take now a one dimensional white noise $\Xi$ with values in $V$ and let $B_t(v) = \Xi(\mathbb{I}_{[0,t]} \otimes v)$ for $v \in V$. For fixed $t \geqslant 0$, $B_t$ extends as homomorphism to all $\Lambda V$ and therefore $B_t \in \mathcal{G}(V)$. Note also that

$$\|B_t(v) - B_s(v)\|_\mathcal{A} \leqslant \|\Xi\|_{\mathcal{G}(L^2(\mathbb{R}^d) \otimes V)} |t-s|^{1/2} \|v\|_V, \qquad t,s \geqslant 0, v \in V, \tag{21}$$

so $B(v) \in C(\mathbb{R}_+, \mathcal{A})$ (here $C(\mathbb{R}_+, \mathcal{A})$ is the space of continuous maps from $\mathbb{R}_+$ to $\mathcal{A}$). Note that $B \in \text{Hom}(\Lambda V, C(\mathbb{R}_+, \mathcal{A}))$: this in particular implies that $B_{t_1}, ..., B_{t_n}$ is a compatible family and that $B = (B_t)_{t \in \mathbb{R}_+}$ is a Gaussian process with continuous trajectories. We have that $B_0(v) = 0$,

$$\omega(B_t(v)) = 0, \qquad \omega(B_t(v)B_s(w)) = \langle v, Cw \rangle (t \wedge s), \qquad t,s \geqslant 0, \quad v,w \in V,$$

where $C$ is the correlation of the Grassmann white noise $\Xi$. Increments of $B$ are independent and higher order moments can be computed via Wick's rule (17). Note in particular that

$$\sup_{0 \leqslant t \leqslant T} \|B_t\|_{\mathcal{G}(V)} < \infty, \tag{22}$$

for all $T>0$. Let us remark here that properties (21) and (22) are very different from the pathwise properties of the commutative (bosonic) Brownian motion where the random variable realizing Brownian motion is unbounded both with respect the probability space and in time (see Remark 18 for other observations on this topic).

# 3 Finite dimensional SDEs

We want to study simple SDEs taking values in $\Lambda V$ as far as it is needed for the purpose of stochastic quantization, that is with additive white noise. We refrain to undertake here a general study of Grassmann SDEs, in particular no stochastic calculus will be needed below. It seems possible to devise such a calculus but we leave it for a future work. More precisely in this section we want to study SDEs driven by an additive Brownian motion $B_t$ when the linear space $V$ is finite dimensional. For any interval $I \subseteq \mathbb{R}$ we denote by $C^0(I; \mathcal{G}(V))$ the space of continuous function from $I$ to $\mathcal{G}(V)$.

**Definition 23.** *Let $F \in \mathcal{L}(V, \Lambda V)$ and assume that $F(v)$ is odd for all $v \in V$, let $\Psi_0 \in \mathcal{G}(V)$ be a random variable compatible and independent of the Brownian motion $B \in C^0(\mathbb{R}_+; \mathcal{G}(V))$. For $T > 0$, $\Psi \in C^0([0,T]; \mathcal{G}(V))$ is a solution in $[0,T]$ to the (additive) SDE driven by $B$ with drift $F$ and initial condition $\Psi_0$, if,*

$$\Psi_t(v) = \Psi_0(v) + \int_0^t \Psi_s(F(v))\mathrm{d}s + B_t(v), \qquad t \in [0,T], v \in V, \tag{23}$$

*where the integral with respect to the variable $s$ is understood in Bochner's sense with values in $\mathcal{A}$.*

**Remark 24.** Note that we do not require that $\Psi_t, \Psi_s$ are compatible for $t \neq s$, in particular

$$C^0([0,T]; \mathcal{G}(V)) \neq \mathrm{Hom}(\Lambda V; C([0,T], \mathcal{A})),$$

nor we require any compatibility of $\Psi$ and $B$. It will turn out that such compatibility holds in fact for the unique solution of (23), but it is not necessary to put such restriction to formulate the notion of solution.

## 3.1 The Grassmann Ornstein–Uhlenbeck motion

In this section we introduce the Grassmann analog of the Ornstein–Uhlenbeck process, that is when $F(v) = Av$ with $A: V \to V$ is a linear operator on $V$. In this case we can write down an explicit formula for the solution to equation (23) extended to $t \geqslant 0$.

**Proposition 25.** *Assume that $A: V \to V$ and $F(v) = Av$. The unique solution to equation (23) is given by*

$$\Psi_t(v) = \Psi_0(e^{At}v) + \Xi(\mathbb{I}_{[0,t]}(\cdot)e^{A(t-\cdot)}v), \qquad t \in \mathbb{R}_+, \tag{24}$$

*where $\Xi$ is the Grassmann Gaussian noise related with the Brownian motion $B$ (see Definition 22 and the discussion that follows it).*

**Proof.** Let $h(t,v) \in L^2(\mathbb{R}; V)$ be given by

$$h(t,v)(s) := \mathbb{I}_{[0,t]}(s)e^{A(t-s)}v, \qquad s \in \mathbb{R}.$$

Then (24) reads $\Psi_t(v) = \Psi_0(e^{At}v) + \Xi(h(t,v))$. An explicit computation using that

$$h(t,v)(s) - \mathbb{I}_{[0,t]}(s)\, v = \int_0^t A\, h(r,v)(s) \mathrm{d}r = \int_0^t h(r, A\, v)(s) \mathrm{d}r, \qquad s \in \mathbb{R},$$

in $L^2(\mathbb{R}; V)$ gives

$$\Xi(h(t,v)) = B_t(v) + \int_0^t \Xi(h(r, A\, v)) \mathrm{d}r.$$

Moreover, by the definition of exponential of a matrix, $\partial_t \Psi_0(e^{At}v) = \Psi_0(e^{At}A\, v)$ so

$$\Psi_t(v) = \Psi_0(v) + \int_0^t \Psi_0(e^{Ar}A\, v) \mathrm{d}r + B_t(v) + \int_0^t \Xi(h(r, A\, v)) \mathrm{d}r$$
$$= \Psi_0(v) + B_t(v) + \int_0^t \Psi_r(A\, v) \mathrm{d}r$$

as required by equation (23) for our assumption on $F$. $\square$

**Remark 26.** *In particular this shows that if $(\Psi_0, \Xi)$ is a Grassmann Gaussian process then also the solution $\Psi$ to the SDE (23), is a Grassmann Gaussian process compatible with $(\Psi_0, \Xi)$. In order to prove this fact it is sufficient to prove that any product of the form*

$$\omega(\Psi_{t_1}(v_1) \cdots \Psi_{t_{2k}}(v_{2k})) \tag{25}$$

*can be computed using the products (weighted with a suitable sign due to the anti-commutation of $\Psi$) of the covariance $\omega(\Psi_{t_i}(v_i)\Psi_{t_j}(v_j))$ for any $i, j = 1, ..., 2k$. This property of the expectation (25) follows from the fact that*

$$\omega((\Psi_{t_1}(v_1) - \Psi_0) \cdots (\Psi_{t_{2k}}(v_n) - \Psi_0)) = \omega(\Xi(h(t_1, v_1)) \cdots \Xi(h(t_{2k}, v_{2k}))) =$$
$$= \sum_\sigma (-1)^\sigma \prod_{i=1}^k \int_0^{\min(t_i, t_j)} \langle e^{A(t_i - s)} v_{\sigma(2i-1)}, C e^{A(t_j - s)} v_{\sigma(2i)} \rangle_V \mathrm{d}s$$

*where $C$ is the covariance of $\Xi$. Using the fact that $\Psi_0$ is a Gaussian random variable independent of $\Xi$, and so $\Psi_0(e^{A(t-s)}v)$ is also a Gaussian random variable independent of $\Xi$, we obtain the Gaussian behavior of $\Psi_t$. The compatibility of $\Psi_t$ with $(\Psi_0, \Xi)$ follows from the fact that $\Psi$ is a linear function of $(\Psi_0, \Xi)$.*

Let us now study the family of random variables

$$\Psi_t^{\mathfrak{s}}(v) = \Xi(e^{A(t-\cdot)}v), \qquad t \in \mathbb{R}.$$

In the next proposition, we will show that (under some specified assumptions on $A$) this represents the stationary solution to the SDE (24).

**Proposition 27.** *Assume that all eigenvalues of $A$ have strictly negative real part less or equal than $-\lambda_A$, where $\lambda_A > 0$. Then, we have $\sup_{t \in \mathbb{R}} \|\Psi_t^{\mathfrak{s}}\|_{\mathcal{G}(V)} < +\infty$, furthermore for any $G \in \Lambda V$ and any $t \in \mathbb{R}$ we get*

$$\omega(G(\Psi_t^{\mathfrak{s}})) = \omega(G(X^A)), \tag{26}$$

*where $X^A \in \mathcal{G}(V)$ is a Gaussian random variable with covariance $C_A$ given by*

$$C_A := \int_0^\infty e^{A^{\mathrm{T}} s} C e^{As} \mathrm{d}s. \tag{27}$$

*where $A^{\mathrm{T}}$ denotes the transpose matrix.*

**Remark 28.** Expression (26) in particular shows the independence on $t$ of the expectation $\omega(G(\Psi_t^{\mathfrak{s}}))$, expressing the stationarity of $\Psi_t^{\mathfrak{s}}$.

**Proof of Proposition 27.** We have easily by the definition of $\Psi_t^{\mathfrak{s}}$:

$$\|\Psi_t^{\mathfrak{s}}\|_{\mathcal{G}(V)} \lesssim \|\mathbb{I}_{(-\infty,t]}(\cdot)e^{A(t-\cdot)}\|_{L^2(\mathbb{R})} \lesssim \frac{1}{\lambda_A},$$

where $\lesssim$ stands for inequality modulo some appropriate positive constant. The random variable $\Psi_t^{\mathfrak{s}}$ is Gaussian so it is completely characterized (in term of moments) by its covariance. Note that the covariance of $\Psi_t^{\mathfrak{s}}$ can be computed as follows

$$\omega(\Psi_t^{\mathfrak{s}}(v)\Psi_t^{\mathfrak{s}}(w)) = \omega(\Xi(\mathbb{I}_{(-\infty,t]}(\cdot)e^{A(t-\cdot)}v)\Xi(\mathbb{I}_{(-\infty,t]}(\cdot)e^{A(t-\cdot)}w)) =$$

$$= \langle \mathbb{I}_{(-\infty,t]}(\cdot)e^{A(t-\cdot)}v, C\mathbb{I}_{(-\infty,t]}(\cdot)e^{A(t-\cdot)}w\rangle = \int_{-\infty}^0 \langle e^{-A^{\mathrm{T}}\tau}v, Ce^{-A\tau}w\rangle \mathrm{d}\tau = \langle v, C_A w\rangle$$

with $C_A$ given by (27). The appearance of the transposition is due to the fact that $V$ is a *real* pre-Hilbert space. □

Let us point out that

$$A^{\mathrm{T}}C_A + C_A A = -C \qquad (28)$$

by a simple integration by parts applied to (27).

**Remark 29.** Suppose that $A$ commutes with $C$, then we have that $\Psi_t^{\mathfrak{s}}$ is for all $t \in \mathbb{R}$ a Grassmann Gaussian process with correlation $C_A = -(A + A^{\mathrm{T}})^{-1}C$ (the operator $A$ is invertible since all its eigenvalues has strictly negative real part).

## 3.2 Existence and uniqueness for general drift

**Theorem 30.** *For any $\Psi_0 \in \mathcal{G}(V)$ compatible with a given Brownian motion $B$, there exists $T > 0$ such that equation (23) admits a unique solution $\Psi \in C^0([0, T]; \mathcal{G}(V))$. Moreover $\Psi$ is compatible with $\Psi_0 \oplus B$.*

**Proof.** We are going to construct a solution via Picard's iteration. Let $\Psi_t^0 = \Psi_0$ for all $t \geqslant 0$ and define the $n + 1$-th Picard's iteration by

$$\Psi_t^{n+1}(v) = \Psi_0(v) + \int_0^t \Psi_s^n(F(v))\mathrm{d}s + B_t(v), \qquad t \in [0, T],$$

for any $v \in V$. Then $\Psi^0 \in C^0([0, T]; \mathcal{G}(V))$ and since $\Psi_0$ is compatible with $B$ we have that, for all $t \geqslant 0$, $\Psi_t^1$ belong to the Grassmann algebra generated by $\Psi_0$ and $B$ in $\mathcal{A}$ and there $\Psi_t^1(v)$ is an odd element for all $v \in V$. Therefore $\Psi_t^1 \in \mathcal{G}(V)$ for all $t \geqslant 0$ and also $\Psi_t^1$ is compatible with $\Psi_0 \oplus B$. Since this is true for all $t$, by approximation in the operator norm we see that $\Psi^1 \in C^0([0, T]; \mathcal{G}(V))$ and that $\Psi^1$ is compatible with $\Psi_0 \oplus B$. By induction we prove then that $\Psi^n \in C^0([0, T]; \mathcal{G}(V))$ for all $n \geqslant 0$ and that it is compatible with $\Psi_0 \oplus B$.

Now observe that $\Psi_s^n(F(v))$ is a polynomial function of $(\Psi_s^n(v))_{v \in V}$ and therefore by a standard contraction argument in the Banach space $\mathcal{A}$ we obtain that $(\Psi^n)_n$ converges, as $n \to +\infty$, in $C^0([0, T]; \mathcal{A})$ for some positive $T$ which depends only on $F$, on $\|\Psi_0\|_{\mathcal{G}(V)}$ and on $\sup_{t \geqslant 0}\|B_t\|_{\mathcal{G}(V)}$. If we call $\Psi$ the limit we have that $\Psi \in C^0([0, T]; \mathcal{G}(V))$ by the first part of the proof and we deduce easily that $\Psi$ is compatible with $\Psi_0 \oplus B$. □

**Theorem 31.** *A solution to equation (23) exists for all times.*

**Proof.** We will have global existence as soon as we can rule out explosions, that is prove that for all $t \geqslant 0$ we have $\|\Psi_t\|_{\mathcal{G}(V)} < \infty$. For $v \in V$, let $\Theta_t(v) = \Psi_t(v) - B_t(v)$ and extend this map to an homeomorphism of $\Lambda V$ into $\mathcal{A}$. Observe that

$$\Theta_t(v) = \Theta_0(v) + \int_0^t m_{\mathcal{A}}[(\Theta_s \otimes B_s)(\Delta F(v))] \mathrm{d}s,$$

where recall that $m_{\mathcal{A}} \colon \mathcal{A} \otimes \mathcal{A} \to \mathcal{A}$ denotes the multiplication in $\mathcal{A}$. In particular, if we consider a (linear) basis of $\Lambda V$ denoted by $(e_A)_A$ and let $\Theta_t^A := \Theta_t(e_A)$ and $B_t^A := B_t(e_A)$ we have

$$\Theta_t^A = \Theta_0^A + \int_0^t \sum_B h_{A,B,C,D} \Theta_s^B \Theta_0^C B_s^D \mathrm{d}s,$$

for real coefficients $(h_{A,B,C,D})_{A,B,C,D}$. This is a finite-dimensional, linear system of non-autonomous ODEs in $\mathcal{A}$ for $(\Theta_t^A)_A$, and

$$\sum_A \|\Theta_t^A\|_{\mathcal{A}} \leqslant \sum_A \|\Theta_0^A\|_{\mathcal{A}} + |h| \int_0^t \sum_B \|\Theta_s^B\|_{\mathcal{A}} \sum_D \|B_s^D\|_{\mathcal{A}} \sum_B \|\Theta_0^B\|_{\mathcal{A}} \mathrm{d}s,$$

with $|h| := \sup_{A,B,C,D} |h_{A,B,C,D}|$. By Gronwall inequality

$$\sum_A \|\Theta_t^A\|_{\mathcal{A}} \leqslant \sum_A \|\Theta_0^A\|_{\mathcal{A}} \exp\left(|h| \int_0^t \sum_C \|B_s^C\|_{\mathcal{A}} \sum_A \|\Theta_0^A\|_{\mathcal{A}} \mathrm{d}s\right) < \infty,$$

for all $t \geqslant 0$. In particular $\|\Psi_t\|_{\mathcal{G}(V)} \leqslant \|\Theta_t\|_{\mathcal{G}(V)} + \|B_t\|_{\mathcal{G}(V)} < \infty$ for all $t \geqslant 0$ and the solution to equation (23) exists for all times. □

## 3.3 An Itô formula for solutions of SDEs

We want to prove an Itô formula for the SDE (23). For $C \in \mathcal{L}(V, V)$, let $Q_C \in \mathcal{L}(\Lambda V \otimes V \otimes V, \Lambda V)$ be given by

$$Q_C(f \otimes v \otimes w) = \langle v, Cw \rangle f, \qquad f \in \Lambda V, v, w \in V. \tag{29}$$

For $F \in \mathcal{L}(V, \Lambda V)$ and $G \in \Lambda V \otimes V$ we define $G \cdot F \in \Lambda V$ by extending linearly

$$(f \otimes v) \cdot F = f F(v), \qquad f \in \Lambda V, v \in V,$$

where on the right hand side we have the multiplication of elements in $\Lambda V$.

Furthermore if $\Psi_1, \Psi_2 \in \mathcal{G}(V)$, $K = K_1 \otimes K_2 \in \Lambda V \otimes V$ and $F \in \mathcal{L}(V, \Lambda V)$ we write

$$\langle K(\Psi_1), F(\Psi_2) \rangle = K_1(\Psi_1) F(K_2)(\Psi_2) \in \mathcal{A}. \tag{30}$$

The definition (30) of the pairing $\langle K(\Psi_1), F(\Psi_2) \rangle$ can be extended by linearity to any $K \in \Lambda V \otimes V$.

**Theorem 32.** *For the global solution to equation (23) we have*

$$\omega(\Psi_t(G)\Psi_0(H)) = \omega(\Psi_0(G)\Psi_0(H)) + \int_0^t \omega\left[\Psi_s\left(\partial_R G \cdot F + \frac{1}{2} Q_C(\partial_R^2 G)\right)\Psi_0(H)\right] \mathrm{d}s \tag{31}$$

*for all $G, H \in \Lambda V$ and $t \geqslant 0$.*

**Proof.** We will prove the statement for $H = 1$ and leave to the reader to generalize it to any $H$. Fix any $T \geqslant 0$. Note that

$$\|B_t - B_s\|_{\mathcal{G}(V)} \lesssim |t-s|^{1/2},$$

for all $0 \leqslant s < t \leqslant T$. From that it follows that, if $F$ in (23) is a polynomial of degree $\deg(F)$ by Theorem 14, we have

$$\|(\Psi_t - \Psi_s) - (B_t - B_s)\|_{\mathcal{G}(V)} = \left\|\int_s^t F(\Psi_s) \mathrm{d}s\right\|_{\mathcal{G}(V)} \lesssim (1 + \|\Psi\|_{[0,T]})^{\deg(F)} |t-s|, \quad (32)$$

where $\|\Psi\|_{[0,T]} := \sup_{s \in [0,T]} \|\Psi_s\|_{\mathcal{G}(V)}$ which, by Theorem 31 and the continuity of $\Psi$ with respect to the time $t$, is finite for every $T \in \mathbb{R}_+$. By Taylor formula (see Lemma 10, Theorem 12 and Theorem 14) applied to polynomial $G$ of degree $\deg(G)$ we have

$$\left\|G(\Psi_r) - G(\Psi_s) - m_{\mathcal{A}}(\Psi_s \otimes (\Psi_r - \Psi_s))\partial_R G + \right.$$
$$\left. + \frac{1}{2}m_{\mathcal{A}}^2(\Psi_s \otimes (\Psi_r - \Psi_s) \otimes (\Psi_r - \Psi_s))(\partial_R \partial_R G)\right\|_{\mathcal{A}} \lesssim (1 + \|\Psi\|_{[0,T]})^{\deg(F)(\deg(G)-2)} |r-s|^{3/2} \quad (33)$$

since $\|\Psi_t - \Psi_s\|_{\mathcal{G}(V)} \lesssim (1 + \|\Psi\|_{[0,T]})^{\deg(F)} |t-s|^{1/2}$ by (32). In a similar way it is possible to obtain

$$\|\langle \partial_R G(\Psi_u), F(\Psi_u)\rangle - \langle \partial_R G(\Psi_s), F(\Psi_u)\rangle\|_{\mathcal{A}} \lesssim (1 + \|\Psi\|_{[0,T]})^{\deg(F)(\deg(G)-2)} |u-s|^{1/2} \quad (34)$$

where $\langle \partial_R G(\Psi_s), F(\Psi_u)\rangle$ is defined by linear extension of the relation (30). On the other hand, by (32), we have

$$\begin{aligned} m_{\mathcal{A}}(\Psi_s \otimes (\Psi_r - \Psi_s))(\partial_R G) &= \int_s^r \langle \partial_R G(\Psi_s), F(\Psi_u)\rangle \mathrm{d}u \\ &\quad + m_{\mathcal{A}}(\Psi_s \otimes (B_r - B_s))(\partial_R G) \\ &= \int_s^r \langle \partial_R G(\Psi_u), F(\Psi_u)\rangle \mathrm{d}u + \langle \partial_R G(\Psi_s), (B_r - B_s)\rangle \\ &\quad + O_{\mathcal{A}}((1 + \|\Psi\|_{[0,T]})^{\deg(F)(\deg(G)-2)} |r-s|^{3/2}) \end{aligned} \quad (35)$$

In a similar way, using inequality (32) we obtain

$$m_{\mathcal{A}}^2(\Psi_s \otimes (\Psi_r - \Psi_s) \otimes (\Psi_r - \Psi_s))(\partial_R^2 G) = m_{\mathcal{A}}^2(\Psi_s \otimes (B_r - B_s) \otimes (B_r - B_s))(\partial_R^2 G) \\ + O_{\mathcal{A}}((1 + \|\Psi\|_{[0,T]})^{\deg(F)(\deg(G)-2)} |r-s|^{3/2}). \quad (36)$$

Furthermore, using the fact that $\partial_R G(\Psi_s)$ is independent of $(B_r - B_s)$ and $((B_r - B_s) \otimes (B_r - B_s))$ is independent of $\partial_R^2 G(\Psi_s)$ (this is due to the fact that $\Psi_s$ is a function of $\{\Psi_0, B_u | u \in [0,s]\}$ and, moreover, $B$ is independent of $\Psi_0$ and has independent increments, we have:

$$\omega(\langle \partial_R G(\Psi_s), (B_r - B_s)\rangle) = \langle \omega((B_r - B_s)), \omega(\partial_R G(\Psi_s))\rangle = 0$$

where on the r.h.s. we understand that averages w.r.t. the state $\omega$ are taken component-wise, and

$$\omega[\langle \partial_R^2 G(\Psi_s), (B_r - B_s) \otimes (B_r - B_s)\rangle_{V \otimes V}]$$
$$=\langle \omega[(\partial_R^2 G(\Psi_s))], \omega[((B_r - B_s) \otimes (B_r - B_s))]\rangle$$
$$=(r-s)\omega[Q_C(\partial_R^2 G(\Psi_s))].$$

By taking a partition $\{t_i\}_{i \in \{0,\ldots,n\}}$, of diameter $\rho_n \to 0$ as $n \to +\infty$, of $[0,t]$ and exploiting inequality (33) we obtain

$$\begin{aligned}
\omega(G(\Psi_t)) &= \omega(G(\Psi_0)) + \sum_{i=0}^{n+1} \omega(G(\Psi_{t_i}) - G(\Psi_{t_{i-1}})) \\
&= \omega(G(\Psi_0)) + \sum_{i=1}^{n+1} \omega(m_{\mathcal{A}}(\Psi_{t_{i-1}} \otimes (\Psi_{t_i} - \Psi_{t_{i-1}}))\partial_R G) \\
&\quad + \frac{1}{2}\omega(m_{\mathcal{A}}^2(\Psi_s \otimes (\Psi_r - \Psi_s) \otimes (\Psi_r - \Psi_s))(\partial_R^2 G)) + O(\rho_n^{1/2}) \\
&= \omega(G(\Psi_0)) + \omega\left(\int_0^t \langle \partial_R G(\Psi_u), F(\Psi_u)\rangle \mathrm{d}u\right) \\
&\quad + \frac{1}{2}\sum_{i=1}^{n+1}(t_i - t_{i-1})\omega[Q_C(\partial_R^2 G(\Psi_{t_i}))] + O(\rho_n^{1/2}),
\end{aligned}$$

where the constants in $O$ are proportional to $(1 + \|\Psi\|_{[0,T]})^{\deg(F)\deg(G)}$ and do not depend on the partition. Taking the limit $n \to +\infty$, and so $\rho_n \to 0$, and using the fact that $\omega[Q_C(\partial_R^2 G(\Psi_s))]$ is continuous in $s$, $(G(\Psi_s))$ being a continuous function from $\mathbb{R}_+$ to $\mathcal{A}$ since $\Psi \in C^0([0,T], \mathcal{G}(V)))$ we obtain the thesis. $\square$

### 3.4 Invariant measures

**Definition 33.** *We say that a continuous linear functional $\rho \colon \mathcal{A} \to \mathbb{C}$ is an invariant measure for the equation (23) if for any $G \in \Lambda V$ and any $t \in \mathbb{R}_+$ we have*

$$\rho(G(\Psi_t)) = \rho(G(\Psi_0)). \tag{37}$$

There is a simple condition, analogous to the Fokker–Planck equation for commutative SDEs, for checking the invariance of the functional $\rho$.

**Lemma 34.** *If $\Psi$ is a solution to (23) and $\rho \colon \mathcal{A} \to \mathbb{C}$ a continuous linear functional for which*

$$\rho(\mathcal{L}G(\Psi_0)) = 0, \qquad G \in \Lambda V, \tag{38}$$

*where*

$$\mathcal{L}G := \partial_R G \cdot F + \frac{1}{2}Q_C(\partial_R^2 G),$$

*then $\rho$ is an invariant measure for eq. (23) in the sense of Definition 33. In other words, eq. (38) is a sufficient condition for $\rho$ to be an invariant measure.*

**Proof.** Let $(v_A)_A$ be a (finite) linear basis of $\Lambda V$ and observe that we have $\mathcal{L}v_A = \sum_B \kappa_A^B v_B$ for a suitable family $(\kappa_A^B)_{A,B}$ of constants. Let $\mathfrak{P}_A(t) := \rho(v_A(\Psi_t))$, by Itô formula (31) we have that $t \in \mathbb{R}_+ \mapsto (\mathfrak{P}_A(t))_A$ is the unique solution of the system of ODEs

$$\begin{aligned}
\partial_t \mathfrak{P}_A(t) &= \partial_t \rho(v_A(\Psi_t)) = \rho((\mathcal{L}v_A)(\Psi_t)) \\
&= \sum_B \kappa_A^B \rho(v_B(\Psi_t)) = \sum_B \kappa_A^B \mathfrak{P}_B(t),
\end{aligned} \tag{39}$$

with initial condition $\mathfrak{P}_A(0) = \rho(v_A(\Psi_0))$. On the other hand condition (38) implies that

$$\sum_B \kappa_A^B \mathfrak{P}_B(0) = \sum_B \kappa_A^B \rho(v_B(\Psi_0)) = \rho(\mathcal{L} v_A(\Psi_0)) = 0,$$

which means that $(t \mapsto \mathfrak{P}_A(0))_A$ is also a solution of (38), and therefore by uniqueness of solutions to this system of ODEs we conclude that $\mathfrak{P}_A(t) = \mathfrak{P}_A(0)$ for all $t \geqslant 0$. Since the expectations of the form (37) are linear combinations of $(\mathfrak{P}_A(t))_A$ the proof is complete. $\square$

The following theorem gives a sufficient condition on the SDE (23) to have the invariant measure (42).

**Theorem 35.** *For any $U \in \Lambda_{\mathrm{even}} V$, consider the SDE*

$$\Psi_t(v) = \Psi_0(v) + \int_0^t \Psi_s(Av - \langle C\partial_R U, v\rangle)\mathrm{d}s + B_t(v) \qquad t \geqslant 0, v \in V, \tag{40}$$

*where, under the state $\omega$, $B$ is a Brownian motion with correlation $C: V \to V$ and $\Psi_0$ a Gaussian initial condition with correlation $C_A$ (defined in equation (27)) independent of $B$. Assume that the condition*

$$A^{\mathrm{T}} C_A - C_A A = 0 \tag{41}$$

*holds. Then linear functional $\tilde{\omega}: \mathcal{A} \to \mathbb{C}$ given by*

$$\tilde{\omega}(\cdot) := \omega(\cdot e^{-2U(\Psi_0)}), \tag{42}$$

*is an invariant measure for (40).*

**Proof.** The proof is a consequence of Lemma 34. Denote $\omega^{\Psi_0}(G) := \omega(G(\Psi_0))$ the law of $\Psi_0$, and observe that for $H \in \Lambda V \otimes V$ we have the integration by parts formula

$$\omega(m_{\mathcal{A}}(\Psi_0 \otimes \Psi_0)(H)) = \omega(\Psi_0(Q_{C_A}(\partial_R H))),$$

where $Q_C$ is defined in (29). Then for all $G \in \Lambda V$ we have

$$\omega[m_{\mathcal{A}}(\Psi_0 \otimes \Psi_0)((\mathbb{I} \otimes A)\partial_R G e^{-2U})] = \omega^{\Psi_0}[Q_{C_A}(\partial_R((\mathbb{I} \otimes A)\partial_R G e^{-2U}))]$$

$$= \omega^{\Psi_0}(Q_{C_A}((\mathbb{I} \otimes A \otimes \mathbb{I})\partial_R^2 G)e^{-2U}) - 2\omega^{\Psi_0}(Q_{C_A}((\mathbb{I} \otimes A \otimes \mathbb{I})(\partial_R G \otimes_2 \partial_R U))e^{-2U})$$

where $\otimes_2$ is defined as $(f \otimes v) \otimes_2 (g \otimes w) = (fg) \otimes v \otimes w$ for $f, g \in \Lambda V$ and $v, w \in V$. Observe that using the trace Tr on $\Lambda V \otimes V \otimes V$ given by $\mathrm{Tr}(f \otimes v \otimes w) = f\langle v, w\rangle$, the antisymmetry of $\partial_R^2 G \in \Lambda V \otimes V \otimes V$ in the last two factors and equality (28) we have

$$Q_{C_A}((\mathbb{I} \otimes A \otimes \mathbb{I})\partial_R^2 G) = \mathrm{Tr}((\mathbb{I} \otimes A \otimes C_A)\partial_R^2 G)$$

$$= \mathrm{Tr}((\mathbb{I} \otimes -C_A A \otimes \mathbb{I})\partial_R^2 G) = \frac{1}{2}\mathrm{Tr}((\mathbb{I} \otimes \mathbb{I} \otimes (A^{\mathrm{T}} C_A + C_A A))\partial_R^2 G)$$

$$= -\frac{1}{2}\mathrm{Tr}((\mathbb{I} \otimes \mathbb{I} \otimes C)\partial_R^2 G) = -\frac{1}{2} Q_C(\partial_R^2 G).$$

Moreover, from (41) we have also

$$Q_{C_A}((\mathbb{I} \otimes A \otimes \mathbb{I})(\partial_R G \otimes_2 \partial_R U)) = \mathrm{Tr}((\mathbb{I} \otimes A \otimes C_A)(\partial_R G \otimes_2 \partial_R U)) = \mathrm{Tr}(\partial_R G \otimes_2 (A^{\mathrm{T}} C_A \partial_R U))$$

$$= \mathrm{Tr}\left(\partial_R G \otimes_2 \left(\frac{A^{\mathrm{T}} C_A + C_A A}{2} \partial_R U\right)\right) + \mathrm{Tr}\left(\partial_R G \otimes_2 \left(\frac{A^{\mathrm{T}} C_A - C_A A}{2} \partial_R U\right)\right)$$

$$= -\frac{1}{2}\mathrm{Tr}(\partial_R G \otimes_2 (C\partial_R U)) = -\frac{1}{2}\langle \partial_R G, C\partial_R U\rangle,$$

since we assumed that $A^{\mathrm{T}} C_A - C_A A = 0$ and used (28) again. Therefore

$$\omega(m_{\mathcal{A}}(\Psi_0 \otimes \Psi_0)((\mathbb{I} \otimes A)(\partial_R G) e^{-2U})) = -\frac{1}{2}\omega^{\Psi_0}(Q_C(\partial_R^2 G) e^{-2U}) + \omega^{\Psi_0}(\langle \partial_R G, C \partial_R U \rangle e^{-2U}).$$

Now observing that

$$\mathcal{L}G(\Psi_0) = m_{\mathcal{A}}(\Psi_0 \otimes \Psi_0)((\mathbb{I} \otimes A)(\partial_R G)) + \Psi_0\left(-\langle \partial_R G, C \partial_R U \rangle + \frac{1}{2} Q_C(\partial_R^2 G)\right)$$

we deduce that

$$\omega(\mathcal{L}G(\Psi_0) e^{-2U(\Psi_0)}) = \omega(m_{\mathcal{A}}(\Psi_0 \otimes \Psi_0)((\mathbb{I} \otimes A)(\partial_R G)) e^{-2U(\Psi_0)}) +$$
$$+ \omega\left[\Psi_0\left(-\langle \partial_R G, C \partial_R U \rangle + \frac{1}{2} Q_C(\partial_R^2 G)\right) e^{-2U(\Psi_0)}\right] = 0. \qquad \square$$

## 3.5 Long-time behavior for small non-linearity

We will now investigate the existence of solutions which are globally bounded in time. In the commutative setting, non-linear equations can have globally bounded solutions only if the non-linearity stays uniformly small or if it shows some coercivity. Like all notions of positivity, also coercivity however does not apply well in the Grassmann setting. The only kind of coercive term we have identified is a linear drift with a negative sign. This is a very mild coercive term, but it turns out to be enough provided the non-linearity is small enough. Consider the equation

$$\Psi_t(v) = \Psi_0(v) + \int_0^t (\Psi_s(Av) + \lambda \Psi_s(F(v))) \mathrm{d}s + B_t(v), \qquad t \geqslant 0, v \in V, \tag{43}$$

for $\Psi \in C([0,T]; \mathrm{Hom}(V, \mathcal{A}))$, with $\lambda \geqslant 0$ and $F \in \mathrm{Hom}(V, \Lambda V)$ where (as in Section 3.1) $A$ is an operator with eigenvalues having strictly negative real part less then $-\lambda_A < 0$. In this setting we introduce the notion of stationary solution to equation (43), extending the one we defined for the linear case in Section 3.1.

**Definition 36.** *We say that $\Psi_t^{\mathfrak{s}} \in C^0(\mathbb{R}, \mathcal{G}(V))$ is a stationary solution to equation (43) (extended to all $t \in \mathbb{R}$) of norm at most $K \in \mathbb{R}_+$ if the following two conditions hold:*

1. $\sup_{t \in \mathbb{R}} \|\Psi_t^{\mathfrak{s}}\|_{\mathcal{G}(V)} \leqslant K$,

2. $(\Psi_t^{\mathfrak{s}})_{t \in \mathbb{R}}$ *solves the following integral equation*

$$\Psi_t^{\mathfrak{s}} = \lambda \int_{-\infty}^t \Psi_\tau^{\mathfrak{s}}(F(e^{A\tau} v)) \mathrm{d}\tau + B_t^A \qquad t \in \mathbb{R} \tag{44}$$

*where $B_t^A(v) = \Xi(e^{A(t-\cdot)}(v))$ is the OU motion introduced in Section 3.1.*

Hereafter we write

$$L = L(A, C) = \sup_{t \in \mathbb{R}} \|B_t^A\|_{\mathcal{G}(V)} < \infty,$$

and we denote by $\mathfrak{K}_K \subset C^0(\mathbb{R}, \mathcal{G}(V))$ the set of $\Phi \in C^0(\mathbb{R}, \mathcal{G}(V))$ such that

$$\sup_{t \in \mathbb{R}} \|\Phi_t\|_{\mathcal{G}(V)} \leqslant K,$$

for some constant $K > 0$.

**Theorem 37.** *Assume $K \geqslant 4L(A,C)$ and suppose that all the eigenvalues of $A$ have negative real part that are less or equal than $-\lambda_A < 0$, then there exists $\lambda_0 = \lambda_0(K, L, F, \lambda_A) > 0$ depending on $K$, $L$ and $F$ such that if $|\lambda| \leqslant 2\lambda_0$ there exists a unique stationary solution of norm at most $K$ to equation (43).*

**Proof.** It is a simple application of the Banach fixed-point theorem by noting that the map

$$\mathcal{K}_\lambda(\Psi)_t := \lambda \int_{-\infty}^t \Psi_\tau(F(e^{A(t-\tau)}v)) d\tau + B_t^A, \qquad t \in \mathbb{R},$$

is a contraction in $\mathfrak{K}_K$ for $\lambda$ small. Indeed, note that $B^A \in \mathfrak{K}_K$ and if $\Psi \in \mathfrak{K}_K$ we have, by Theorem 14,

$$\|\mathcal{K}_\lambda(\Psi)_t\|_{\mathcal{G}(V)} \leqslant |\lambda| \|F\|_{\Lambda_\pi V}(1+K)^{\deg(F)} \int_{-\infty}^t e^{-\lambda_A(t-\tau)} d\tau + L(A,C),$$

where $\|F\|_{\Lambda_\pi V}$ is defined by (15). Therefore $\mathcal{K}_\lambda$ maps $\mathfrak{K}_K$ into itself provided $|\lambda| \leqslant 2\lambda_0 := \lambda_0(K, L, F, \lambda_A)$ with $\lambda_0$ such that

$$2\lambda_0 \frac{\|F\|_{\Lambda_\pi V}(1+K)^{\deg(F)}}{\lambda_A} + 2L \leqslant K. \tag{45}$$

Furthermore

$$\|\mathcal{K}_\lambda(\Psi_1)_t - \mathcal{K}_\lambda(\Psi_2)_t\|_{\mathcal{G}(V)} \leqslant \frac{|\lambda|}{\lambda_A}(\deg(F))\|F\|_{\Lambda_\pi V}(1+K)^{\deg(F)-1}\|\Psi_{1,t} - \Psi_{2,t}\|_{\mathcal{G}(V)}.$$

Therefore, if we take $\lambda_0$ small enough to satisfy also

$$\frac{2\lambda_0}{\lambda_A}\|F\|_{\Lambda_\pi V}(\deg(F))(1+K)^{\deg(F)-1} < 1 \tag{46}$$

the map $\mathcal{K}_\lambda$ is a strict contraction. □

**Remark 38.** The function $\lambda_0$ in Theorem 37 can be taken to be a decreasing function of the parameters $L$ and $K$.

We introduce an approximation for the solution to equation (44). Let $X \in \mathcal{G}(V)$ independent of $B_t^A$ and consider the element $\Psi^X_{-T,t} \in C^0(\mathbb{R}, \mathcal{G}(V))$ which is $\Psi^X_{-T,t} = X$ if $t \leqslant -T$, for some $T > 0$, and is a solution to the SDE

$$\Psi^X_{-T,t} = X(e^{A(t+T)}v) + \lambda \int_{-T}^t F(\Psi^X_{-T,t}(e^{A(t-\tau)}v)) + B_t^A(e^{A(t+T)}v) - B_{-T}^A(e^{A(t+T)}v),$$

when $t > -T$ and extended as $\Psi^X_{-T,t} = X$ for $t \leqslant -T$.

**Lemma 39.** *Under the assumptions of Theorem 37, take $\|X\|_{\mathcal{G}(V)} \leqslant K/8$, then, for any $|\lambda| \leqslant \lambda_0$, we have that $\Psi^X_{-T,\cdot} \in \mathfrak{K}_K$. Moreover $\Psi^X_{-T,\cdot}$ converges to $\Psi^{\mathfrak{s}}$ in $\mathcal{G}(V)$ as $T \to +\infty$, uniformly on compact subsets of $\mathbb{R}$.*

**Proof.** Let $\tau \in \mathbb{R}$ and $R_\tau := \sup_{t \in (-\infty, \tau]} \|\Psi^X_{-T,t}\|_{\mathcal{G}(V)}$. Theorem 14 gives the bound

$$\begin{aligned} R_\tau &\leqslant \|X\|_{\mathcal{G}(V)} + \|B_{-T}^A\|_{\mathcal{G}(V)} + \\ &\quad + \|F\|_{\Lambda_\pi V}\frac{\lambda_0}{\lambda_A}(1+R_\tau)^{\deg(F)} + L(A,C) \\ &\leqslant \frac{3K}{8} + \|F\|_{\Lambda_\pi V}\frac{\lambda_0}{\lambda_A}(1+R_\tau)^{\deg(F)} + L(A,C) \end{aligned}$$

where we used that $\|B^A_{-T}\|_{\mathcal{G}(V)} \leqslant L \leqslant K/4$ and $\|X\|_{\mathcal{G}(V)} \leqslant K/8$. On the other hand by definition of $\lambda_0$, we have

$$\frac{3K}{8} + \|F\|_{\Lambda_\pi V} \frac{\lambda_0}{\lambda_A}(1+K)^{\deg(F)} + L(A,C) \leqslant \frac{7}{8}K < K$$

for any $\tau > 0$. Consider the set $Z = \{\tau \in \mathbb{R} : R_\tau \leqslant K\}$. The set $Z$ is non-empty, since $R_\tau \leqslant \frac{1}{2}K$ when $\tau \leqslant -T$, closed since $\tau \mapsto R_\tau$ is continuous and open since $\tau \mapsto R_\tau$ is an increasing continuous function of $\tau$ and if $\tau \in Z$ then $R_{\tau+\varepsilon} < K$ for $\varepsilon > 0$ small enough so $(\tau-\varepsilon, \tau+\varepsilon) \in Z$. Then we must have $Z = \mathbb{R}$ and

$$\sup_{t \in \mathbb{R}} \|\Psi^X_{-T,t}\|_{\mathcal{G}(V)} \leqslant K,$$

for $|\lambda| \leqslant \lambda_0$. For the convergence we argue as follows. If $t \geqslant -T$, by Theorem 14, we have that

$$\|\Psi^X_{-T,t} - \Psi^{\mathfrak{s}}_t\|_{\mathcal{G}(V)} \lesssim e^{-\lambda_A(t+T)}(\|X\|_{\mathcal{G}(V)} + \|B^A_{-T}\|_{\mathcal{G}(V)}) +$$
$$+ |\lambda| \|F\|_{\Lambda_\pi V}(1+K)^{\deg(F)} \int_{-\infty}^{-T} e^{-\lambda_A(t-\tau)} \mathrm{d}\tau$$
$$+ \alpha \int_{-T}^{t} e^{-\lambda_A(t-\tau)} \|\Psi^X_{-T,\tau} - \Psi^{\mathfrak{s}}_\tau\| \mathrm{d}\tau.$$

where $\alpha = |\lambda| \|F\|_{\Lambda_\pi V}(\deg(F))(1+K)^{\deg(F)-1}$. By Gronwall inequality, we obtain

$$\|\Psi^X_{-T,t} - \Psi^{\mathfrak{s}}_t\|_{\mathcal{G}(V)} \leqslant \left(\|X\|_{\mathcal{G}(V)} + \frac{K}{3} + \frac{|\lambda| \|F\|_{\Lambda_\pi V}(1+K)^{\deg(F)}}{\lambda_A}\right) \times$$
$$\times \left(e^{-\lambda_A(t+T) + \alpha t} \int_{-T}^{t} e^{-\alpha \tau} d\tau\right)$$
$$\lesssim e^{-\lambda_A(t+T)} e^{\alpha t}(e^{\alpha T} - e^{-\alpha t}) \lesssim e^{-(\lambda_A - \alpha)(t+T)}.$$

Since $\lambda_A - \alpha > 0$, by the condition (46), we have that

$$\|\Psi^X_{-T,t} - \Psi^{\mathfrak{s}}_t\|_{\mathcal{G}(V)} \lesssim e^{-(\lambda_A - \alpha)T},$$

uniformly in $T$ when $t > P > -T$ for any such fixed $P \in \mathbb{R}$. This proves that $\Psi^X_{-T,t}$ converges to $\Psi^{\mathfrak{s}}_t$ uniformly on compact sets. $\square$

We are now in position to prove the main result of this section, the stochastic quantization of a Gibbs measure on Grassmann Gaussian variables.

**Theorem 40.** *Assume that under the state $\omega$ the random variable $X$ is a Grassmann Gaussian with covariance $C_A$ and $B^A$ an OU motion independent of $X$. Assume that $K$ is a constant such that $L \leqslant K/3$ and $\|X\|_{\mathcal{G}(V)} \leqslant K/8$. Let $U \in \Lambda V_{\text{even}}$ and let*

$$F(v) = -\frac{\lambda}{2}\langle C\partial_R U, v\rangle$$

*for $v \in V$. Then $\tilde{\omega}(\cdot) = \omega(\cdot e^{-\lambda U(X)})$ is an invariant measure for the equation equation (43) and for all $|\lambda| \leqslant \lambda_0(K, L, F, \lambda_A)$, for any $t \in \mathbb{R}$ and for any $G \in \Lambda V$, we have that*

$$\omega(G(X)e^{-\lambda U(X)}) = \omega(G(\Psi^{\mathfrak{s}}_t))\omega(e^{-\lambda U(X)}).$$

**Proof.** Note that $(\Psi^X_{-T,t})_{t \geqslant -T}$ is a solution of the SDE (40) starting at $-T$ with the value $X$ and with potential $\lambda U/2$. Therefore by Theorem 35 we have that, for any $t \geqslant -T$,

$$\omega(G(X)e^{-\lambda U(X)}) = \omega(G(\Psi^X_{-T,-T})e^{-\lambda U(X)}) = \omega(G(\Psi^X_{-T,t})e^{-\lambda U(X)}).$$

By Lemma 39 we can take the limit at the right hand side as $T \to +\infty$, obtaining

$$\omega(G(X)e^{-U(X)}) = \omega(G(\Psi_t^{\mathfrak{s}})e^{-U(X)}).$$

Since (by construction) $X$ is independent of $B_t^A$ under $\omega$ (and therefore under $\tilde{\omega}$) and since $\Psi_t^{\mathfrak{s}}$ is a function only of $B^A$ we get

$$\omega(G(X)e^{-U(X)}) = \omega(G(\Psi_t^{\mathfrak{s}}))\omega(e^{-U(X)}). \qquad \square$$

**Remark 41.** Note that by taking $G = \lambda U$ we obtain

$$1 = \omega(e^{\lambda U(\Psi_t^{\mathfrak{s}})})\omega(e^{-\lambda U(X)})$$

which implies in particular that both $\omega(e^{\lambda U(\Psi_t^{\mathfrak{s}})})$ and $\omega(e^{-\lambda U(X)})$ are nonzero and that

$$\omega(e^{-\lambda U(X)}) = [\omega(e^{\lambda U(\Psi_t^{\mathfrak{s}})})]^{-1}$$

$$\frac{\omega(G(X)e^{-\lambda U(X)})}{\omega(e^{-\lambda U(X)})} = \omega(G(\Psi_t^{\mathfrak{s}})) = \omega(e^{\lambda U(\Psi_t^{\mathfrak{s}})})\omega(G(X)e^{-\lambda U(X)}).$$

**Remark 42.** Let now $\Psi^{\mathfrak{s},\sigma}$ be the stationary solution corresponding to the potential $\sigma U$ for $\sigma \in [0, \lambda]$. Then

$$\frac{\mathrm{d}}{\mathrm{d}\sigma}\log \omega(e^{-\sigma U(X)}) = -\frac{\omega(U(X)e^{-\sigma U(X)})}{\omega(e^{-\sigma U(X)})} = -\omega(U(\Psi_t^{\mathfrak{s},\sigma}))$$

and integrating $\sigma$ gives a remarkable formula for the logarithm of the partition function $\omega(e^{-\lambda U(X)})$ in terms of the parametric stochastic quantization $(\Psi_t^{\mathfrak{s},\sigma})_{\sigma \in [0,1]}$:

$$-\log \omega(e^{-\lambda U(X)}) = \int_0^\lambda \omega(U(\Psi_t^{\mathfrak{s},\sigma}))\mathrm{d}\sigma.$$

## 4 Infinite dimensional SDEs

We want to study SDEs of the form

$$\mathrm{d}\Psi_t = (A\Psi_t + \lambda F(\Psi_t))\mathrm{d}t + \mathrm{d}B_t \qquad (47)$$

in the case where $\Psi$ takes values in $\mathcal{G}(V)$ with $V = \mathcal{S}(\mathbb{R}^d) \otimes \mathbb{R}^n$, or $\mathcal{S}(\mathbb{R}^d) \otimes \mathbb{C}^n$, or $C^\infty(\mathbb{T}^d) \otimes \mathbb{R}^n$ for some $n \in \mathbb{N}_0$ and with $A$ a (possibly unbounded) linear operator on $V$. Heuristically, $d$ represents the dimension of base space and $n$ is the number of Grassmannian fields considered. Here, we discuss in detail the case $V = \mathcal{S}(\mathbb{R}^d) \otimes \mathbb{R}^n$, but the considerations in this section can be generalized to the other cases, or more generally, to any nuclear space of smooth functions (see Remark 43).

The main difference between the finite and infinite dimensional case for Grassmann SDEs is that in general there is no natural topology on $\mathcal{G}(V)$. The space $V = \mathcal{S}(\mathbb{R}^d) \otimes \mathbb{R}^n$ is a Fréchet nuclear space. Since not every element of $\mathcal{G}(V)$ is continuous with respect to this topology, we consider here the space $\mathcal{G}^{-\infty}(V) \subset \mathcal{G}(V)$, the subset of $\mathcal{G}(V)$ for which $X \in \mathcal{G}^{-\infty}(V)$ if and only if $X$ restricted to $V \subset \Lambda V$ is a continuous linear map from $V$ to $\mathcal{A}$. Hereafter we always consider $X \in \mathcal{G}^{-\infty}(V)$ both as an homomorphism from the algebra $\Lambda V$ into $\mathcal{A}$, and as a continuous linear map from $V$ into $\mathcal{A}$.

Since $V$ is a Fréchet nuclear space whose dual is nuclear, the set of continuous linear maps from $V$ into $\mathcal{A}$ can be identified with the tensor product $V^* \hat{\otimes} \mathcal{A} = (\mathcal{S}'(\mathbb{R}^d) \otimes \mathbb{R}^n) \hat{\otimes} \mathcal{A}$, see, e.g., Proposition 50.5 in [174]. Here $V^*$ is the topological dual of $V$, the Schwartz space of tempered distribution $\mathcal{S}'(\mathbb{R}^d) \otimes \mathbb{R}^n$ and the symbol $\hat{\otimes}$ in the product $V^* \hat{\otimes} \mathcal{A}$ means that we take the closure of $V^* \otimes \mathcal{A}$ with respect to any natural topology on $V^* \otimes \mathcal{A}$ being all of them equivalent since $V^*$ is nuclear (see, e.g., Theorem 50.1 in [174]). In other words we can look upon $\mathcal{G}^{-\infty}(V)$ as the subset of elements of $(\mathcal{S}'(\mathbb{R}^d) \otimes \mathbb{R}^n) \otimes \mathcal{A}$ that can be extended to a homomorphism from $\Lambda V$ into $\mathcal{A}$.

**Remark 43.** The identification of the set of continuous linear maps from $V$ into $\mathcal{A}$ with $V^* \hat{\otimes} \mathcal{A}$ does not hold if $V$ is a generic topological vector space, for example if $V$ is a Banach space. Indeed in the case where $V$ is a infinite dimensional Banach space the space $V^* \otimes \mathcal{A}$ has not unique natural norm. In particular the completion $V^* \otimes_\varepsilon \mathcal{A}$ with respect the injective norm (which is the weaker of them) we have $V^* \otimes_\varepsilon \mathcal{A} \subsetneq \mathcal{L}(V, \mathcal{A})$ (see, e.g., [148]). For this reason we prefer to work in the case where $\mathcal{L}(V, \mathcal{A}) = V^* \hat{\otimes} \mathcal{A}$ which is for example when $V$ is a Fréchet nuclear space.

Using a reasoning similar to the one exploited in Lemma 7, it is possible to prove that $\mathcal{G}^{-\infty}(V)$ is a closed subset of $V^* \otimes \mathcal{A}$ with respect to the strong topology (of both $V^* = \mathcal{S}'(\mathbb{R}^d) \otimes \mathbb{R}^n$ and $\mathcal{A}$), i.e. if $(X_\nu)_\nu$ is a net in $\mathcal{G}^{-\infty}(V)$ converging to $X$ in the topology of $V^* \otimes \mathcal{A}$, then $X$ can be extended to an homomorphism from $\Lambda V$ into $\mathcal{A}$. Working on $\mathcal{G}^{-\infty}(V)$ is not so easy since $(\mathcal{S}'(\mathbb{R}^d) \otimes \mathbb{R}^n) \otimes \mathcal{A}$ is not even a Fréchet space. For this reason we introduce now some subsets of $\mathcal{G}^{-\infty}(V)$ (which are analogous to Besov spaces) which admit a stronger metric topology.

**Remark 44.** In what fallow we use the notation if $X \in \mathcal{G}^{-\infty}(V)$ we write, for any $k \in \{1, ..., n\}$, the map $X^k \in \mathcal{G}^{-\infty}(\mathcal{S}(\mathbb{R}^n))$ as follows
$$X^k(f) = X(f \otimes e_k)$$
where $\{e_k\}_{k=1,...,n}$ is the standard orthonormal basis of $\mathbb{R}^n$.

We recall here only the definition of the Besov norm for functions taking values in a Banach space, further references and details are given in Appendix B (see also [9, 10]). Let $K_i \in \mathcal{S}(\mathbb{R}^d)$, with $i \geqslant -1$, be the functions corresponding to the Littlewood–Paley blocks used in the definition of Besov spaces $B^s_{p,q} := B^s_{p,q}(\mathbb{R}^d, \mathbb{R})$, where $s \in \mathbb{R}$ and $p, q \in [1, +\infty]$. If $X \in \mathcal{G}^{-\infty}(V)$ we define
$$\Delta_i X(x) = (\Delta_i X^1(x), ..., \Delta_i X^k(x)) := (X(K_i(x - \cdot) \otimes e_k))_{k=1,...,n}, \qquad x \in \mathbb{R}^d, i \geqslant -1$$
where $\{e_k\}_{k=1,...n}$ is the standard basis of $\mathbb{R}^n$. The function $\Delta_i X^k : \mathbb{R}^d \to \mathcal{A}$ can be identified with an element $\Delta_i X^k \in \mathcal{S}(\mathbb{R}^d) \hat{\otimes} \mathcal{A} = \mathcal{S}(\mathbb{R}^d, \mathcal{A})$ and we define
$$\|X\|_{B^s_{p,q,n}(\mathbb{R}^d, \mathcal{A})} := \sum_{k=1}^n \left[ \sum_{i \geqslant -1} 2^{qis} \left( \int_{\mathbb{R}^d} \|\Delta_i X^k(x)\|^p_\mathcal{A} dx \right)^{q/p} \right]^{1/q}.$$

Given $s \in \mathbb{R}$, we say that $X \in \mathcal{G}^s(V)$ when $X \in \mathcal{G}^{-\infty}(V)$ and $\|X\|_{B^s_{\infty,\infty}(\mathbb{R}^d, \mathcal{A})} < +\infty$. Thus we introduce the natural metric $d_{\mathcal{G}^s(V)}(X, Y) := \|X - Y\|_{\mathcal{G}^s(V)} := \|X - Y\|_{B^s_{\infty,\infty,n}(\mathbb{R}^d, \mathcal{A})}$. Hereafter we always use the notation $\mathcal{C}^s_n(\mathbb{R}^d, \mathcal{A}) = B^s_{\infty,\infty,n}(\mathbb{R}^d, \mathcal{A})$ and $\mathcal{C}^s = \mathcal{C}^s_1(\mathbb{R}^d, \mathbb{R})$. Moreover $\mathcal{C}^{s*}$ will denote the topological dual of $\mathcal{C}^s$.

**Lemma 45.** *The set $\mathcal{G}^s(V)$ with the metric $d_{\mathcal{G}^s(V)}$ is a complete metric space.*

**Proof.** The proof is analogous to that of Lemma 7. □

If $X \in \mathcal{G}^s(V)$ we can extend $X$ to a linear map $\mathcal{L}((\mathcal{C}^{s*} \otimes \mathbb{R}^n), \mathcal{A})$, i.e. to a linear map from the space $\mathcal{C}^{s*} \otimes \mathbb{R}^n$ (where $\mathcal{C}^{s*} = \mathcal{C}^{s*}(\mathbb{R}^d, \mathbb{R})$ is the (strong topological) dual of $\mathcal{C}^s$ and we identify $\mathbb{R}^n$ with its dual using its standard basis) taking values into the algebra $\mathcal{A}$. This extension can be obtained in the following way: If $g \in \mathcal{C}^{s*} \otimes \mathbb{R}^n$, for any sequence $(g_\ell)_\ell \subset \mathcal{S}(\mathbb{R}^d) \otimes \mathbb{R}^n$ converging to $g$ with respect to the weak* topology of $\mathcal{C}^{s*} \otimes \mathbb{R}^n$ we have that $(X(g_\ell))_\ell$ converges in $\mathcal{A}$ to a unique element

$$X(g) := \lim_{n \to +\infty} X(g_n).$$

The map $X$ now defined on all of $\mathcal{C}^{s*} \otimes \mathbb{R}^n$ is continuous. In other words, if $X \in \mathcal{G}^s(V)$ then it can be identified with an element of $\mathcal{L}((\mathcal{C}^{s*} \otimes \mathbb{R}^n), \mathcal{A})$. Furthermore it is simple to prove, using suitable $\mathcal{S}(\mathbb{R}^d) \otimes \mathbb{R}^n$ approximations, that, if $g, h \in \mathcal{C}^{s*} \otimes \mathbb{R}^n$, then we have

$$X(g)X(h) = -X(h)X(g). \tag{48}$$

In the case $s > 0$, for any $x \in \mathbb{R}^d$ the Dirac delta distribution $\delta_x$ is in $\mathcal{C}^{s*}$, this implies that we can identify the linear maps $X^k$ with some function in $\mathcal{C}^s(\mathbb{R}^d, \mathcal{A})$ defined as, using an abuse of notation,

$$X^k(x) := X(\delta_x \otimes e_k).$$

Using this identification we have that for any $f_n \in \mathcal{S}(\mathbb{R}^d)$, $X(f \otimes e_k) = \int_{\mathbb{R}^d} X^k(x) f(x) \mathrm{d}x$, where the integral is taken in Bochner sense.

The relation (48) has important consequences. For any positive integer $\ell \in \mathbb{N}_0$, let

$$\Lambda_\pi^\ell(\mathcal{C}^{s*} \otimes \mathbb{R}^n) \subset \underbrace{(\mathcal{C}^{s*} \otimes \mathbb{R}^n) \otimes_\pi (\mathcal{C}^{s*} \otimes \mathbb{R}^n) \otimes_\pi \cdots \otimes_\pi (\mathcal{C}^{s*} \otimes \mathbb{R}^n)}_{\ell \text{ times}}$$

be the closure of $\Lambda_\pi^\ell(\mathcal{C}^{s*} \otimes \mathbb{R}^n)$ in $\otimes_\pi^\ell(\mathcal{C}^{s*} \otimes \mathbb{R}^n)$ using the natural inclusion $i_{\Lambda^\ell(\mathcal{C}^{s*} \otimes \mathbb{R}^n)}$ of $\Lambda^\ell(\mathcal{C}^{s*} \otimes \mathbb{R}^n)$ in $\otimes^\ell(\mathcal{C}^{s*} \otimes \mathbb{R}^n)$ (see the discussion after Lemma 10 for more details). Here $\otimes_\pi$ is the projective tensor product of Banach spaces whose norm is defined as follows: If $f \in \otimes^\ell(\mathcal{C}^{s*} \otimes \mathbb{R}^n)$ then

$$\|f\|_{\otimes_\pi^\ell(\mathcal{C}^{s*} \otimes \mathbb{R}^n)} := \inf \left\{ \sum_{i=1}^k \prod_{j=1}^\ell \|f_j^i\|_{(\mathcal{C}^{s*} \otimes \mathbb{R}^n)}, \text{where } f = \sum_{i=1}^k f_1^i \otimes \cdots \otimes f_\ell^i \right\}.$$

**Lemma 46.** *If $X \in \mathcal{G}^s(V)$, then, for every $\ell \in \mathbb{N}_0$, it can be extended to a continuous homomorphism from $\bigoplus_{i=1}^\ell \Lambda_\pi^i(\mathcal{C}^{s*} \otimes \mathbb{R}^n)$ into $\mathcal{A}$. Here $\bigoplus_{i=1}^\ell \Lambda_\pi^i(\mathcal{C}^{s*} \otimes \mathbb{R}^n)$ is considered as equipped with the product obtained from the standard product of $\bigoplus_{i=1}^\infty \Lambda_\pi^i(\mathcal{C}^{s*} \otimes \mathbb{R}^n)$ via the projection onto $\bigoplus_{i=1}^\ell \Lambda_\pi^i(\mathcal{C}^{s*} \otimes \mathbb{R}^n)$.*

**Proof.** The proof is based on the fact that $X$ can be identified with an element of $\mathcal{L}(\mathcal{C}^{s*}, \mathcal{A})$ and on the relation (48). □

We want now to define a natural notion of function $F$ of an element $X \in \mathcal{G}^s(V)$. For example, if $s > 0$, and writing $Y = F(X)$, we want to be able to take $F = P_\ell$ of the form

$$Y^j(x) := Y(\delta_x \otimes e_j) =$$

$$= \sum_{k_1, \ldots, k_\ell = 1}^n \sum_{(h_1, \ldots, h_\ell) \in D} K^j_{(h_1, \ldots, h_\ell), k_1, \ldots, k_\ell} X^{k_1}(x + h_1) X^{k_2}(x + h_2) \cdots X^{k_\ell}(x + h_\ell), \tag{49}$$

where $D \subset \mathbb{R}^{d\ell}$ is some finite set and $K^j_{(h_1,...,h_\ell),k_1,...,k_\ell} \in \mathbb{R}$ and $j \in \{1,...,k\}$ ($\ell \in \mathbb{N}_0$ in this case represent the degree of the homogeneous polynomial $P_\ell$). When $\ell$ is odd, the function $P_\ell$ defined above maps $P_\ell(\mathcal{G}^s(V)) \subset \mathcal{C}^s_n(\mathbb{R}^d, \mathcal{A})$ and furthermore for any $g, h \in \mathcal{C}^{s*}$ and $k, p \in \{1,...,n\}$ we have

$$Y^k(g)Y^p(h) = -Y^p(h)Y^k(g)$$
$$Y^k(g)X^p(h) = -X^p(h)Y^k(g). \tag{50}$$

We want to consider some more general setting than the described in equation (49) that maintains the important property (50).

Consider $P_\ell \in \mathcal{L}(\mathcal{S}(\mathbb{R}^d) \otimes \mathbb{R}^n, \Lambda^\ell_\pi(\mathcal{C}^{s*} \otimes \mathbb{R}^n))$, then $P_\ell$ can be seen as an element of $(\mathcal{S}'(\mathbb{R}^d) \otimes \mathbb{R}^n) \otimes \Lambda^\ell_\pi(\mathcal{C}^{s*} \otimes \mathbb{R}^n)$. Using the same technique as before we can introduce the Banach space $\mathcal{C}^s_n(\mathbb{R}^d, \Lambda^\ell_\pi(\mathcal{C}^{s*} \otimes \mathbb{R}^n)) \subset (\mathcal{S}'(\mathbb{R}^d) \otimes \mathbb{R}^n) \otimes \Lambda^\ell_\pi(\mathcal{C}^{s*} \otimes \mathbb{R}^n)$ equipped with its natural norm. In this case we will use the notation

$$P_\ell(X) := X \circ P_\ell, \tag{51}$$

where $(X \circ P_\ell)(v) = X(P_\ell(v)) \in \mathcal{A}$ for any $v \in \mathcal{S}(\mathbb{R}^d) \otimes \mathbb{R}^n$. Consider $u \in \mathbb{N}_0$, if $F \in \mathcal{C}^s_n(\mathbb{R}^d, \bigoplus_{j=0}^u \Lambda^{2j+1}_\pi(\mathcal{C}^{s*} \otimes \mathbb{R}^n))$ there are some $P_\ell \in \mathcal{C}^s_n(\mathbb{R}^d, \Lambda^\ell_\pi(\mathcal{C}^{s*} \otimes \mathbb{R}^n))$ (for odd $\ell \in \mathbb{N}$) such that

$$F := \sum_{j=0}^u P_{2j+1}.$$

Thus we define $F(X) = \sum_{j=0}^u P_{2j+1}(X)$. The function $F$ defined above can be interpreted as an odd polynomial of degree $2u+1$ in the Grassmannian random variable $X$.

**Definition 47.** *We say that the linear function $F \in \mathcal{C}^s_n(\mathbb{R}^d, \bigoplus_{j=0}^u \Lambda^{2j+1}_\pi(\mathcal{C}^{s*} \otimes \mathbb{R}^n))$ is a standard odd function of degree $2u+1$ if there is a sequence of functions $F_r \in \mathcal{C}^s_n(\mathbb{R}^d, \bigoplus_{j=0}^k \Lambda^{2j+1}(\mathcal{C}^{s*} \otimes \mathbb{R}^n))$ (where $\Lambda^{2j+1}(\mathcal{C}^{s*} \otimes \mathbb{R}^n) = \underbrace{(\mathcal{C}^{s*} \otimes \mathbb{R}^n) \otimes \cdots \otimes (\mathcal{C}^{s*} \otimes \mathbb{R}^n)}_{2j+1 \text{ times}}$ and $\otimes$ is the algebraic tensor product) such that for any $g \in V$ we have $F_r(g) \to F(g)$, as $r \to +\infty$, in $\bigoplus_{j=0}^k \Lambda^{2j+1}_\pi(\mathcal{C}^{s*} \otimes \mathbb{R}^n)$ equipped with its standard topology.*

**Remark 48.** When $s > 0$, an example of a function $F = (F^1, ..., F^n)$ satisfying Definition 47 is

$$F^k(X)(x) := F(X)(\delta_x \otimes e_k) = X^k(x) \sum_{k_1, k_2 = 1}^n \int_{\mathbb{R}^d} \mathcal{G}_{k_1, k_2}(x-y) X^{k_1}(y) X^{k_2}(y) \mathrm{d}y,$$

where $\mathcal{G}_{k_1, k_2} \in \mathcal{C}^{s+\epsilon}$ are functions that decrease faster than any polynomial at infinity. We note that a finite dimensional approximation of $F^k(X)$ is given by

$$F^k_r(X)(x) = \sum_{k_1, k_2 = 1}^n \sum_{D_i \in \mathcal{P}_r} X^k(x) X^{k_1}(x - x_i) X^{k_2}(x - x_j) \int_{D_i} \mathcal{G}_{k_1, k_2}(y) \mathrm{d}y,$$

where $(\mathcal{P}_r)_{r \in \mathbb{N}}$ is sequence of increasing partitions of $\mathbb{R}^d$ for which

$$\max_{k_1, k_2 = 1,...,n} \sup_{D_i \in \mathcal{P}_r} \left| \int_{D_i} \mathcal{G}_{k_1, k_2}(y) \mathrm{d}y \right| \to 0$$

as $r \to \infty$.

**Theorem 49.** *Let $s \in \mathbb{R}$, let $F \in \mathcal{C}_n^s(\mathbb{R}^d, \bigoplus_{j=0}^u \Lambda_\pi^{2j+1}(\mathcal{C}^{s*} \otimes \mathbb{R}^n))$ be a standard odd function of degree $2u+1$ (see Definition 47), and consider $X \in \mathcal{G}^s(\mathcal{S}(\mathbb{R}^d))$. Then $F(X) \in \mathcal{C}_n^s(\mathbb{R}^d, \mathcal{A})$, furthermore for any $v_1, v_2 \in \mathcal{C}^{s*} \otimes \mathbb{R}^n$ we have*

$$\begin{aligned} F(X)(v_1)F(X)(v_2) &= -F(X)(v_2)F(X)(v_1) \\ F(X)(v_1)X(v_2) &= -X(v_2)F(X)(v_1). \end{aligned} \tag{52}$$

**Proof.** The fact that $F(X) \in \mathcal{C}_n^s(\mathbb{R}^d, \mathcal{A})$ follows from Lemma 46. The relations (52) are obvious when $F \in \mathcal{C}_n^s(\mathbb{R}^d, \bigoplus_{j=0}^k \Lambda^{2j+1}(\mathcal{C}^{s*} \otimes \mathbb{R}^n))$. Using a finite dimensional approximation $F_r$ converging to $F$, required by Definition 47, the relations (52) are extended to the more general case. $\square$

**Remark 50.** An important consequence of Theorem 49 is that if $X \in \mathcal{G}^s(V)$ and $F$ satisfies the hypotheses of Theorem 49, $F(X)$ can be extended to a Grassmann random variable (which in the following we denote by the same symbol $F(X)$) in $\mathcal{G}^s(V)$. Furthermore, $F(X)$ and $X$ are compatible (in the sense of the corresponding definition in Section 2).

We are now in the position to give a precise meaning to equation (47) and to the various hypotheses on its data. Let the operator $A$ be the (generally unbounded) generator of a strongly continuous and exponentially contractive semigroup $(e^{At})_{t \geq 0}$ on $L^2(\mathbb{R}^d)$. Moreover consider the stationary Ornstein-Uhlenbeck motion

$$B_t^A(v) = \Xi(\mathbb{I}_{(-\infty,t]} \otimes e^{A(t-\cdot)}v), \qquad t \in \mathbb{R}, v \in \mathcal{S}(\mathbb{R}^d) \otimes \mathbb{R}^n, \tag{53}$$

where $\Xi$ is a Grassmann Gaussian noise on $\mathcal{G}(C^\infty(\mathbb{R}_+) \otimes \mathcal{S}(\mathbb{R}^d) \otimes \mathbb{R}^n)$ with covariance $\mathbb{I} \otimes C$ (we consider on $C_c^\infty(\mathbb{R}_+) \otimes \mathcal{S}(\mathbb{R}^d) \otimes \mathbb{R}^n$ the pre-Hilbert space structure given by the natural inclusion of $C_c^\infty(\mathbb{R}_+) \otimes \mathcal{S}(\mathbb{R}^d) \otimes \mathbb{R}^n$ in $L^2(\mathbb{R}_+ \times \mathbb{R}^d) \otimes \mathbb{R}^n$).

**Definition 51.** *We say that the process $\Psi \in C^0([0,T], \mathcal{G}^s(V))$ satisfies equation (47) with initial condition $\Psi_0$ if $F$ is a standard odd function of finite degree, $\Psi_0 \in \mathcal{G}^s(V)$ is independent of $\Xi$, and for any $t \in [0,T]$ and any $v \in V$ we have*

$$\Psi_t(v) = \Psi_0(e^{At}v) + \lambda \int_0^t F(\Psi_s)(e^{A(t-s)}v)\mathrm{d}s + B_t^A(v) - B_0^A(e^{At}v), \tag{54}$$

*where the integral in the variable $s$ is taken in Bochner sense with respect the norm of $\mathcal{A}$.*

Hereafter, we introduce the following notations if $F \in \mathcal{C}_n^s(\mathbb{R}^d, \bigoplus_{j=0}^u \Lambda_\pi^{2j+1}(\mathcal{C}^{s*} \otimes \mathbb{R}^n))$:

$$\|F\|_{\mathcal{C}_n^s,\pi} := \sup_{k=1,\ldots,n} \|F^k\|_{\mathcal{C}^s(\mathbb{R}^d, \bigoplus_{j=0}^u \Lambda_\pi^{2j+1}\mathcal{C}^{s*})},$$

$$\deg(F) := 2u+1.$$

**Theorem 52.** *Suppose that $F \in \mathcal{C}_n^s(\mathbb{R}^d, \bigoplus_{j=0}^u \Lambda_\pi^{2j+1}(\mathcal{C}^{s*} \otimes \mathbb{R}^n))$ and that $F$ satisfies Definition 47, then for any $X, Y \in \mathcal{G}^s(V)$ compatible we have*

$$\|F(X)\|_{\mathcal{G}^s(V)} \leq \|F\|_{\mathcal{C}_n^s,\pi}(1 + \|X\|_{\mathcal{G}^s(V)})^{\deg(F)} \tag{55}$$

$$\begin{aligned} \|F(X) - F(Y)\|_{\mathcal{G}^s(V)} &\leq \deg(F)\|F\|_{\mathcal{C}_n^s,\pi}\|X-Y\|_{\mathcal{G}^s(V)} \times \\ &\quad \times (1 + \max(\|X\|_{\mathcal{G}^s(V)}, \|Y\|_{\mathcal{G}^s(V)}))^{\deg(F)-1}. \end{aligned} \tag{56}$$

**Proof.** Let $F_r^k \in \mathcal{C}^s(\mathbb{R}^d, \bigoplus_{j=0}^{u} \Lambda^{2j+1}\mathcal{C}^{s*})$ (where $\bigoplus_{j=0}^{u} \Lambda^{2j+1}\mathcal{C}^{s*}$ is the usual algebraic antisymmetric tensor product) be the approximation of $F^k$ required in Definition 47, then using the same reasoning of Theorem 12 and Theorem 14, if $K_i$ is the kernel corresponding to the Littlewood–Paley block $\Delta_i$, we obtain

$$\sup_{k=1,\ldots,n} \|F_r^k(X)(K_i(x-\cdot))\|_{\mathcal{A}} \leqslant \sup_{k=1,\ldots,n} \|F_r^k(K_i(x-\cdot))\|_{\Lambda_\pi \mathcal{C}^{s*}}(1+\|X\|_{\mathcal{G}^s(V)})^{\deg(F)}. \tag{57}$$

If we multiply both sides of (57) by $2^{is}$, take the limit as $r \to +\infty$ and afterwords the supremum on $x$ and $j$ we get inequality (55). In a similar way we get (56) by generalizing the proof of Theorem 12 and Theorem 14. □

**Remark 53.** Hereafter we use the notations, for $x \geqslant 0$,

$$f_F(x) = \|F\|_{\mathcal{C}_{n,\pi}^s}(1+x)^{\deg(F)}, \tag{58}$$

$$g_F(x) = \deg(F)\|F\|_{\mathcal{C}_{n,\pi}^s}(1+x)^{\deg(F)-1}. \tag{59}$$

We introduce some conditions on $A$ and $\Xi$ in order to prove a result on existence and uniqueness of solutions to (54).

**Definition 54.** *We say that the unbounded operator $A$ on $L^2(\mathbb{R}^d)$, and the Grassmann Gaussian noise $\Xi$ with covariance $C$ are adapted to the space $\mathcal{G}^s(V)$, for a given $s \in \mathbb{R}$, if*

  i. *The operator $A$ generates a contraction semigroup on $\mathcal{C}^s(\mathbb{R}^d)$ such that $\|e^{At}\|_{\mathcal{L}(\mathcal{C}^s(\mathbb{R}^d), \mathcal{C}^s(\mathbb{R}^d))} \lesssim e^{-\lambda_A t}$, $t \geqslant 0$, for some strictly positive constant $\lambda_A > 0$.*
  ii. *The Gaussian process $B^A$ defined by (53) belongs to $C(\mathbb{R}, \mathcal{G}^s(V))$ and furthermore we have $\sup_{t \in \mathbb{R}} \|B_t^A\|_{\mathcal{C}^s(\mathbb{R}^d, \mathcal{A})} < +\infty$.*

**Theorem 55.** *Suppose that $F$ satisfies Definition 47 and $A, \Xi$ satisfy Definition 54 and let $\Psi_0$ be a variable in $\mathcal{G}^s(V)$ independent of the Grassmann Gaussian noise $\Xi$. Then there exists a $T > 0$ such that there is a unique solution $\Psi \in C([0,T], \mathcal{G}^s(V))$ to equation (54). Furthermore $\Psi_t$ is compatible with $(\Psi_0, \Xi, \Psi_s)$ for every $t, s \in [0,T]$.*

**Proof.** The proof is a straightforward generalization of that of Theorem 30 for the case of our infinite dimensional setting. □

In order to obtain solutions for all times it is not possible to extend directly Theorem 31, since its proof was essentially based on the hypothesis that $V$ is finite dimensional. On the other hand the proofs of Theorem 37 and Lemma 39 do not depend on the dimension. For this reason we introduce the space $\mathfrak{K}_K \subset C(\mathbb{R}, \mathcal{G}^s(V))$ defined as the set of functions $\Psi_t \in C(\mathbb{R}, \mathcal{G}^s(V))$ such that

$$\sup_{t \in \mathbb{R}} \|\Psi_t\|_{\mathcal{C}_n^s(\mathbb{R}^d, \mathcal{A})} \leqslant K,$$

for some constant $K > 0$, the function $L$ given by

$$L(A, C) = \sup_{t \in \mathbb{R}} \|B_t^A\|_{\mathcal{C}_n^s(\mathbb{R}^d, \mathcal{A})},$$

and the equation for stationary solutions corresponding to (54) as

$$\Psi_t^{\mathfrak{s}} = \lambda \int_{-\infty}^{t} F(\Psi_\tau^{\mathfrak{s}}(e^{A(t-\tau)}))\mathrm{d}\tau + B_t^A \qquad t \in \mathbb{R}. \tag{60}$$

We also denote by $\Psi_{-T,t}^X$ the Grassmann process defined by $\Psi_{-T,t}^X = X$ if $t \leqslant -T$, $T \geqslant 0$, and such that $\Psi_{-T,t}^X$ is the solution to equation (54) with initial condition $X$ at time $t = -T$.

**Theorem 56.** *Suppose that $F$ satisfies Definition 47 and $A, \Xi$ satisfy Definition 54 and let $K \geqslant 4L(A,C)$. Then there exists $\lambda_0 = \lambda_0(K,L,F,\lambda_A) > 0$ such that for any $|\lambda| \leqslant 2\lambda_0$ there exists a unique (stationary) solution to equation (60) in $\mathfrak{K}_K$. Furthermore for any $|\lambda| \leqslant \lambda_0$ and for any random variable $X \in \mathcal{G}^s(V)$, independent of $B_t^A$ and such that $\|X\|_{\mathcal{C}_n^s(\mathbb{R}^d,\mathcal{A})} < K/8$, we have that $\Psi_{-T,t}^X$ converges, as $T \to +\infty$, in $\mathcal{G}^s(V)$ to $\Psi_t^{\mathfrak{s}}$ uniformly on compact subsets of $\mathbb{R}$.*

**Proof.** Due to our appropriate infinite dimensional setting the proof is a straightforward generalization of Theorem 37 and Lemma 39. □

**Remark 57.** The quantity $\lambda_0$ defined in Theorem 56 has to satisfy the following inequalities

$$\frac{2\lambda_0 f_F(K)}{\lambda_A} + 2L \leqslant K, \qquad \frac{2\lambda_0}{\lambda_A} g_F(K) < 1, \tag{61}$$

where $f_F$, $g_F$ and $\lambda_A$ are the quantities introduced in Remark 53 and Definition 54 respectively.

## 5 The Yukawa₂ model

We are now in position to study the stochastic quantization for an Euclidean fermionic QFT. We consider, as an example, the Euclidean counterpart of the Yukawa model describing the interaction of a scalar field with a spin $1/2$ field in a $1+1$ Minkowski space-time ("Yukawa₂ model") [130, 45, 21, 113, 172]. We will not attempt here to remove the UV cutoff and therefore the dimension of the space will not play any fundamental role. The Yukawa₂ model contains already the main features of a large class of Euclidean fermionic theories. It will be clear from the considerations in this section that similar results for other fermionic models are easily derivable as long as we keep an ultraviolet cut-off and we limit ourselves to a perturbative and massive regime.

### 5.1 Definition of the model

We define the Euclidean Yukawa₂ quantum field theory following Osterwalder and Schrader [131]. The Euclidean free Fermions, introduced by Osterwalder and Schrader (compare with [131], page 282 and equation (3.3), page 284), are given by a quadruplet $\tilde{\chi} = (\chi, \bar{\chi}) \in \mathcal{G}(\mathcal{S}(\mathbb{R}^2) \otimes \mathbb{C}^4)$ of complex Grassmann Gaussian fields. Here we employ the convention of denoting by $\chi, \bar{\chi} \in \mathcal{G}(\mathcal{S}(\mathbb{R}^2) \otimes \mathbb{C}^2)$ two independent fields (in particular the bar over $\bar{\chi}$ does not denote complex conjugation). By considering two fields $\chi, \bar{\chi}$ we are, "doubling the dimension". This "dimension doubling", also employed by Osterwalder and Schrader, has the purpose of obtaining an Euclidean covariance with the right symmetry properties.

We adapt to the present case some notation for complex Grassmann Gaussian random variables which we introduced in Section 2.3. We consider $V = \mathcal{S}(\mathbb{R}^2) \otimes \mathbb{C}^4 \subseteq L^2(\mathbb{R}^2) \otimes \mathbb{C}^4$ as a complex pre-Hilbert space with an Hermitian scalar product denoted $(\cdot,\cdot)_V$ equipped with the standard real structure $c$ given by the complex conjugation (i.e. $cf = \bar{f}$) determining the bilinear form $\langle f,g\rangle_c := (cf,g)_V$. Given an operator $A$ on $V$, we denote by $A^c$ the complex-conjugated Hilbert adjoint, i.e. $A^c = c A^* c$, $A^*$ being the Hilbert adjoint of $A$. Finally, we say that an operator $K$ is $c$-antisymmetric when $K^c = -K$. (recall Definition 19 and Definition 20).

With these notations, the fields $\tilde{\chi}$ are defined to be complex Grassmann Gaussian variables with covariance

$$\omega(\tilde{\chi}(f)\tilde{\chi}(g)) = \langle f, Kg\rangle_c, \qquad f,g \in \mathcal{S}(\mathbb{R}^2) \otimes \mathbb{C}^4. \tag{62}$$

where $K\colon L^2(\mathbb{R}^2)\otimes \mathbb{C}^4 \to L^2(\mathbb{R}^2)\otimes \mathbb{C}^4$ is given by the following $c$-antisymmetric bounded operator

$$K = (m_f^2 - \Delta)^{-1} C_f, \qquad C_f := \begin{pmatrix} 0 & (\slashed{\nabla} + m_f) \\ -(\slashed{\nabla} + m_f)^c & 0 \end{pmatrix}. \tag{63}$$

Here $\Delta$ is the Laplacian on $\mathbb{R}^2$, $m_f > 0$ is the fermion mass and $\slashed{\nabla}$ is the Euclidean Dirac operator in 2 Euclidean space-time dimensions. Explicitly, $\slashed{\nabla} := \gamma_1 \nabla_1 + \gamma_2 \nabla_2$ with $\nabla_i$ the partial derivatives with respect to $x_i$, $i = 1, 2$, and with $\gamma_1, \gamma_2 \colon \mathbb{C}^2 \to \mathbb{C}^2$ the Euclidean Dirac matrices, i.e. Hermitian matrices satisfying

$$\gamma_a \gamma_b + \gamma_b \gamma_a = 2\delta_{ab}, \qquad a, b = 1, 2. \tag{64}$$

Concretely, we can take $\gamma_1, \gamma_2$ to be simply the first two Pauli matrices:

$$\gamma_1 = \begin{pmatrix} 0 & 1 \\ 1 & 0 \end{pmatrix}, \quad \gamma_2 = \begin{pmatrix} 0 & -\mathrm{i} \\ \mathrm{i} & 0 \end{pmatrix}.$$

Note that $K = C_f^{-1}$ and that

$$\tilde{C}_f = (m_f^2 - \Delta)^{1/2} K = (m_f^2 - \Delta)^{-1/2} C_f$$

is also bounded and antisymmetric in $L^2(\mathbb{R}^2) \otimes \mathbb{C}^4$, hence the covariance matrix of the field $\tilde{\psi}$ can also be written as

$$\omega(\tilde{\chi}(f)\tilde{\chi}(g)) = \langle (m_f^2 - \Delta)^{-1/4} f, \tilde{C}_f (m_f^2 - \Delta)^{-1/4} g \rangle_c.$$

Therefore we can realize $\tilde{\chi}$ as

$$\tilde{\chi}(f) := \tilde{\zeta}((m_f^2 - \Delta)^{-1/4} f), \qquad f \in \mathcal{S}(\mathbb{R}^2) \otimes \mathbb{C}^4,$$

where $\tilde{\zeta} \in \mathcal{G}(\mathcal{S}(\mathbb{R}^2) \otimes \mathbb{C}^4)$ is the complex Grassmann Gaussian field specified by the bounded covariance

$$\omega(\tilde{\zeta}(f)\tilde{\zeta}(g)) = \langle f, \tilde{C}_f g \rangle_c, \qquad f, g \in \mathcal{S}(\mathbb{R}^2) \otimes \mathbb{C}^4.$$

As a consequence, in this specific realization, we have

$$\|\tilde{\chi}(f)\| = \|\tilde{\zeta}((m_f^2 - \Delta)^{-1/4} f)\| \lesssim \|(m_f^2 - \Delta)^{-1/4} f\|_{L^2(\mathbb{R}^2) \otimes \mathbb{C}^4} < \infty, \tag{65}$$

where the implicit constant in the inequality depends only on $m_f$. This realization gives better spatial regularity properties to the Euclidean fields $\tilde{\chi}$.

In order to describe the Yukawa$_2$ model we need also to introduce an Euclidean boson field $\varphi$ of mass $m_b > 0$: it is the centered Gaussian field $\varphi$ on $\mathbb{R}^2$ with covariance $(m_b^2 - \Delta)^{-1}$ where $m_b > 0$ is looked upon as a fixed parameter and $\Delta$ is the Laplacian on $\mathbb{R}^2$. We can realize $\varphi$ as the canonical multiplication operator on the Hilbert space $\mathcal{H}_\varphi = L^2(\mathcal{C}_{\mathrm{loc}}^{-\kappa}(\mathbb{R}^2), \nu)$ where $\nu$ is the centered Gaussian measure on $\mathcal{C}_{\mathrm{loc}}^{-\kappa}(\mathbb{R}^2)$ (the local version of the $\mathcal{C}^{-\kappa}$ Besov space) with covariance $(m_b^2 - \Delta)^{-1}$ and $\kappa > 0$ is any small positive quantity. We shall denote by $\omega_\varphi$ the state associated with $\nu$ and $\mathcal{C}_{\mathrm{loc}}^{-\kappa}(\mathbb{R}^2)$.

In addition, large scale and small scale cutoffs are introduced as follows. Fix a smooth positive functions $h \in C^\infty(\mathbb{R}^2)$ with compact support and let $\mathfrak{a}$ be a regular enough (see below for the precise regularity) continuous function with compact support such that $\mathfrak{a}(0) = 1$. Take $\varepsilon > 0$ and let $\mathfrak{a}_\varepsilon(x) = \mathfrak{a}(x/\varepsilon)\varepsilon^{-2}$. Then we define new (regularized) Grassmann variables

$$\tilde{\psi}(x) := \tilde{\chi}(\mathfrak{a}_\varepsilon(\cdot - x)), \qquad x \in \mathbb{R}^2.$$

Then $\tilde{\psi} \in \mathcal{G}^s(\mathcal{S}(\mathbb{R}^2) \otimes \mathbb{C}^4)$ for some $s > 0$ (in general depending on the regularity of $\mathfrak{a}$ see below).

Following the argument of Section 6 of [131] the Euclidean Yukawa$_2$ model is defined as follows.

**Definition 58.** *The approximate fermionic Schwinger functions $(\rho^n_{\varepsilon,h})_n$ of the Yukawa$_2$ model on $\mathbb{R}^2$ with infrared (IR) cutoff $h$ and ultraviolet (UV) fermion cutoff $\varepsilon$ are given by the Euclidean averages*

$$\rho^n_{\varepsilon,h}(f_1,...,f_n) = \frac{\omega \otimes \omega_\varphi\big(\tilde{\psi}(f_1)\cdots\tilde{\psi}(f_n)e^{-V_h(\varphi,\tilde{\psi})}\big)}{\omega \otimes \omega_\varphi\big(e^{-V_h(\varphi,\tilde{\psi})}\big)} \qquad (66)$$

*with potential $V_h(\varphi,\tilde{\psi})$ given by the unbounded (non-self-adjoint) operator on $\mathcal{H} \otimes \mathcal{H}_\varphi$, where $\mathcal{H}$ is the fermionic Hilbert space,*

$$V_h(\varphi,\tilde{\psi}) = \lambda \int_{\mathbb{R}^2} \mathrm{d}x\, h(x)(\bar{\psi}(x)\cdot\psi(x))\varphi(x), \qquad (67)$$

$\lambda \in \mathbb{R}$ *is a coupling constant and $\omega$, resp. $\omega_\varphi$ are the fermion, resp. boson, states described above.*

Note that there is no difficulty in making sense of (66), due to the following observations (see also [131]). The map $x \in \mathbb{R}^2 \mapsto \tilde{\psi}_\varepsilon(x) \in \mathcal{A} = \mathcal{L}(\mathcal{H})$ is $\mathcal{C}^s$ due to the presence of the regularization $\varepsilon > 0$ and the regularity of the Euclidean field $\tilde{\psi}$ given in (65). The massive Gaussian free field $\varphi$ has (local) regularity $\mathcal{C}^{-\kappa}_{\mathrm{loc}}$ for arbitrarily small $\kappa > 0$, therefore, taking $s > \kappa$, by a standard Besov estimates for Banach valued functions (see e.g. [9, 10]), we have, for any *fixed* $\varphi \in \mathcal{C}^{-\kappa}_{\mathrm{loc}}(\mathbb{R}^2)$

$$\|V_h(\varphi,\tilde{\psi})\|_\mathcal{A} \lesssim \lambda \|\tilde{\psi}\|^2_{\mathcal{C}^\kappa(\mathcal{A})} \|h\varphi\|_{B^{-\kappa}_{1,1}(\mathbb{R}^2)} \lesssim \lambda \|h\varphi\|_{\mathcal{C}^{-\kappa}(\mathbb{R}^2)}.$$

It is also not difficult to see that the map $V_h: \mathcal{C}^{-\kappa} \to \mathcal{A}$ given by $V_h: \varphi \mapsto V_h(\varphi,\tilde{\psi})$ is linear and continuous in $\varphi$, therefore measurable and that

$$\omega \otimes \omega_\varphi\big(\tilde{\psi}(f_1)\cdots\tilde{\psi}(f_n)e^{-V_h(\varphi,\tilde{\psi})}\big) = \omega\big(\tilde{\psi}(f_1)\cdots\tilde{\psi}(f_n)\omega_\varphi\big(e^{-V_h(\varphi,\tilde{\psi})}\big)\big)$$

e.g. by proving this first for a piecewise approximation of $\varphi$ and then taking the limit. Moreover we have (by the properties of the Bochner integral)

$$\big\|\omega_\varphi\big(e^{-V_h(\varphi,\tilde{\psi})}\big)\big\|_\mathcal{A} \leqslant \omega_\varphi(\exp(C|\lambda|\|h\varphi\|_{\mathcal{C}^{-\kappa}})) < \infty,$$

for all $\lambda \in \mathbb{R}$. We conclude that the averages in (66) are well defined since

$$\big|(\omega \otimes \omega_\varphi)\big(\tilde{\psi}(f_1)\cdots\tilde{\psi}(f_n)e^{-V_h(\varphi,\tilde{\psi})}\big)\big| \leqslant \big\|\tilde{\psi}(f_1)\cdots\tilde{\psi}(f_n)\omega_\varphi\big(e^{-V_h(\varphi,\tilde{\psi})}\big)\big\|_\mathcal{A} < \infty.$$

The strict positivity $Z = \omega \otimes \omega_\varphi\big(e^{-V_h(\varphi,\tilde{\psi})}\big) > 0$ has been proved in [131] via the Feynman–Kac formula and reflection positivity for cutoffs $h$ which are symmetric with respect to. reflection in one spatial coordinate. In [163] the positivity of $Z$ is proved in Lemma 4.1 there using the Matthew–Salam formula. Alternatively, provided $\lambda$ is small enough one can use Remark 41 to assert that $Z \neq 0$.

By expanding the quantity $\omega_\varphi\big(e^{-V_h(\varphi,\tilde{\psi})}\big) \in \mathcal{A}$ in power series and integrating away the Gaussian field we obtain

$$\omega_\varphi\big(e^{-V_h(\varphi,\tilde{\psi})}\big) = e^{\mathcal{V}_h(\tilde{\psi})},$$

with a purely fermionic potential given by

$$\mathcal{V}_h(\tilde{\psi}) = \lambda^2 \int_{(\mathbb{R}^2)^2} (h\,\bar{\psi}\cdot\psi)(x)\mathcal{G}(x-y)(h\,\bar{\psi}\cdot\psi)(y)\mathrm{d}x\mathrm{d}y. \tag{68}$$

Here $\mathcal{G} = (m_b^2 - \Delta)^{-1}$ is the Green's function of the operator $(m_b^2 - \Delta)$ and, as before, $m_b$ is the boson mass. This gives us the purely fermionic representation of the Schwinger functions (66) given by:

$$\rho_{\varepsilon,h}^n(f_1,...,f_n) := \omega_{rY}^h(\tilde{\psi}(f_1)\cdots\tilde{\psi}(f_n)) = \frac{\omega\bigl(\tilde{\psi}(f_1)\cdots\tilde{\psi}(f_n)e^{\mathcal{V}_h(\tilde{\psi})}\bigr)}{\omega\bigl(e^{\mathcal{V}_h(\tilde{\psi})}\bigr)}. \tag{69}$$

Perturbation theory suggests that the boson field requires a mass renormalization $\mathfrak{a}_\varepsilon$ which depends on the fermion cutoff $\varepsilon > 0$. In theory one would need to modify the bosonic Gaussian field covariance to include a mass renormalization and let $\varphi$ have covariance $(m_b^2 - \Delta + a_\varepsilon)^{-1}$ instead of the original $(m_b^2 - \Delta)^{-1}$. However in this paper we will not discuss the removal of the UV cutoff and keep $\varepsilon > 0$ fixed. Therefore, for simplicity, we ignore this mass renormalization.

## 5.2 Stochastic quantization

In order to stochastically quantize the Yukawa model, we consider the Grassmannian differential equation

$$\mathrm{d}\tilde{\Psi}_t^h = ((\Delta - m_f^2)\tilde{\Psi}_t^h + \lambda^2 F_{h,rY}(\tilde{\Psi}_t^h))\mathrm{d}t + \mathrm{d}\tilde{B}_{\varepsilon,t}, \quad \tilde{\Psi}_t^h := (\Psi_t^h, \bar{\Psi}_t^h) \colon \mathbb{R}^2 \to \mathcal{A} \otimes \mathbb{C}^4 \tag{70}$$

where $\Delta$ is the Laplacian on $\mathbb{R}^2$, the function $F_{h,rY}$ is defined as

$$F_{h,rY}(\tilde{\Psi}_t^h) = C_f\,\mathfrak{a}_\varepsilon^{*2} * \left[\frac{h}{2}\left(\int_{\mathbb{R}^2} \mathcal{G}(y-\cdot)h(y)\Psi_t^h(y)\bar{\Psi}_t^h(y)\mathrm{d}y\right)J\tilde{\Psi}_t^h\right],$$

with $\mathfrak{a}_\varepsilon(x) = \varepsilon^{-2}\mathfrak{a}\bigl(\frac{x}{\varepsilon}\bigr)$, $\varepsilon > 0$, $\mathfrak{a}\colon \mathbb{R}^2 \to \mathbb{R}^2$ is regular enough function with compact support (see below for the more precise requirements), $\mathfrak{a}_\varepsilon^{*2} = \mathfrak{a}_\varepsilon * \mathfrak{a}_\varepsilon$, $h\colon \mathbb{R}^2 \to \mathbb{R}$ is a smooth function (in the beginning with compact support, later on we will consider $h \equiv 1$) and $\mathcal{G}$ is the Green's function of the elliptic operator $(m_b^2 - \Delta)$ and $J\colon \mathbb{C}^4 \to \mathbb{C}^4$ is given by

$$J = \begin{pmatrix} 0 & 1 \\ -1 & 0 \end{pmatrix}.$$

Finally, $\tilde{B}_{\varepsilon,t}$ is a Brownian motion given by

$$\tilde{B}_{\varepsilon,t}(v) = \tilde{\Xi}(\mathbb{I}_{(-\infty,t]} \otimes (\mathfrak{a}_\varepsilon * v)), \qquad t \in \mathbb{R}, v \in \mathcal{S}(\mathbb{R}^2) \otimes \mathbb{C}^4,$$

where $\tilde{\Xi}$ is a white noise on the pre-Hilbert space $L^2(\mathbb{R}_+) \otimes \mathcal{S}(\mathbb{R}^2) \otimes \mathbb{C}^4 \subseteq L^2(\mathbb{R}_+) \otimes B_{2,2}^{-1/2}(\mathbb{R}^2) \otimes \mathbb{C}^4$ equipped with the scalar product

$$((f,g)) = (f, (m_f^2 - \Delta)^{-1/2}g)_{L^2(\mathbb{R}_+ \times \mathbb{R}^2) \otimes \mathbb{C}^4} \qquad f, g \in L^2(\mathbb{R}_+) \otimes \mathcal{S}(\mathbb{R}^2) \otimes \mathbb{C}^4, \tag{71}$$

and having the bounded covariance $\tilde{C}_f = (m_f^2 - \Delta)^{-1/2}C_f$.

We claim that the stationary solution to the Grassmann SPDEs (70) is associated to the Yukawa measure $\omega_{rY}^h(\cdot) = \omega(\cdot e^{\mathcal{V}_h})$. In particular it is the solution to the integral equation

$$\tilde{\Psi}_t^h(v) = \lambda^2 \int_{-\infty}^t F_{h,rY}(\tilde{\Psi}_\tau^h)(e^{A(t-\tau)}v)\mathrm{d}\tau + \tilde{B}_{\varepsilon,t}^A(v), \qquad t \in \mathbb{R}, v \in \mathbb{C}^4, \tag{72}$$

where $A = \Delta - m_f^2$ and the regularized noise $\tilde{B}_{\varepsilon,t}^A$ is given by the expression

$$\tilde{B}_{\varepsilon,t}^A(v) = \tilde{\Xi}(\mathbb{I}_{(-\infty,t]} \otimes e^{A(t-\cdot)}(\mathfrak{a}_\varepsilon * v)), \qquad t \in \mathbb{R}, v \in \mathcal{S}(\mathbb{R}^2) \otimes \mathbb{C}^4, \tag{73}$$

where $\tilde{\Xi}$ is the white noise introduced above.

**Remark 59.** From now on we will make the following assumption. We fix $\delta > 0$ and $\gamma > 0$, let $\alpha = \gamma + 1/2 + 2\delta$ and take $\mathfrak{a} \in B_{1,\infty}^\alpha(\mathbb{R}^2) \cap \mathcal{C}^\alpha(\mathbb{R}^2)$. Again we recall that $\varepsilon > 0$ will be kept fixed.

The results which we will prove is the following:

**Theorem 60.** *With the above notation and conditions on $F_{h,rY}$ consider $\alpha, \gamma, \delta$ and $\mathfrak{a}$ as in Remark 59. Then for $0 \leqslant |\lambda| < \lambda_{\varepsilon,h}$ small enough equation (70) admits a unique stationary solution $\tilde{\Psi}^h$ in $\mathcal{C}^\gamma(\mathbb{R}^d, \mathcal{A}) \otimes \mathbb{C}^4$ given by (72). Furthermore we have that, for any $k \geqslant 1$, $f_1, ..., f_k \in C_c^\infty(\mathbb{R}^2) \otimes \mathbb{C}^4$ and $t \in \mathbb{R}$,*

$$\omega_{rY}^h(\tilde{\psi}(f_1)\cdots\tilde{\psi}(f_k)) = \omega(\tilde{\Psi}_t^h(f_1)\cdots\tilde{\Psi}_t^h(f_k)) \tag{74}$$

*where $(\rho_{rY}^h, \tilde{\psi})$ is a regularized Yukawa$_2$ model with spatial cut-off $h$ given by (69).*

*Furthermore if $0 \leqslant |\lambda| < \lambda_{\varepsilon,1}$ is small enough, and if we denote by $\tilde{\Psi}_t$ the stationary solution to equation (70) with $h \equiv 1$, we have that for any $f_1, ..., f_k \in C_c^\infty(\mathbb{R}^2) \otimes \mathbb{C}^4$ and for any $t \in \mathbb{R}$*

$$\lim_{h_n \to 1} \omega_{rY}^{h_n}(\tilde{\psi}(f_1)\cdots\tilde{\psi}(f_k)) = \omega(\tilde{\Psi}_t(f_1)\cdots\tilde{\Psi}_t(f_k)). \tag{75}$$

*where $(h_n)_{n \geqslant 1}$ is any sequence of smooth functions with compact support converging to 1 uniformly in $\mathcal{C}^{\gamma+\delta}$ on compact subsets of $\mathbb{R}^2$.*

**Proof.** We outline the idea of the proof while the details will constitute the rest of this section. The existence and uniqueness of stationary solution to equation (70) will be obtained in Lemma 64 thanks to the general theory of infinite dimensional Grassmannian SDE described in Section 4. The main difficulty is to prove equality (74) and the convergence of the limit (75). We proved the explicit "law" of the stationary solutions to Grassmannian SDE only for the finite dimensional case (see Theorem 40). For this reason we have to approximate the solutions $\tilde{\Psi}_t^h$ with the stationary solution of finite dimensional equations. We can apply Theorem 40, and, then, obtain (74) and (75) by removing the various approximations and proving the convergence of the finite dimensional equation to the infinite dimensional one. This will be the subject of Lemma 70 and Lemma 73 which together with equation (79) will conclude the proof. □

**Remark 61.** At a heuristic level the stochastic quantization equation is not unique and it could be of the form

$$d\tilde{\Psi}_t^h = L\left(C_f \tilde{\Psi}_t^h + \frac{\delta V^h}{\delta \tilde{\Psi}}(\tilde{\Psi}_t^h)\right)dt + d\tilde{B}_t^L$$

where $V^h$ is the (fermionic effective) potential, $L$ is some linear operator and $\tilde{B}$ is a suitable Gaussian noise whose covariance depends on $L$. In order to apply the theory of Section 3.5 and Section 4, and so obtain the invariant measure as the "law" of the solution to the stationary equation, we need that the operator $LC_f$ is (strictly) positive definite. This is the main reason of our choice $L = C_f$.

As we discussed, we start by proving that, under the hypotheses of Theorem 60, equation (72) admits a unique solution.

**Lemma 62.** *We have $\tilde{B}_\varepsilon^A \in \mathcal{C}^\delta(\mathbb{R}, \mathcal{C}^\gamma(\mathbb{R}^2, \mathcal{A}) \otimes \mathbb{C}^4)$ and*

$$\sup_{t \in \mathbb{R}} \|\tilde{B}_{\varepsilon,t}^A\|_{\mathcal{C}^{\gamma+2\delta}(\mathbb{R}^2, \mathcal{A}) \otimes \mathbb{C}^4} < +\infty.$$

The proof of Lemma 62 is given in Section 5.4.

**Lemma 63.** *For any $s > 0$, the function $F_{h,rY}$ is well defined on $\mathcal{C}^s(\mathbb{R}^d, \mathcal{A})$ into itself, it satisfies Definition 47. Furthermore we have $\|F_{h,rY}\|_{\mathcal{C}^s, \pi} \lesssim \|h\|_{\mathcal{C}^s}^2$.*

**Proof.** First we prove that $F \in \mathcal{C}^s(\mathbb{R}^2, \Lambda_\pi^3(\mathcal{C}^{s*}(\mathbb{R}^2) \otimes \mathbb{C}^4))$. We start by noting that, if we denote by $\delta_x$ the Dirac delta with unit mass at the point $x \in \mathbb{R}^2$, we have $\delta_x \otimes e_i \in \mathcal{C}^{0*}(\mathbb{R}^2)$, for $i = 1, ..., 4$, where $e_i$ is the projection on the component $i$-th of $\mathbb{C}^4$. Furthermore, since for $s \notin \mathbb{N}$, $s > 0$, $\mathcal{C}^s(\mathbb{R}^2)$ is the space of $s$ continuous Hölder functions, we have that for any $0 < s' < s$ and $s' \notin \mathbb{N}$ the map $x \mapsto \delta_x \otimes e_i$ is in $\mathcal{C}^{s'}(\mathbb{R}^2, \mathcal{C}^{s*}(\mathbb{R}^2) \otimes \mathbb{C}^4)$. This and the fact that, by Theorem 82, $\mathcal{C}^{s'}(\mathbb{R}^2, \Lambda_\pi(\mathcal{C}^{s*}(\mathbb{R}^2) \otimes \mathbb{C}^4))$ is a Banach algebra with respect to the multiplication $\wedge$, imply that the map

$$\tilde{F}(x) := h(x)((\delta_x \otimes e_1) \wedge (\delta_x \otimes e_3) + (\delta_x \otimes e_2) \wedge (\delta_x \otimes e_4))$$

is $s'$-Hölder as a function from $\mathbb{R}^2$ into $\Lambda_\pi^2(\mathcal{C}^{s*}(\mathbb{R}^2) \otimes \mathbb{C}^4)$, with

$$\|\tilde{F}\|_{\mathcal{C}^{s'}(\mathbb{R}^2, \Lambda_\pi^2(\mathcal{C}^{s*}(\mathbb{R}^2) \otimes \mathbb{C}^4))} \lesssim \|h\|_{\mathcal{C}^s} \tag{76}$$

where the constants in the symbol $\lesssim$ do not depend on $h$. On the other hand, by Theorem 80, $\mathcal{G} * \tilde{F} \in \mathcal{C}^{s'+2}(\mathbb{R}^2, \Lambda_\pi(\mathcal{C}^{s*}(\mathbb{R}^2) \otimes \mathbb{C}^4))$ which, if $s - s' > 2$, is contained in $\mathcal{C}^s(\mathbb{R}^2, \Lambda_\pi(\mathcal{C}^{s*}(\mathbb{R}^2) \otimes \mathbb{C}^4))$. Since every component of $F_{\varepsilon,h,Y}$ is of the form $((\slashed{\nabla} + m_f)\mathfrak{a}_\varepsilon^{*2}) * ((h(x)\delta_x \otimes e_i) \wedge \mathcal{G}*\tilde{F})$, using the fact that, by Theorem 80 and Theorem 81, $(\slashed{\nabla} + m_f)\mathfrak{a}_\varepsilon^{*2} * \cdot : \mathcal{C}^s(\mathbb{R}^2, \mathcal{A}) \to \mathcal{C}^{s-1+2\alpha}(\mathbb{R}^2, \mathcal{A}) \subset \mathcal{C}^s(\mathbb{R}^2, \mathcal{A})$ (since $\alpha > \frac{1}{2}$ for the conditions in Remark 59) is a continuous operator, and by Theorem 82 we get the thesis. As byproduct of the previous reasoning we obtain also that $\|F_{h,rY}\|_{\mathcal{C}^s, \pi} \lesssim \|h\|_{\mathcal{C}^s}^2$. This can be seen using the fact that, when $s > 0$, $\tilde{\Psi}_{\varepsilon,t}$ is Hölder continuous in space and that $\mathcal{G}$ decreases exponentially at infinity.

In order to finish the proof we have to show that $F_{h,rY}$ admits an approximation of the form required by Definition 47. This can be obtained by approximating the convolution by $\mathcal{G}$ and $\mathfrak{a}_\varepsilon$ (the regularizing function in the definition of $\tilde{\Psi}_\varepsilon$) using a finite sum. Indeed, if $\mathcal{P}_n$ is a sequence of increasing (finite) partition of the $\mathbb{R}^2$, we can approximate the integral $\mathcal{G}*\tilde{F}(x)$ by the sum

$$S_n(x) = \sum_{D_i \in \mathcal{P}_n} \tilde{F}(x - x_i) \int_{D_i} \mathcal{G}(y) \mathrm{d}y, \qquad x \in \mathbb{R}^2,$$

where $x_i \in D_i$ is any point in $D_i$. If we choose the sequence of partitions $\mathcal{P}_n$ in such a way that $\lim_{n \to +\infty} \{\sup_i (\int_{D_i} \mathcal{G}(y) \mathrm{d}y)\} = 0$, we have that, as $n \to \infty$, $S_n \to \mathcal{G}*\tilde{F}$ in $\mathcal{C}^{s'}(\mathbb{R}^2, \Lambda_\pi(\mathcal{C}^{s*}(\mathbb{R}^2) \otimes \mathbb{C}^4))$, whenever $s' < s$. The convolution with $(C_f \mathfrak{a}_\varepsilon^{*2})$ can be handled in a similar way, completing the proof. □

Using the previous two lemmas we are in a position to prove the following.

**Lemma 64.** *There exist constants $K, \lambda_* > 0$ depending on $m_f, m_b$ and $\|h\|_{\mathcal{C}^s}^2$, such that if $|\lambda| \leqslant \lambda_*$ there exists a unique solution to the stationary equation (72) which is uniformly bounded in $\mathcal{C}^\gamma(\mathbb{R}^d, \mathcal{A}) \otimes \mathbb{C}^4$ with bound $K$.*

**Proof.** The proof is based on Theorem 55 and Theorem 56, whose hypotheses are satisfied in this setting because of Lemma 62 and Lemma 63. The uniformity on the constant $K$ and the bound $\lambda_*$ on the admissible $\lambda$ follows from Lemma 63, Theorem 52 and the inequalities (61). □

## 5.3 Approximation

In this subsection we introduce the finite dimensional approximations of both the Yukawa measure $\omega_{rY}$ and the stochastic differential equation in order to establish the stochastic quantization connection. These approximation consist in making the spatial dependence periodic with period $R$ and the performing a Galerkin approximation with the first $N$ Fourier modes. Removal of these approximation will be done by first passing to the $N \to \infty$ limit and then to the $R \to \infty$ limit.

Consider the torus
$$\mathbb{T}_R^2 = \mathbb{R}^2 / (2\pi R \mathbb{Z})^2$$
for $R > 0$. Hereafter we often identify $\mathbb{T}_R^2$ with the square $(-\pi R, \pi R]^2 \subset \mathbb{R}^2$. We defined the periodic version $\tilde{\psi}_R$ of the Fermion field $\tilde{\psi}$, in the following way: consider $f \in C^\infty(\mathbb{T}_R^2)$ then we write
$$\tilde{\psi}_R(f) := \tilde{\chi}(K^{-1}(\mathbb{I}_{\mathbb{T}_R^2} K_{\mathbb{T}_R^2}(\mathfrak{a}_\varepsilon *_{\mathbb{T}_R^2} f))),$$
where $*_{\mathbb{T}_R^2}$ is the convolution in $\mathbb{T}_R^2$, $\mathfrak{a}_\varepsilon$ is identified with its restriction to $(-\pi R, \pi R]^2$, the operator $K_{\mathbb{T}_R^2}: C^\infty(\mathbb{T}_R^2) \to C^\infty(\mathbb{T}_R^2)$ has the form
$$K_{\mathbb{T}_R^2} = (m_f^2 - \Delta)^{-1} C_f = (m_f^2 - \Delta)^{-1} \begin{pmatrix} 0 & (\slashed{\nabla} + m_f) \\ -(\slashed{\nabla} + m_f)^c & 0 \end{pmatrix},$$
we identify the function $\mathbb{I}_{\mathbb{T}_R^2}(\cdot) g(\cdot): \mathbb{T}_R^2 \to \mathbb{R}$ (where $g: \mathbb{T}_R^2 \to \mathbb{R}$) with the function which is equal to $g$ on $\mathbb{T}_R^2$ and 0 otherwise, and $K^{-1}: L^2(\mathbb{R}^2) \otimes \mathbb{C}^4 \to H^1(\mathbb{R}^2) \otimes \mathbb{C}^4$ is the inverse of the operator $K$ defined in equation (63). The map $x \mapsto \tilde{\psi}_R(x)$ is in the space $\mathcal{C}^{\gamma+2\delta}(\mathbb{T}_R^2, \mathcal{A})$ and it is uniformly bounded in $\mathcal{A}$ (where $\gamma, \delta \in \mathbb{R}_+$ are the constants of Remark 59). Furthermore, if we denote by
$$\tilde{\psi}_R^p(f) := \tilde{\psi}_R(Q_R(f)), \quad Q_R(f)(x) := \sum_{n \in \mathbb{Z}^2} f(x - 2\pi R n), x \in \mathbb{T}_R^2,$$
i.e. the periodic extension of $\tilde{\psi}_R$ to the whole $\mathbb{R}^2$, we have
$$\sup_{x \in \mathcal{K}} \|\tilde{\psi}_R^p(x) - \tilde{\psi}(x)\|_{\mathcal{A} \otimes \mathbb{C}^4} \to 0,$$
uniformly on compacts $\mathcal{K} \subset \mathbb{R}^2$ as $R \to \infty$. A proof of this statement can be obtained in a way similar to Lemma 72 (see Section 5.4) and is based on the fact that the correlation function of $\tilde{\psi}_R$ converges to the one of $\tilde{\psi}$. As a consequence, if we define the functional $\mathcal{V}_{R,h}: \mathcal{C}^0(\mathbb{T}_R^2, \mathcal{A} \otimes \mathbb{C}^4) \to \mathbb{R}$ as
$$\mathcal{V}_{R,h}(\tilde{\psi}_R) := \int_{\mathbb{T}_R^4} h(x) h(y) (\psi_R(x) \bar{\psi}_R(x)) \mathcal{G}_R(x-y)(\psi_R(y) \bar{\psi}_R(y))(y) \mathrm{d}x \mathrm{d}y, \tag{77}$$
where $\mathcal{G}_R: \mathbb{T}_R^2 \to \mathbb{R}$ is the function $\mathcal{G}_R(x) := \sum_{n \in \mathbb{Z}^2} \mathcal{G}(x - 2\pi R n)$, we have
$$\begin{aligned}\omega\big(\tilde{\psi}(f_1), ..., \tilde{\psi}(f_n) e^{\mathcal{V}_h(\tilde{\psi})}\big) &= \lim_{R \to \infty} \omega\big(\tilde{\psi}_R^p(f_1) \cdots \tilde{\psi}_R^p(f_n) e^{\mathcal{V}_{h_R}(\tilde{\psi}_R^p)}\big) \\ &= \lim_{R \to \infty} \omega\big(\tilde{\psi}_R(f_1) \cdots \tilde{\psi}_R(f_n) e^{\mathcal{V}_{R,h}(\tilde{\psi}_R)}\big)\end{aligned} \tag{78}$$
for all compactly supported $f_1, ..., f_n, h \in \mathcal{S}(\mathbb{R}^2)$ (where we used the fact that if the support of $f_1, ..., f_n, h$ is contained in $\mathbb{T}_R^2 = (-\pi R, \pi R]^2$, then $Q_R(f_1) = f_1, ..., Q_R(f_n) = f_n, Q_R(h) = h$), and $h_R$ is the periodic extension of $h|_{\mathbb{T}_R^2}$ defined in Remark 68 below.

Let now $(e_k)_{k \in \mathbb{Z}^2}$ be the orthonormal Fourier basis of $L^2(\mathbb{T}_R^2)$ and $P_N \colon L^2(\mathbb{T}_R^2) \to L^2(\mathbb{T}_R^2)$ be the orthogonal projection to the closed subspace generated by $\{e_k \colon k \in \mathbb{Z}^2, |k| \leqslant N\}$. Let $\tilde{\psi}_{N,R} := P_N \tilde{\psi}_R$ and observe again that $\|\tilde{\psi}_{N,R}(x) - \tilde{\psi}_R(x)\|_{\mathcal{C}^\alpha(\mathbb{T}_R^2)} \to 0$ as $N \to \infty$, for any $0 < \alpha < \gamma + 2\delta$ (where $\gamma, \delta \in \mathbb{R}_+$ are the constants of Remark 59) and therefore

$$\omega\big(\tilde{\psi}(f_1) \cdots \tilde{\psi}(f_n) e^{\mathcal{V}_h(\tilde{\psi})}\big) = \lim_{R \to \infty} \lim_{N \to \infty} \omega\big(\tilde{\psi}_{N,R}(f_1) \cdots \tilde{\psi}_{N,R}(f_n) e^{\mathcal{V}_{R,h}(\tilde{\psi}_{N,R})}\big). \tag{79}$$

Consider the finite set $\Lambda_N = \{k \in \mathbb{Z}^2 \colon |k| \leqslant N\}$ and let

$$\tilde{\theta}(k) := \int_{\mathbb{T}_R^2} \overline{e_k(x)} \tilde{\psi}_{N,R}(x) \mathrm{d}x, \qquad k \in \Lambda_N.$$

We can show easily that $(\tilde{\theta}(k))_{k \in \Lambda_N}$ is a Grassmann Gaussian vector with covariance given by

$$\omega(\tilde{\theta}_\alpha(k) \tilde{\theta}_\beta(\ell)) = \delta_{k+\ell=0} \begin{pmatrix} 0 & (i\slashed{k} + m_f)^{-1} \\ -((i\slashed{k} + m_f)^{-1})^c & 0 \end{pmatrix}_{\alpha, \beta} \tag{80}$$

where $\alpha, \beta = 1, \ldots, 4$ denote the canonical coordinates in $\mathbb{C}^4$ and $k, \ell \in \Lambda_N$. Moreover

$$\tilde{\psi}_{N,R}(x) = \sum_{k \in \Lambda_N} e_k(x) \tilde{\theta}(k), \qquad x \in \mathbb{T}_R^2$$

and $\tilde{\psi}_{N,R,\varepsilon} = \mathfrak{a}_\varepsilon *_{\mathbb{T}_R^2} \tilde{\psi}_{N,R}$ where the convolution is done on $\mathbb{T}_R^2$ as indicated. After the approximation of the measure we need also to approximate accordingly equation (70). First we approximate the noise $B_{\varepsilon,t}^A$ by a periodic noise

$$B_{\varepsilon,R,t}^A(v) = \tilde{\Xi}\Big(\mathbb{I}_{(-\infty,t]}(\cdot) \otimes \Big(\mathbb{I}_{\mathbb{T}_R^2} e^{(\Delta - m_f^2)(t-\cdot)}(\mathfrak{a}_\varepsilon *_{\mathbb{T}_R^2} v)\Big)\Big), \quad v \in C^\infty(\mathbb{T}_R^2),$$

and we define the stationary solution to the *equation with periodic* noise

$$\tilde{\Psi}_{R,t}^h(v) = \lambda^2 \int_{-\infty}^t F_{h,rY}^{\mathbb{T}_R^2}(\tilde{\Psi}_{R,\tau}^h)(e^{(\Delta - m_f^2)(t-\tau)} v) \mathrm{d}\tau + \tilde{B}_{\varepsilon,R,t}^A, \tag{81}$$

where the nonlinearity $F_{h,rY}^{\mathbb{T}_R^2}$ is given by

$$F_{h,rY}^{\mathbb{T}_R^2}(\tilde{\Psi}_{R,\tau}^h)(v) = -\left\langle C_f \mathfrak{a}_\varepsilon^{*2} * \left[\frac{h}{2} \left(\int_{\mathbb{T}_R^2} \mathcal{G}_R(y - \cdot) h(y) \Psi_{R,t}^h(y) \bar{\Psi}_{R,t}^h(y) \mathrm{d}y\right) J \tilde{\Psi}_{R,t}^h\right], v\right\rangle$$

and $\mathcal{G}_R(x) := \sum_{n \in \mathbb{Z}^2} \mathcal{G}(x - 2\pi R n)$. We can finally introduce a *finite dimensional approximation* to equation (81) as

$$\tilde{\Psi}_{N,R,t}^h(v) = \lambda^2 \int_{-\infty}^t F_{h,rY}^{\mathbb{T}_R^2}(P_N \tilde{\Psi}_{N,R,\tau}^h)(P_N e^{(\Delta - m_f^2)(t-\tau)} v) \mathrm{d}\tau + \tilde{B}_{\varepsilon,R,t}^A(v). \tag{82}$$

**Remark 65.** We call equation (82) finite dimensional not because its solution $\tilde{\Psi}_{N,R,t}^{h_R}$ is a finite dimensional Grassmannian process but because $P_N(\tilde{\Psi}_{N,R,t}^{h_R})$ solves the finite dimensional nonlinear SDE and

$$(I - P_N)(\tilde{\Psi}_{N,R,t}^h) = (I - P_N)(\tilde{B}_{\varepsilon,R,t}^A).$$

**Remark 66.** The drift $P_N\Big(F_{h,rY}^{\mathbb{T}_R^2}(P_N(\tilde{\Psi}_{N,R,t}^h))\Big)$ of the finite dimensional process $P_N(\tilde{\Psi}_{N,R,t}^{h_R})$ is given by

$$P_N\Big(F_{h,rY}^{\mathbb{T}_R^2}(P_N(\tilde{\Psi}_{N,R,t}^h))\Big)(v) = -\frac{1}{2} \langle P_N \mathfrak{a}_\varepsilon^{*2} * C_f \partial_R U(P_N \tilde{\Psi}_{N,R,t}^h), v\rangle$$

where
$$U(P_N(\tilde\Psi))=\mathcal{V}_{R,h}(P_N(\tilde\Psi)),$$
where $\mathcal{V}_{R,h}$ is defined in equation (77).

**Remark 67.** Lemma 63 holds for the drift of equation (81), and it is easy to generalize Lemma 62 to the OU process $(\tilde B^A_{\varepsilon,R,t})_{t\in\mathbb{R}}$. This implies that it is possible to generalize Lemma 64 to equation (81) proving existence and uniqueness of its global solution (when $|\lambda|$ is small enough). Finally, it is important to note that the bound $K$ on the solution $\tilde\Psi^{h_R}_{R,t}$ to equation (81), and the bounds on the constant $|\lambda|$ can be chosen independent of $h$ when $\|h\|_{\mathcal{C}^\gamma}<C$ for some constant $C>0$.

**Remark 68.** The processes $\tilde\Psi^h_{R,t}$ and $\tilde\Psi^h_{N,R,t}$ take values in $\mathcal{C}^{\gamma+2\delta}(\mathbb{T}^2_R,\mathcal{A})$. On the other hand we can look at the processes $\tilde\Psi^h_{R,t}$ and $\tilde\Psi^h_{N,R,t}$ as *functions in $\mathcal{C}^{\gamma+2\delta}(\mathbb{R}^2,\mathcal{A})$ which are periodic (in both the coordinates) of period $2\pi R$*. More precisely, we denote by $\tilde\Psi^{p,h_R}_{R,t}:=\tilde\Psi^h_{R,t}\circ Q_R$ and $\tilde\Psi^{p,h_R}_{N,R,t}:=\tilde\Psi^h_{N,R,t}\circ Q_R$ the periodic extensions of $\tilde\Psi^h_{R,t}$ and $\tilde\Psi^h_{N,R,t}$ to the whole $\mathbb{R}^2$. The process $\tilde\Psi^{p,h_R}_{R,t}$ solve the Grassmannian SDE (on $\mathbb{R}^2$)

$$\tilde\Psi^{p,h_R}_{R,t}(v)=\lambda^2\int_{-\infty}^t F_{h_R,rY}(\tilde\Psi^{p,h_R}_{R,t})(e^{(\Delta-m_f^2)\tau}v)\mathrm{d}\tau+\tilde B^A_{\varepsilon,R,t}\circ Q_R(v)$$

where $h_R\colon\mathbb{R}^2\to\mathbb{R}$ is the $2\pi R$ periodic function defined as

$$h_R(x)=\sum_{n\in\mathbb{Z}^2}h(x-2\pi Rn),$$

and $F_{h_R,rY}$ is the same nonlinearity as $F_{h,rY}$ where the function $h$ is replaced by $h_R$.

At this point we are in position to establish the stochastic quantization connection between the approximate measure and the approximate equation, as follows.

**Lemma 69.** *Suppose that the support of $h$ is contained in $\mathbb{T}^2_R$ for some $R>0$, then for $|\lambda|$ small enough (that can be chosen in a way independent of the spatial cut-off $h$ and of $N,R$), $N\in\mathbb{N}$, $t\in\mathbb{R}$, we have that for any $k\geqslant 1$ and $f_1,...,f_k\in P_NL^2(\mathbb{T}^2_R)$,*

$$\omega\Big(\tilde\psi_R(P_N(f_1))\cdots\tilde\psi_R(P_N(f_k))e^{\mathcal{V}_{R,h}(P_N(\tilde\psi_R))}\Big)=\omega(\tilde\Psi^h_{N,R,t}(f_1)\cdots\tilde\Psi^h_{N,R,t}(f_k))\omega\big(e^{\mathcal{V}_{R,h}(P_N(\tilde\psi))}\big)$$

*where $\mathcal{V}_{R,h}$ is defined in equation (77).*

**Proof.** The proof is based on the fact that for any $f_k\in P_NL^2(\mathbb{T}^2_R)$ we have $\tilde\Psi^{h_R}_{N,R,t}(f_k)=P_N(\tilde\Psi^{h_R}_{N,R,t})(f_k)$ and on Theorem 40, Remark 65 and Remark 66. □

It remains now to remove the approximations, first letting $N\to\infty$ and then $R\to\infty$.

**Lemma 70.** *Suppose that $h$ has compact support contained in $\mathbb{T}^2_R$ for some $R>0$, then for $|\lambda|$ small enough, any $0<q<\gamma$, we have that*

$$\tilde\Psi^h_{N,R,t}\to\tilde\Psi^h_{R,t},$$

*in $\mathcal{G}^q(C^\infty(\mathbb{T}^2_R)\otimes\mathbb{C}^4)$ when $N\to+\infty$, uniformly in $t\in\mathbb{R}$. Therefore, for any $k\geqslant 1$, $f_1,...,f_k\in C^\infty(\mathbb{T}^d_R)\otimes\mathbb{C}^4$, $R>0$, and any $t\in\mathbb{R}$ we have*

$$\omega\Big(\tilde\psi_R(f_1)\cdots\tilde\psi_R(f_k)e^{\mathcal{V}_{R,h}(\tilde\psi_R)}\Big)=\omega(\tilde\Psi^h_{R,t}(f_1)\cdots\tilde\Psi^h_{R,t}(f_k))\omega\big(e^{\mathcal{V}_{R,h}(\tilde\psi_R)}\big). \qquad(83)$$

**Proof.** Recalling Definition 54 we see that we can take $\lambda_A = m_f^2$ in our case. Then we have that the difference $\Psi_{N,R,t}^h - \Psi_{R,t}^h$ satisfies the following inequality:

$$\frac{\|\tilde{\Psi}_{N,R,t}^h - \tilde{\Psi}_{R,t}^h\|_{\mathcal{C}^q}}{\lambda^2} \lesssim \int_{-\infty}^{t} e^{-\lambda_A(t-\tau)} \left\| (I - P_N) F_{h,rY}^{\mathbb{T}_R^2}(P_N(\tilde{\Psi}_{N,R,\tau}^h)) \right\|_{\mathcal{C}^q} \mathrm{d}\tau +$$

$$+ \int_{-\infty}^{t} e^{-\lambda_A(t-\tau)} g_F(\max(\|\tilde{\Psi}_{N,R,\tau}^h\|_{\mathcal{C}^q}, \|P_N \tilde{\Psi}_{N,R,\tau}^h\|_{\mathcal{C}^q})) \times$$

$$\times (\|\tilde{\Psi}_{N,R,\tau}^h - P_N \tilde{\Psi}_{N,R,\tau}^h\|_{\mathcal{C}^q}) \mathrm{d}\tau +$$

$$+ \int_{-\infty}^{t} e^{-\lambda_A(t-\tau)} g_F(\max(\|\tilde{\Psi}_{N,R,\tau}^h\|_{\mathcal{C}^q}, \|\tilde{\Psi}_{R,\tau}^h\|_{\mathcal{C}^q})) \times$$

$$\times \|\tilde{\Psi}_{N,R,\tau}^h - \tilde{\Psi}_{R,\tau}^h\|_{\mathcal{C}^q} \mathrm{d}\tau.$$

Using Gronwall inequality we obtain that

$$\|\tilde{\Psi}_{N,R,t}^h - \tilde{\Psi}_{R,t}^h\|_{\mathcal{C}^q} \lesssim (\|I - P_N\|_{\mathcal{L}(\mathcal{C}^\gamma, \mathcal{C}^q)}) \left(1 + \int_{-\infty}^{t} \exp(-(\lambda_A - \lambda^2 g_F(K))(t-\tau)) \mathrm{d}\tau\right).$$

Since $\lambda$ is small we have that $\lambda_A - \lambda^2 g_F(K) > 0$ (this is exactly the request (61) for $\lambda^2 \leqslant \lambda_0$). This means that

$$\|\tilde{\Psi}_{N,R,t}^h - \tilde{\Psi}_{R,t}^h\|_{\mathcal{C}^q} \lesssim (\|I - P_N\|_{\mathcal{L}(\mathcal{C}^\gamma, \mathcal{C}^q)}).$$

On the other hand, since $\gamma > q$ by assumption, we have $\|I - P_N\|_{\mathcal{L}(\mathcal{C}^\gamma, \mathcal{C}^q)} \to 0$ as $N \to +\infty$. The thesis is proved noting that the convergence $\tilde{\Psi}_{N,R,t}^h \to \tilde{\Psi}_{R,t}^h$ is uniform in $t$. Equation. (83) is then an easy consequence of this convergence. $\square$

In order to tackle the $R \to \infty$ limit and the infinite volume limit and prove Theorem 60 we have to introduce weighted Besov spaces (see the Appendix B and [175] for the details). Let $\ell \in \mathbb{R}$ and consider $\rho_\ell: \mathbb{R}^d \to \mathbb{R}_+$ defined by $\rho_\ell(x) = (1 + |x|^2)^{-\ell/2}$, $x \in \mathbb{R}^d$. Consider the following norm on $\mathcal{S}'(\mathbb{R}^d) \otimes \mathcal{A}$

$$\|u\|_{\mathcal{C}_\ell^s(\mathbb{R}^d, \mathcal{A})} = \sup_{j \geqslant -1} \sup_{x \in \mathbb{R}^d} 2^{js} \|\Delta_j u(x) \rho_\ell(x)\|_{\mathcal{A}}, \qquad s \in \mathbb{R}.$$

Define the space $\mathcal{C}_\ell^s(\mathbb{R}^d, \mathcal{A})$ as the subspace of $\mathcal{S}'(\mathbb{R}^d) \otimes \mathcal{A}$ where $\|\cdot\|_{\mathcal{C}_\ell^s(\mathbb{R}^d, \mathcal{A})}$ is finite. We define also $\mathcal{G}_\ell^s(V) = \mathcal{C}_\ell^s(\mathbb{R}^d, \mathcal{A}) \cap \mathcal{G}^{-\infty}(V)$ (recall that $V = \mathcal{S}(\mathbb{R}^d) \otimes \mathbb{R}^n$ or $\mathcal{S}(\mathbb{R}^d) \otimes \mathbb{C}^n$ for some $n \geqslant 1$) and the distance $d_{\mathcal{G}_\ell^s(V)}(X - Y) := \|X - Y\|_{\mathcal{C}_\ell^s(\mathbb{R}^d, \mathcal{A})}$.

**Remark 71.** It is important to note that, by Theorem 79 below,

$$\sup_{x \in \mathbb{R}^d} \|\Delta_j u(x) \rho_\ell(x)\|_{\mathcal{A}} \sim \sup_{x \in \mathbb{R}^d} \|\Delta_j(\rho_\ell u)(x)\|_{\mathcal{A}},$$

where $\sim$ means that both quantities considered can be bounded from above and from below by some positive constant times the other quantity.

Hereafter we identify $\tilde{B}_{\varepsilon, R, t}^A$ (that in principle is defined on $\mathbb{T}_R^2$) with the periodic extension $\tilde{B}_{\varepsilon, R, t}^{p, A} := \tilde{B}_{\varepsilon, R, t}^A \circ Q_R$ of $\tilde{B}_{\varepsilon, R, t}^A$ to $\mathbb{R}^2$.

**Lemma 72.** *Let $\gamma, \delta$ be as stated in Remark 59, for any $s < \gamma + 2\delta$ and any $\ell > 0$ we have that*

$$\sup_{t \in \mathbb{R}} d_{\mathcal{G}_\ell^s(V)}(\tilde{B}_{\varepsilon, t}^A, \tilde{B}_{\varepsilon, R, t}^A) \to 0$$

when $R \to +\infty$.

The proof of Lemma 72 is postponed to Section 5.4.

**Lemma 73.** *Let $|\lambda|$ be sufficiently small and $0 < q < \gamma$, then*

$$\tilde{\Psi}_{R,t}^{p,h_R} \to \tilde{\Psi}_t^h,$$

*in $\mathcal{G}_\ell^q(\mathcal{S}(\mathbb{R}^2) \otimes \mathbb{C}^4)$ uniformly in $t \in \mathbb{R}$ when $R \to +\infty$ (where $\tilde{\Psi}_{R,t}^{p,h_R} := \tilde{\Psi}_{R,t}^h \circ Q_R$ is defined in Remark 68). Moreover, let $(h_n \colon \mathbb{R}^2 \to \mathbb{R})_{n \geqslant 1}$ be a sequence of smooth functions with compact supports converging to $1$ uniformly on compacts in the norm of $\mathcal{C}^\gamma$ and with a uniform bound in $\mathcal{C}^\gamma$, $\gamma > 0$. We have also*

$$\tilde{\Psi}_t^{h_n} \to \tilde{\Psi}_t,$$

*in $\mathcal{G}_\ell^q(\mathcal{S}(\mathbb{R}^2) \otimes \mathbb{C}^4)$, for any $\ell > 0$ and uniformly in $t \in \mathbb{R}$ when $n \to +\infty$.*

**Remark 74.** Since $\tilde{\Psi}_{R,t}^{p,h_R}(x) = \tilde{\Psi}_{R,t}^h(x)$ for any $x \in \mathbb{T}_R^2 = (-\pi R, \pi R]^2$, Lemma 73 implies that, for any $\mathcal{K} \subset \mathbb{R}^2$ compact set, which is a subset of $\mathbb{T}_R^2$ for $R$ big enough, we get

$$\sup_{x \in \mathcal{K}} \|\tilde{\Psi}_{R,t}^h(x) - \tilde{\Psi}_t^h(x)\|_{\mathcal{A} \otimes \mathbb{C}^4} \to 0$$

as $R \to +\infty$.

**Proof of Lemma 73.** The proof is based on the fact that for any $\tilde{\Psi}_1, \tilde{\Psi}_2 \in \mathcal{G}_0^\gamma(\mathcal{S}(\mathbb{R}^2) \otimes \mathbb{C}^4)$ (it is important to note that in these spaces we take the constant weight) we have, by Theorem 82 in Appendix B,

$$\|F_{h_R,rY}(\tilde{\Psi}_1) - F_{h,rY}(\tilde{\Psi}_2)\|_{\mathcal{C}_\ell^q(\mathbb{R}^2, \mathcal{A} \otimes \mathbb{C}^4)} \lesssim g_F(\max(\|h_R \tilde{\Psi}_1\|_{\mathcal{C}^q}, \|h \tilde{\Psi}_1\|_{\mathcal{C}^q}, \|h \tilde{\Psi}_2\|_{\mathcal{C}^q})) \times$$
$$\times \big(\|h\|_{\mathcal{C}^q} \|\tilde{\Psi}_1 - \tilde{\Psi}_2\|_{\mathcal{C}_\ell^q(\mathbb{R}^2, \mathcal{A} \otimes \mathbb{C}^4)} + \|\Psi_1\|_{\mathcal{C}_0^q(\mathbb{R}^2, \mathcal{A} \otimes \mathbb{C}^4)} \|h_R - h\|_{\mathcal{C}_\ell^q}\big), \tag{84}$$

where $g_F$ is the function in Remark 53, $h$ the IR cut-off in Definition 58, and $h_R$ is the periodicization of $h$ (see Remark 68). Using a reasoning similar to the one exploited in the proof of Lemma 70 we get, for any $t \in \mathbb{R}$,

$$\|\tilde{\Psi}_{R,t}^{p,h_R} - \tilde{\Psi}_t^h\|_{\mathcal{C}_\ell^q(\mathbb{R}^2, \mathcal{A} \otimes \mathbb{C}^4)} \lesssim \bigg(\|h_R - h\|_{\mathcal{C}_\ell^q} + \sup_{\tau \in \mathbb{R}} \|\tilde{B}_{\varepsilon,\tau}^A - \tilde{B}_{\varepsilon,R,\tau}^A\|_{\mathcal{C}_\ell^q(\mathbb{R}^2, \mathcal{A} \otimes \mathbb{C}^4)}\bigg) \times$$
$$\times \bigg(1 + \int_{-\infty}^t \exp(-(\lambda_A - \lambda^2 g_F(K))(t-\tau)) \mathrm{d}\tau\bigg).$$

By $|\lambda|$ small enough, the bound converges to $0$ as $R \to +\infty$ uniformly in $t \in \mathbb{R}$, since, by Lemma 72,

$$\lim_{R \to +\infty} \sup_{\tau \in \mathbb{R}} \|\tilde{B}_{\varepsilon,\tau}^A - \tilde{B}_{\varepsilon,R,\tau}^A\|_{\mathcal{C}_\ell^q(\mathbb{R}^2, \mathcal{A} \otimes \mathbb{C}^4)} \to 0$$

and $\|h_R - h\|_{\mathcal{C}_\ell^q} \to 0$, as $R \to +\infty$, by the definition of $h_R$.

For the convergence $\tilde{\Psi}_t^{h_n} \to \tilde{\Psi}_t$ when $h_n \to 1$, the proof is similar considering that in the inequality corresponding to (84) we will have $\|h_n - 1\|_{\mathcal{C}_\ell^q}$ instead of $\|h_R - h\|_{\mathcal{C}_\ell^q}$ which converges to $0$, as $n \to +\infty$, by hypothesis, when $\ell > 0$. □

## 5.4 Some auxiliary results

We close this section with the proofs of Lemma 62 on the regularity of the noise $B_{\varepsilon,t}^A$ and of Lemma 72 on the convergence of the periodic noise $B_{\varepsilon,R,t}^A$ as $R \to \infty$.

**Proof of Lemma 62.** Since, by the requirement of Remark 59 and by Theorem 81 in Appendix B, the map $\mathfrak{a}_\varepsilon * \cdot$ improves the Besov regularity by $\alpha$ (namely $\mathfrak{a}_\varepsilon *: \mathcal{C}^s \to \mathcal{C}^{s+\alpha}$ for any $s \in \mathbb{R}$) it is enough to prove that (see equation (73))

$$\tilde{B}_{0,t}^A := \tilde{\Xi}(\mathbb{I}_{(-\infty,t]} \otimes e^{A(t-\cdot)}(v))$$

is in $\mathcal{C}^\delta(\mathbb{R}, \mathcal{C}^{-\frac{1}{2}-2\delta}(\mathbb{R}^2, \mathcal{A}) \otimes \mathbb{C}^4)$. By definition of Grassmann Gaussian random variables, if $K_i \in \mathcal{S}(\mathbb{R}^2)$ is the function corresponding to the Littlewood–Paley block $\Delta_i$ (i.e. $K_i = \check{\varphi}_i$ the inverse Fourier transform of the dyadic partition of the unity $\{\varphi_i\}_{i \geqslant -1}$, see Appendix B for the precise definition), we have

$$\sup_{x \in \mathbb{R}^2} \|\Delta_i \tilde{B}_{0,t} - \Delta_i \tilde{B}_{0,s}\|_{\mathcal{A} \otimes \mathbb{C}^4} \lesssim$$
$$\lesssim \left\| \mathbb{I}_{(-\infty,t]}(\tau) e^{(\Delta - m_f^2)(t-\tau)} (-\Delta + m^2)^{1/4}(K_i) - \right.$$
$$\left. - \mathbb{I}_{(-\infty,s]}(\tau) e^{(\Delta - m_f^2)(s-\tau)} (m_f^2 - \Delta)^{1/4}(K_i) \right\|_{L^2(\mathbb{R}^2 \times \mathbb{R})}$$
$$\lesssim \left\| \mathbb{I}_{(s,t]}(\tau) e^{(\Delta + m_f^2)(t-\tau)} (m_f^2 - \Delta)^{1/4}(K_i) \right\|_{L^2(\mathbb{R}^2 \times \mathbb{R})} +$$
$$+ \left\| \mathbb{I}_{(-\infty,s]}(\tau) \big( e^{(\Delta - m_f^2)(t-s)} - 1 \big) e^{-(\Delta - m_f^2)(s-\tau)} (m_f^2 - \Delta)^{1/4}(K_i) \right\|_{L^2(\mathbb{R}^2 \times \mathbb{R})}$$
$$\lesssim \left( \int_{\mathbb{R}^2} \int_s^t \phi_i(p)^2 e^{-2(|p|^2 + m_f^2)(t-\tau)} |p| \mathrm{d}p \mathrm{d}\tau \right)^{1/2} +$$
$$+ \left( \int_{\mathbb{R}^2} \int_{-\infty}^s \left( \int_0^{t-s} |p|^2 e^{-(|p|^2 + m_f^2)k} \mathrm{d}k \right)^2 \phi_i(p)^2 e^{-2(|p|^2 + m_f^2)(s-\tau)} |p| \mathrm{d}p \mathrm{d}\tau \right)^{1/2}$$
$$\lesssim \left( \int_{2^{i-1} \lesssim |p| \lesssim 2^i} \left( \int_0^{t-s} \frac{1}{\tau^{1-2\delta}} \mathrm{d}\tau \right) |p|^{-1+4\delta} \mathrm{d}p \right)^{1/2} +$$
$$+ \left( \int_{2^{i-1} \lesssim |p| \lesssim 2^i} \left( \int_0^{t-s} \frac{1}{k^{1-\delta}} \mathrm{d}k \right)^2 |p|^{-1+4\delta} \mathrm{d}p \right)^{1/2}$$
$$\lesssim 2^{i\left(\frac{1}{2}+2\delta\right)} |t-s|^\delta.$$

From this, the lemma follows by the definition of the norm of the Besov space $\mathcal{C}^{-\frac{1}{2}-2\delta}(\mathbb{R}^2, \mathcal{A})$. $\square$

**Proof of Lemma 72.** Using a method similar to the one in the proof of Lemma 62, it is simple to see that

$$\sup_{R>1} \left( \sup_{t \in \mathbb{R}} \|(\Delta+1)\tilde{B}_{\varepsilon,R,t}^A\|_{\mathcal{C}^{\gamma+2\delta-2}(\mathbb{R}^2, \mathcal{A}) \otimes \mathbb{C}^4} \right) < +\infty,$$

which implies that, for any $s < \gamma + \delta$, we have

$$\|\Delta_i(\tilde{B}_{\varepsilon,R,t}^A - \tilde{B}_{\varepsilon,t}^A)\|_\infty \leqslant C_1 2^{-i(\gamma+2\delta)} \tag{85}$$

where the constant $C_1 \in \mathbb{R}_+$ does not depend on $i$ and $R > 1$. This means that, fixing $s < s' < \gamma + 2\delta$, we have

$$\|\tilde{B}_{\varepsilon,R,t}^A - \tilde{B}_{\varepsilon,t}^A\|_{\mathcal{C}_\ell^s(\mathbb{R}^2, \mathcal{A}) \otimes \mathbb{C}^4} \lesssim \sum_{i \geqslant -1} 2^{is'} \|\Delta_i(\tilde{B}_{\varepsilon,R,t}^A - \tilde{B}_{\varepsilon,t}^A) \rho_\ell\|_\infty \tag{86}$$

and, thus by the uniform estimate (85) and Lebesgue's dominated convergence theorem, if, for any $i \geqslant -1$, $\|\Delta_i(\tilde{B}_{\varepsilon,R,t}^A - \tilde{B}_{\varepsilon,t}^A)\rho_\ell\|_\infty \to 0$ as $R \to +\infty$, then $\|\tilde{B}_{\varepsilon,R,t}^A - \tilde{B}_{\varepsilon,t}^A\|_{\mathcal{C}_\ell^s(\mathbb{R}^2, \mathcal{A}) \otimes \mathbb{C}^4} \to 0$ as $R \to +\infty$, and thus, by inequality (86), $\|\tilde{B}_{\varepsilon,R,t}^A - \tilde{B}_{\varepsilon,t}^A\|_{\mathcal{C}_\ell^s(\mathbb{R}^2, \mathcal{A}) \otimes \mathbb{C}^4} \to 0$ as $R \to +\infty$.

Fix $i \geqslant 1$, then we have

$$\|\Delta_i(\tilde{B}^A_{\varepsilon,R,t} - \tilde{B}^A_{\varepsilon,t})(x)\rho_\ell(x)\|_{\mathcal{A}\otimes\mathbb{C}^4} \leqslant$$
$$\leqslant \|\rho_\ell(x)\Xi_R(\tilde{K}_i(x-\cdot)\mathcal{G}_t(\cdot)\mathbb{I}_{[-\kappa R,\kappa R]^2}) - \rho_\ell(x)\Xi(\tilde{K}_i(x-\cdot)\mathcal{G}_t(\cdot)\mathbb{I}_{[-\kappa R,\kappa R]^2})\|_{\mathcal{A}\otimes\mathbb{C}^4} +$$
$$+ \|\rho_\ell(x)\Xi_R(\tilde{K}_i(x-\cdot)\mathcal{G}_t(\cdot)\mathbb{I}_{([-\kappa R,\kappa R]^2)^c}) - \rho_\ell(x)\Xi(\tilde{K}_i(x-\cdot)\mathcal{G}_t(\cdot)\mathbb{I}_{([-\kappa R,\kappa R]^2)^c})\|_{\mathcal{A}\otimes\mathbb{C}^4} \quad (87)$$
$$\leqslant \|\rho_\ell(x)\Xi_R(\tilde{K}_i(x-\cdot)\mathcal{G}_t(\cdot)\mathbb{I}_{([-\kappa R,\kappa R]^2)^c})\|_{\mathcal{A}\otimes\mathbb{C}^4} +$$
$$+ \|\rho_\ell(x)\Xi(\tilde{K}_i(x-\cdot)\mathcal{G}_t(\cdot)\mathbb{I}_{([-\kappa R,\kappa R]^2)^c})\|_{\mathcal{A}\otimes\mathbb{C}^4},$$

where $\tilde{K}_i$ is given by

$$\tilde{K}_i(x) := (m_f^2 - \Delta)^{1/4}(K_i),$$

$K_i$ being the function related to the block $i$-th Littlewood-Paley, i.e. $\Delta_i f = K_i * f$, $\mathcal{G}_t$ is the Green's function of the operator $e^{-(m_f^2-\Delta)t}$, and where we use that $\Xi_R(\tilde{K}_i(x-\cdot)\mathcal{G}_t(\cdot)\mathbb{I}_{[-\kappa R,\kappa R]^2}) = \Xi(\tilde{K}_i(x-\cdot)\mathcal{G}_t(\cdot)\mathbb{I}_{[-\kappa R,\kappa R]^2})$ since they coincide on $[-\kappa R, \kappa R]^2$ for any $0 \leqslant \kappa < \frac{1}{2}$. We focus on the term $J(x) := \|\rho_\ell(x)\Xi(\tilde{K}_i(x-\cdot)\mathcal{G}_t(\cdot)\mathbb{I}_{([-\kappa R,\kappa R]^2)^c}(\cdot))\|_{\mathcal{A}\otimes\mathbb{C}^4}$. In this case we can write for $|x| > \kappa_1 R$

$$J(x) \lesssim \frac{1}{(\kappa_1 R)^\ell} \sup_{x\in\mathbb{R}^2}\left|\int_0^{+\infty}\int_{\mathbb{R}^4}\tilde{K}_i(x-y_1)\tilde{K}_i(x-y_2) \times \right.$$
$$\left.\left(\int_{([-\kappa R,\kappa R]^2)^c}\mathcal{G}_t(y_1-z)\mathcal{G}_t(y_2-z)\mathrm{d}z\right)\mathrm{d}y_1\mathrm{d}y_2\mathrm{d}t\right|$$
$$\lesssim \frac{1}{(\kappa_1 R)^\ell}\left|\int_{\mathbb{R}_+\times\mathbb{R}^4}\exp(-\beta_i|y_1|^{\alpha/2} - \beta_i|y_2|^{\alpha/2}) \times \right.$$
$$\left.\times\left(\int_{\mathbb{R}^2}\mathcal{G}_t(y_1-z)\mathcal{G}_t(y_2-z)\mathrm{d}z\right)\mathrm{d}y_1\mathrm{d}y_2\mathrm{d}t\right|$$
$$\lesssim \frac{1}{(\kappa_1 R)^\ell}\int_{\mathbb{R}^4}e^{(-\beta_i|y_1|^{\alpha/2}-\beta_i|y_2|^{\alpha/2})}\left|\int_0^{+\infty}\exp\left(-\frac{|y_1-y_2|^2}{2t} - 2m^2 t\right)\frac{\mathrm{d}t}{t}\right|\mathrm{d}y_1\mathrm{d}y_2$$
$$\lesssim \frac{1}{(\kappa_1 R)^\ell}\int_{\mathbb{R}^4}e^{(-\beta_i|y_1|^{\alpha/2}-\beta_i|y_2|^{\alpha/2})}(|\log(|y_1-y_2|)|+1)\mathrm{d}y_1\mathrm{d}y_2$$
$$\lesssim \frac{1}{(\kappa_1 R)^\ell},$$

where we used the fact that there are $\beta_i$ (depending on $i \geqslant -1$) and $\alpha \in (0,1)$ such that $|\tilde{K}_i(y)| \lesssim e^{-\beta_i|y|^\alpha}$ (see Remark 78 below). Consider now $|x| \leqslant \kappa_1 R$, then we have

$$J(x) \lesssim \left|\int_0^{+\infty}\int_{|y_1|,|y_2|>\kappa_2 R}\tilde{K}_i(x-y_1)\tilde{K}_i(x-y_2)P(y_1,y_2,t)\mathrm{d}y_1\mathrm{d}y_2\mathrm{d}t\right| +$$
$$\left|\int_0^{+\infty}\int_{|y_1|,|y_2|<\kappa_2 R}\tilde{K}_i(x-y_1)\tilde{K}_i(x-y_2)P(y_1,y_2,t)\mathrm{d}y_1\mathrm{d}y_2\mathrm{d}t\right| +$$
$$2\left|\int_0^{+\infty}\int_{|y_1|>R,|y_2|<\kappa_2 R}\tilde{K}_i(x-y_1)\tilde{K}_i(x-y_2)P(y_1,y_2,t)\mathrm{d}y_1\mathrm{d}y_2\mathrm{d}t\right|$$
$$\lesssim J_1(x) + J_2(x) + J_3(x),$$

where $\kappa_1 < \kappa_2 < \kappa < \frac{1}{2}$ and

$$P(y_1,y_2,t) = \left(\int_{([-\kappa R,\kappa R]^2)^c}\mathcal{G}_t(y_1-z)\mathcal{G}_t(y_2-z)\mathrm{d}z\right)$$

Then we get

$$\begin{aligned}
J_1(x) &\lesssim \int_0^{+\infty}\int_{|y_1|,|y_2|>\kappa_2 R} \exp(-\beta_i|x-y_1|^\alpha-\beta_i|x-y_1|^\alpha)P(y_1,y_2,t)\mathrm{d}y_1\mathrm{d}y_2\mathrm{d}t \\
&\lesssim \int_0^{+\infty}\int_{|y_1|,|y_2|>\kappa_2 R} \exp(-\beta_i(|y_1|-|x|)^\alpha-\beta_i(|y_2|-|x|)^\alpha) \times \\
&\qquad \left(\int_{\mathbb{R}^2}\mathcal{G}_t(y_1-z)\mathcal{G}_t(y_2-z)\mathrm{d}z\right)\mathrm{d}y_1\mathrm{d}y_2\mathrm{d}t \\
&\lesssim (e^{-\beta_i(\kappa_2-\kappa_1)^\alpha R^\alpha}(1+R^2)) \times \\
&\qquad \int_{\mathbb{R}^2} e^{(-\beta_i|y_1'|^\alpha-\beta_i|y_2'|^\alpha)}\left|\int_0^{+\infty}\exp\left(-\frac{|y_1-y_2|^2}{2t}-2m^2t\right)\frac{\mathrm{d}t}{t}\right|\mathrm{d}y_1\mathrm{d}y_2 \\
&\lesssim (e^{-\beta_i(\kappa_2-\kappa_1)^\alpha R^\alpha}(1+R^2))\int_{\mathbb{R}^2} e^{(-\beta_i|y_1'|^\alpha-\beta_i|y_2'|^\alpha)}(|\log(|y_1-y_2|)|+1)\mathrm{d}y_1\mathrm{d}y_2 \\
&\lesssim (e^{-\beta_i(\kappa_2-\kappa_1)^\alpha R^\alpha}(1+R^2))
\end{aligned}$$

For $J_2(x)$ we have

$$\begin{aligned}
J_2(x) &\lesssim (\kappa_2 R^2)^2 \sup_{|y_1|<\kappa_2 R,|y_1|<\kappa_2 R}\left|\int_0^{+\infty}\left(\int_{([-\kappa R,\kappa R]^2)^c}\mathcal{G}_t(y_1-z)\mathcal{G}_t(y_2-z)\mathrm{d}z\right)\mathrm{d}t\right| \\
&\lesssim (\kappa_2 R^2)^2 \int_0^{+\infty} e^{\frac{(\kappa-\kappa_2)^2 R^2}{2t}-2m^2 t}\frac{\mathrm{d}t}{t} \\
&\lesssim \frac{(\kappa-\kappa_2)^{2n}}{R^{2n}}(\kappa_2 R^2)^2\int_0^{+\infty} t^{n-1}e^{-2m^2 t}\mathrm{d}t \lesssim \frac{(\kappa-\kappa_2)^2\kappa_2^2}{R^{2n-4}}.
\end{aligned}$$

where we used that $\exp(-x)\lesssim x^{-n}$ for any $n\in\mathbb{N}$. The term $J_3$ can be estimated in the same way combining the techniques used for $J_1$ and $J_2$, In conclusion

$$\sup_{x\in\mathbb{R}^2} J(x) \lesssim \left(\frac{(\kappa-\kappa_2)^2\kappa_2^2}{R^{2n-4}}+\frac{1}{(\kappa_1 R)^\ell}+e^{-\beta_i(\kappa_2-\kappa_1)^\alpha R^\alpha}(1+R^2)\right)\to 0$$

when $n>2$ and $\ell>0$.

For $\|\Xi_R(\tilde{K}_i(x-\cdot)\mathcal{G}_t(\cdot)\mathbb{I}_{([-\kappa R,\kappa R]^2)^c})\|_{\mathcal{A}\otimes\mathbb{C}^4}$ we have

$$\begin{aligned}
\|\Xi_R(\tilde{K}_i(x-\cdot)\mathcal{G}_t(\cdot)&\mathbb{I}_{([-\kappa R,\kappa R]^2)^c})\|_{\mathcal{A}\otimes\mathbb{C}^4}\lesssim \\
&\lesssim \left|\int_0^{+\infty}\int_{\mathbb{R}^4}\tilde{K}_i(x-y_1)\tilde{K}_i(x-y_2)P_R(y_1,y_2,t)\mathrm{d}y_1\mathrm{d}y_2\mathrm{d}t\right|,
\end{aligned} \qquad (88)$$

where

$$P_R(y_1,y_2,t):=\int_{\mathbb{T}_R^2\setminus([-\kappa R,\kappa R]^2)^c}\mathcal{G}_t(y_1-z)\mathcal{G}_t(y_2-z)\mathrm{d}z+\sum_{z'\in R\mathbb{Z}^2}\int_{\mathbb{T}_R^2}\mathcal{G}_t(y_1-z)\mathcal{G}_t(y_2-z)\mathrm{d}z.$$

By noting that for $|y_1|,|y_2|>\kappa_2 R$

$$\int_0^{+\infty} P_R(y_1,y_2,t)\mathrm{d}t \leqslant C_2(|\log(|y_1-y_2|)|+1)$$

for some constant $C_2\in\mathbb{R}_+$ independent of $R>1$ and, if one of the inequalities $|y_1|<\kappa_2 R$ or $|y_2|<\kappa_2 R$ holds, then

$$P_R(y_1,y_2,t)\leqslant \frac{\exp\left(-\beta'\frac{(\kappa-\kappa_2)^2 R^2}{t}+2m_f^2 t\right)}{t}.$$

Then using inequality (88), it is possible to prove that

$$\sup_{x \in \mathbb{R}^2} \|\rho_\ell(x) \Xi_R(\tilde{K}_i(x-\cdot)\mathcal{G}_t(\cdot)\mathbb{I}_{([-\kappa R, \kappa R]^2)^c})\|_{\mathcal{A} \otimes \mathbb{C}^4} \to 0$$

as $R \to +\infty$ by exploiting an argument analogous to the one mentioned above for proving $\sup_{x \in \mathbb{R}^2} \|\rho_\ell(x) \Xi(\tilde{K}_i(x-\cdot)\mathcal{G}_t(\cdot)\mathbb{I}_{([-\kappa R, \kappa R]^2)^c})\|_{\mathcal{A} \otimes \mathbb{C}^4} \to 0$. $\square$

# Appendix A Convergence of the perturbative series: finite dimensional case

In this appendix we sketch some implications of the Pauli exclusion principle for the existence of global solutions to finite-dimensional Grassmann SDEs. In principle some of these considerations also apply to some more realistic models like the Yukawa$_2$ model; however, to our surprise, stochastic quantization, at least in the perturbative regime, can be carried on without establishing them, as we demonstrated in Section 5.

The solutions to the non-linear equation (43) can be represented via a series. We investigate now the properties of this series and its absolute convergence. For simplicity let us assume that $(v_\alpha)_{\alpha=1,\ldots,N}$ is a finite basis of $V$ and that $F \in \mathrm{Hom}(V, \Lambda V)$ is given by a sum of cubic monomials

$$F(v_\alpha) = \sum_{\alpha_1, \alpha_2, \alpha_3} \lambda^\alpha_{\alpha_1 \alpha_2 \alpha_3} \psi(v_{\alpha_1}) \psi(v_{\alpha_2}) \psi(v_{\alpha_2}),$$

where $\psi \in \mathrm{Hom}(V, \Lambda V)$ is the canonical injection of $V$ into $\Lambda V$ and $(\lambda^\alpha_{\alpha_1 \alpha_2 \alpha_3})$ a family of coefficients in $\mathbb{R}$. Moreover we take $A = -\mathbb{I}$. Equation (43) has then the integral formulation given by

$$\Psi_t(v) = \Phi_t(v) + \int_0^t \Psi_s(e^{-(t-s)} F(v)) \mathrm{d}s, \qquad t \geqslant 0, v \in V, \tag{89}$$

where

$$\Phi_t(v) := \Psi_0(e^{-t} v) + B_t^A(v) - B_0^A(e^{-t} v), \qquad t \geqslant 0, v \in V.$$

By iteratively expanding $\Psi_s$ on the right hand side of (89) we obtain a series expansion for $\Psi_t$ of the form

$$\Psi_t = \sum_\tau J_\tau(\Phi)(t) \tag{90}$$

$$= J_\bullet(\Phi)(t) + J_{[\bullet\bullet\bullet]}(\Phi)(t) + J_{[[\bullet\bullet\bullet]\bullet\bullet]}(\Phi)(t) + \cdots + J_{[[\bullet\bullet\bullet][\bullet\bullet[\bullet\bullet\bullet]]\bullet]\bullet]}(\Phi)(t) + \cdots$$

The series is indexed by planar trees $\tau$ which have branches of order 3 and where $J$ is a multilinear integral operator such that

$$J_\bullet(\Phi)(t)^\alpha = \Phi_t(v_\alpha)$$

$$J_{[\tau_1 \tau_2 \tau_3]}(\Phi)(t)^\alpha = \sum_{\alpha_1, \alpha_2, \alpha_3} \int_0^t e^{-(t-s)} \lambda^\alpha_{\alpha_1, \alpha_2, \alpha_3} J_{\tau_1}(\Phi)(s)^{\alpha_1} J_{\tau_2}(\Phi)(s)^{\alpha_2} J_{\tau_3}(\Phi)(s)^{\alpha_3} \mathrm{d}s$$

where $\bullet$ denotes the simple tree and $[\tau_1, \ldots, \tau_3]$ the tree with branches $\tau_1, \ldots, \tau_3$. Our goal is to prove that the above series converges for all times and derive estimates for the norm $\|\Psi_t\|$ of the solution $\Psi$. Expanding the expression for $J_\tau(\Phi)$ we have

$$J_\tau(\Phi) = \int_{\mathbb{R}_+^{|I(\tau)|}} \left( \prod_{p \in I(\tau)} \lambda^{\alpha_p}_{\alpha_{p_1} \alpha_{p_2} \alpha_{p_3}} \right) \left( \prod_{p \in L(\tau)} \Phi^{\alpha_p}(s_{p_-}) \right) \left( \prod_{p \in I(\tau)} \mathbb{1}_{s_p \leqslant s_{p_-}} e^{-(s_{p_-} - s_p)} \mathrm{d}s_p \right)$$

where $I(\tau)$ denotes the internal nodes, $L(\tau)$ the leaves, $p$ the parent node and $p_1, p_2, p_3$ the three children of every internal node, $p_-$ denotes the parent of the node $p$ with the convention that $s_{p_-} = t$ if $p$ in the root, the products are taken with the order induced by visiting the nodes in a deep-first manner, and $\Phi^\alpha(t) = \Phi_t(v_\alpha)$. We have $|\tau| = |L(\tau)| + |I(\tau)|$, $|\tau| = 1 + 3|I(\tau)|$, $|L(\tau)| = 1 + 3|I(\tau)| - |I(\tau)| = 1 + 2|I(\tau)|$.

Note that $(\Phi_t^\alpha)^2 = 0$ and that the increments $\Phi^\alpha(t,s) := \Phi_t^\alpha - \Phi_s^\alpha$ have the norm bound

$$\|\Phi^\alpha(t,s)\|^2 \leqslant \int_s^t e^{-(t-r)} \mathrm{d}r + |e^{-(t-s)} - 1|^2 \leqslant 2(1 - e^{-(t-s)}) =: H(t-s)$$

which can be estimated as

$$\|\Phi^\alpha(t,s)\| \leqslant |t-s|^{1/2}.$$

The following Lemma uses $(\Phi_t^\alpha)^2 = 0$ (Pauli exclusion principle) to derive good estimates for products of fields. We observe that

$$\|\|\Phi\|\| := \sup_\alpha \left[ \sup_{t \geqslant 0} \|\Phi_t^\alpha\| + \sup_{t > s \geqslant 0} \frac{\|\Phi_t^\alpha - \Phi_s^\alpha\|}{|t-s|^{1/2}} \right] < \infty.$$

**Lemma 75.** *For any $n \geqslant 1$ and any $t_1, ..., t_n \in [0, T]$, $\alpha_1, ..., \alpha_n \in \{1, ..., N\}$ we have*

$$\|\Phi_{t_1}^{\alpha_1} \cdots \Phi_{t_n}^{\alpha_n}\| \leqslant \frac{C^{n+1} T^{n/8}}{(n!)^{1/8}} \|\|\Phi\|\|^n$$

*for all $n \geqslant 4N$, where $C$ is a universal constant depending only on $N$.*

$$\frac{n}{2N} \leqslant r \leqslant \frac{n+2N}{2N}$$

$$\leqslant \frac{T}{r} \leqslant \frac{2NT}{n}$$

$$\frac{n}{4} \leqslant \frac{n}{2} - \frac{n}{4} \leqslant \frac{n}{2} - N \leqslant n - Nr \leqslant \frac{n}{2}$$

**Proof.** Partition $[0,T]$ in $r = \lceil n/(2N) \rceil$ equal intervals $(I_k)_k$ of size $|I_k| \leqslant 2NT/n$ and let $(s_k)_k$ be the centers of those intervals. Now in the product replace each $\Phi_{t_1}^{\alpha_i}$ by $\Phi_{t_1}^{\alpha_i} = \Phi_{s_k}^{\alpha_i} + \Phi^{\alpha_i}(t_i, s_k)$ where $s_k$ is the nearest to $t_i$ of the centers. By doing so we rewrite $\Phi_{t_1}^{\alpha_1} \cdots \Phi_{t_n}^{\alpha_n}$ as a sum $S$ of $2^n$ products of fields which are either $\Phi_{s_k}^{\alpha_i}$ or $\Phi^{\alpha_i}(t_i, s_k)$, moreover for a given interval $I_k$ there cannot be more than $N$ factors $\Phi_{s_k}^{\alpha_i}$ by the exclusion principle so overall we have at most $Nr$ factors $\Phi_{s_k}^{\alpha_i}$ and $n - Nr \geqslant n/2 - N \geqslant n/4$ factors $\Phi^{\alpha_i}(t_i, s_k)$. Then use the fact that

$$\|\Phi_{s_k}^{\alpha_i}\| \leqslant \|\|\Phi\|\|, \quad \|\Phi^{\alpha_i}(t_i, s_k)\| \leqslant \|\|\Phi\|\| |t_i - s_k|^{1/2} \leqslant \|\|\Phi\|\| (2NT/r)^{1/2}$$

to estimate

$$\frac{\|S\|}{\|\|\Phi\|\|^n} \leqslant 2^n \left( \frac{2NT}{n} \right)^{\frac{n}{8}} = 2^n \left( \frac{2NT}{n} \right)^{n/8} \leqslant \frac{C^{n+1} T^{n/8}}{(n!)^{1/8}}. \qquad \square$$

Using this lemma above we can estimate

$$\|J_\tau(\Phi)(t)\| \leqslant C |\lambda|^{|I(\tau)|} (CN\|\|\Phi\|\|)^{|L(\tau)|} \left[ \frac{t^{|L(\tau)|}}{|L(\tau)|!} \right]^{1/8},$$

for all $\tau$ with $|L(\tau)| \geqslant 4N$, where we estimated all the integrals by constants uniformly in $t$ due to the presence of the exponential factors. The number of trees $\tau$ with a given number of nodes $|\tau| = n \geqslant 4N$ is no more than $D^n$ for some $n$ and $D > 0$ and $|L(\tau)| = 1 + 2|I(\tau)| = 1 + (2/3)(|\tau| - 1) = 1/3 + 2|\tau|/3$ so we obtain at the end for the solution $\Psi_t$ of (89)

$$\begin{aligned} \|\Psi_t\| &\leqslant \sum_\tau \|J_\tau(\Phi)(t)\| \lesssim (1 + \|\!|\Phi|\!\|)^{4N-1} \\ &+ \sum_{n \geqslant 4N} D^n |\lambda|^{(2/3)(n-1)} N^{2n/3} \frac{C^{2n/3} t^{n/12+1/24}}{(n/12+1/24)!} \|\!|\Phi|\!\|^n \end{aligned} \quad (91)$$

which is a series which converges for all $t$ (with some stretched exponential behavior).

**Theorem 76.** *There exists an increasing function $E(t)$ depending on $N, |\lambda|, \|\!|\Phi|\!\|$ such that*

$$\sup_{s \leqslant t} \|\Psi_s\| \leqslant E(t),$$

*where we recall that $(\Psi_t)_{t \geqslant 0}$ is the unique solution of (89).*

In particular this shows again that explosion is not possible for this kind of SDEs, unlike their bosonic analogs.

Our goal now is to remove the time dependence on our estimate for the solution exploiting the exponential decay. Let us partition the interval $[0, \infty]$ into intervals $(I_h)_{h \geqslant 0}$ of size 1. In each of these intervals we will use the finer partition of the previous lemma to estimate products, which will then not depend on $T$ anymore. Of course now the problem is that we have a bound of the form

$$\|\Phi^{\alpha_1}(t_1) \cdots \Phi^{\alpha_n}(t_n)\| \leqslant \|\!|\Phi|\!\|^n \prod_{h \geqslant 0} \frac{C^{n_h+1}}{(n_h!)^{1/8}}$$

where $n_h$ is the number of time variables in $I_h$ and $\sum_{h \geqslant 0} n_h = n \geqslant 4N$, and so far there is no useful upper bound on this (since we observe that we can have the situation where $n_h = 1$ for all $h$ and so we loose the factorial contribution). Let $Q$ be the number of intervals with $n_h > 0$. In the expression for $J_\tau(\Phi)(t)$ there are at least $Q$ factors of the form $e^{-(s-s')}$, $s' < s$, for which the sum of the quantities $s - s'$ is at least $Q$. So we have the bounds

$$\|J_\tau(\Phi)(t)\| \leqslant \|\!|\Phi|\!\|^n |\lambda|^{(2/3)(n-1)} N^n \sum_{Q \geqslant 1} e^{-Q/2} \prod_{h: n_h > 0} \frac{C^{n_h+1}}{(n_h!)^{1/8}}$$

$$\leqslant \|\!|\Phi|\!\|^n \frac{C^n |\lambda|^{(2/3)(n-1)}}{(n!)^{1/8}} \sum_{Q \geqslant 1} e^{-Q/2} Q^{n/8}$$

where we used that

$$\prod_{h: n_h > 0} \frac{1}{n_h!} \leqslant \frac{1}{n!} \sum_{\substack{k_1, \ldots, k_Q \geqslant 1 \\ k_1 + \cdots + k_Q = n}} \frac{n!}{k_1! \cdots k_Q!} = \frac{1}{n!} Q^n.$$

Now, by Jensens' inequality

$$\sum_{Q \geqslant 1} e^{-Q/2} Q^{n/8} \lesssim \left[\sum_{Q \geqslant 1} e^{-Q/2} Q^n\right]^{1/8} \lesssim c^n (n!)^{1/8}$$

and we obtain that

$$\|J_\tau(\phi)(t)\| \lesssim \|\!|\Phi|\!\|^n |\lambda|^{(2/3)(n-1)} C^n \quad (92)$$

possibly with a different value for $C$. So provided $|\lambda|$ is small enough (depending on $N$, $|\!|\!|\Phi|\!|\!|$) we conclude, from the representation (90), the uniform bound

$$\sup_{t\geqslant 0}\|\Psi_t\| \lesssim (1+|\!|\!|\Phi|\!|\!|)^{4N-1} + \sum_{n\geqslant 4N} C^n D^n |\!|\!|\Phi|\!|\!|^n |\lambda|^{(2/3)(n-1)} < \infty.$$

## Appendix B  Besov spaces of Banach algebras

In this section we want to recall some results about Besov spaces of functions (or distributions) from $\mathbb{R}^d$ taking values in a Banach algebra $\mathcal{A}$. All the results of this section can be found in [9, 10] for the theory of Besov spaces taking values in a Banach space and [175] for weighted Besov spaces. We present here only the $\mathbb{R}^d$ case. The definition of Besov spaces on $\mathbb{T}^d$ is similar, and all the results announced here also hold on $\mathbb{T}^d$.

We denote by $\mathcal{S}(\mathbb{R}^d)$ the set of smooth functions $f \in \mathcal{S}(\mathbb{R}^d)$ such that

$$\|f\|_{\ell,\alpha} := \|(1+|x|)^\ell |D^\alpha f|\|_{L^\infty(\mathbb{R}^d)} < +\infty$$

where $\ell \in \mathbb{R}_+$ and $\alpha \in \mathbb{N}^d$. We denote by $\mathcal{S}'(\mathbb{R}^d)$ the strong dual of $\mathcal{S}(\mathbb{R}^d)$ (equipped with the topology of the semi-norms $\|\cdot\|_{\ell,\alpha}$) and by $\mathcal{O}_M(\mathbb{R}^d)$ the set of smooth functions $f \in \mathcal{O}_M(\mathbb{R}^d)$ such that for any $\alpha \in \mathbb{N}^d$ there exists $\ell_\alpha \in \mathbb{R}_+$ for which

$$\|(1+|x|)^{-\ell_\alpha} |D^\alpha f|\|_{L^\infty(\mathbb{R}^d)} < +\infty.$$

We use the notations

$$\mathcal{S}(\mathbb{R}^d, \mathcal{A}) = \mathcal{S}(\mathbb{R}^d) \hat{\otimes} \mathcal{A}, \quad \mathcal{S}'(\mathbb{R}^d, \mathcal{A}) = \mathcal{S}'(\mathbb{R}^d) \hat{\otimes} \mathcal{A}, \quad \mathcal{O}_M(\mathbb{R}^d, \mathcal{A}) = \mathcal{O}_M(\mathbb{R}^d) \hat{\otimes} \mathcal{A},$$

where $\mathcal{A}$ is a Banach algebra and $A \hat{\otimes} B$ is the completion with respect to the natural metric of the algebraic tensor product $A \otimes B$ when $A$ is a Fréchet nuclear space and $B$ is a Banach space, see Remark 43. Note that these tensor products are well defined since the spaces $\mathcal{S}(\mathbb{R}^d), S'(\mathbb{R}^d), \mathcal{O}(\mathbb{R}^d)$ are nuclear. It is important to note that

$$\mathcal{S}(\mathbb{R}^d, \mathcal{A}) \xhookrightarrow{d} \mathcal{O}_M(\mathbb{R}^d, \mathcal{A}) \xhookrightarrow{d} \mathcal{S}'(\mathbb{R}^d, \mathcal{A})$$

where the arrows mean that each space is continuously embedded and dense in the following one.

**Theorem 77.** *It is possible to define uniquely* $\langle\cdot,\cdot\rangle \colon \mathcal{S}(\mathbb{R}^d, \mathcal{A}) \times \mathcal{S}'(\mathbb{R}^d, \mathcal{A}) \to \mathcal{A}$, $\cdot \colon \mathcal{S}(\mathbb{R}^d, \mathcal{A}) \times \mathcal{S}'(\mathbb{R}^d, \mathcal{A}) \to \mathcal{S}'(\mathbb{R}^d, \mathcal{A})$, $\cdot \colon \mathcal{S}(\mathbb{R}^d) \times \mathcal{S}'(\mathbb{R}^d, \mathcal{A}) \to \mathcal{S}'(\mathbb{R}^d, \mathcal{A})$, $* \colon \mathcal{S}(\mathbb{R}^d) \times \mathcal{S}'(\mathbb{R}^d, \mathcal{A}) \to \mathcal{O}_M(\mathbb{R}^d, \mathcal{A})$ *and* $D^\alpha \colon \mathcal{S}'(\mathbb{R}^d, \mathcal{A}) \to \mathcal{S}'(\mathbb{R}^d, \mathcal{A})$ *(where $\alpha \in \mathbb{N}^d$) which extend in an epicontinuous way the following operations: any $f \in \mathcal{S}(\mathbb{R}^d)$, $u \in \mathcal{S}'(\mathbb{R}^d)$ and $a_1, a_2 \in \mathcal{A}$ we have*

$$\begin{aligned}
\langle f \otimes a_1, u \otimes a_2 \rangle &= \langle f, u \rangle a_1 a_2 \\
(f \otimes a_1) \cdot (u \otimes a_2) &= (fu) \otimes (a_1 a_2) \\
f \cdot (u \otimes a_2) &= (fu) \otimes a_1 \\
f * (u \otimes a_2) &= (f*u) \otimes a_1 \\
D^\alpha(u \otimes a_1) &= (D^\alpha u) \otimes a_1
\end{aligned}$$

*where $\langle f, u \rangle$ is the normal pairing in $\mathcal{S}(\mathbb{R}^d) \times \mathcal{S}'(\mathbb{R}^d)$, $(fu)$ is the product in $\mathcal{S}(\mathbb{R}^d) \times \mathcal{S}'(\mathbb{R}^d)$, $(f*u)$ is the convolution in $\mathcal{S}(\mathbb{R}^d) \times \mathcal{S}'(\mathbb{R}^d)$ and $D^\alpha$ is the $\alpha$ derivative in $\mathcal{S}'(\mathbb{R}^d)$.*

**Proof.** The proof of this theorem can be found in [10] Appendix 1. $\square$

We recall the definition of Littlewood–Paley blocks: let $\chi, \varphi$ be smooth non-negative functions from $\mathbb{R}^n$ into $\mathbb{R}$ such that

- $\text{supp}(\chi) \subset B_{\frac{4}{3}}(0)$ and $\text{supp}(\varphi) \subset B_{\frac{8}{3}}(0) \setminus B_{\frac{3}{4}}(0)$,
- $\chi, \varphi \leqslant 1$ and $\chi(y) + \sum_{j \geq 0} \varphi(2^{-j} y) = 1$ for any $y \in \mathbb{R}^n$,
- $\text{supp}(\chi) \cap \text{supp}(\varphi(2^{-i} \cdot)) = \emptyset$ for $i \geqslant 1$,
- $\text{supp}(\varphi(2^{-j} \cdot)) \cap \text{supp}(\varphi(2^{-i} \cdot)) = \emptyset$ if $|i - j| > 1$,

where by $B_r(x)$ we denote the ball of center $x \in \mathbb{R}^n$ and radius $r > 0$. We introduce the following notations: $\varphi_{-1} = \chi$, $\varphi_j(\cdot) = \varphi(2^{-j} \cdot)$, $K_j = \hat{\varphi}_j$ and, for any $\ell > 0$, we define $\rho_\ell(x) := (1 + |x|^2)^{-\ell/2}$.

**Remark 78.** It is always possible to choose $\chi$ and $\varphi$ in such a way that there exists some constants $\gamma_{-1}, \gamma > 0$ and $0 < \theta_{-1}, \theta < 1$ such that

$$|K_{-1}(y)| \lesssim \exp(-\gamma_{-1} |y|^{\theta_{-1}}) \quad \text{and} \quad |K_0(y)| \lesssim \exp(-\gamma |y|^\theta),$$

(see, i.e., [155] Section 1.2.2 Proposition 1 or [124]).

If $v \in \mathcal{S}'(\mathbb{R}^d, \mathcal{A})$ and if $i \in \mathbb{Z}$, $i \geqslant -1$ we define $i$-th Littlewood–Paley block as follows

$$\Delta_i v = K_i * v \in \mathcal{O}_M(\mathbb{R}^d, \mathcal{A}).$$

Then, if $s \in \mathbb{R}$, $p, q \in [1, +\infty]$ and $\ell \in \mathbb{R}$ we define the function

$$\|v\|_{B^s_{p,q,\ell}(\mathbb{R}^d, \mathcal{A})} = \left( \sum_{j=-1}^{+\infty} 2^{jsq} \|\Delta_j v\|^q_{L^p_\ell(\mathbb{R}^d, \mathcal{A})} \right)^{1/q},$$

when $q \in [1, +\infty)$ and $\|v\|_{B^s_{p,+\infty,\ell}(\mathbb{R}^d, \mathcal{A})} = \sup_j (2^{js} \|\Delta_j v\|_{L^p_\ell(\mathbb{R}^d, \mathcal{A})})$, where $\|\cdot\|_{L^p_\ell(\mathbb{R}^d, \mathcal{A})}$ is the norm in the space $L^p_\ell(\mathbb{R}^d, \mathcal{A})$ that is

$$\|f\|_{L^p_\ell(\mathbb{R}^d, \mathcal{A})} = \left( \int_{\mathbb{R}^n} (\|f(y)\|_\mathcal{A} \rho_\ell(y))^p \mathrm{d}y \right)^{1/p},$$

if $p \in [1, +\infty)$ and

$$\|f\|_{L^\infty_\ell(\mathbb{R}^d, \mathcal{A})} = \sup_{y \in \mathbb{R}^d} (\|f(y)\|_\mathcal{A} \rho_\ell(y)),$$

if $p = +\infty$. For any $v \in \mathcal{S}(\mathbb{R}^d, \mathcal{A})$ the norm $\|v\|_{B^s_{p,q,\ell}(\mathbb{R}^d, \mathcal{A})} < +\infty$ is finite. Then we look at $B^s_{p,q,\ell}(\mathbb{R}^d, \mathcal{A})$ as the closure of $\mathcal{S}(\mathbb{R}^d, \mathcal{A})$ in $\mathcal{S}'(\mathbb{R}^d, \mathcal{A})$ with respect to the norm $\|\cdot\|_{B^s_{p,q,\ell}(\mathbb{R}^d, \mathcal{A})}$. Hereafter, if $s \in \mathbb{R}$, $p, q \in [1, +\infty]$ and $\ell \in \mathbb{R}$, we use the following notations $\mathcal{C}^s_\ell(\mathbb{R}^d, \mathcal{A}) := B^s_{\infty,\infty,\ell}(\mathbb{R}^d, \mathcal{A})$, $B^s_{p,q}(\mathbb{R}^d, \mathcal{A}) = B^s_{p,q,0}(\mathbb{R}^d, \mathcal{A})$, $\mathcal{C}^s(\mathbb{R}^d, \mathcal{A}) := B^s_{\infty,\infty,0}(\mathbb{R}^d, \mathcal{A})$, $B^s_{p,q,\ell} := B^s_{p,q,\ell}(\mathbb{R}^d, \mathbb{R})$ etc.

In this paper we need only the next three results.

**Theorem 79.** *For any $s, \ell \in \mathbb{R}$ and $v \in \mathcal{C}^s_\ell(\mathbb{R}^d, \mathcal{A})$ we have that $\rho_\ell v \in \mathcal{C}^s(\mathbb{R}^d, \mathcal{A})$ and*

$$\|v\|_{\mathcal{C}^s_\ell(\mathbb{R}^d, \mathcal{A})} \sim \|\rho_\ell v\|_{\mathcal{C}^s(\mathbb{R}^d, \mathcal{A})}.$$

**Proof.** The estimate holds by Theorem 6.5 in [175] for weighted Besov spaces on $\mathbb{R}^d$ with polynomial-like weights. The theorem for a general Banach algebra $\mathcal{A}$ is a straightforward generalization. $\square$

**Theorem 80.** *Consider $m > 0$, $\alpha, s, \ell \in \mathbb{R}$ such that $s, s + \alpha \notin \mathbb{N}$ then we have that $(-\Delta + m)^{-\alpha}$, where $\Delta$ is the standard Laplacian on $\mathbb{R}^d$, is a continuous linear map from $\mathcal{C}^s_\ell(\mathbb{R}^d, \mathcal{A})$ into $\mathcal{C}^{s+\alpha}_\ell(\mathbb{R}^d, \mathcal{A})$*

**Proof.** This is exactly Theorem 5.3.2 of [10] for the unweighted case $\ell = 0$. The theorem can be easily extended using the techniques of Chapter 6 of [175] to Besov spaces with weights $\rho_\ell$. □

**Theorem 81.** *Let $\alpha > 0$ and consider $\mathfrak{a} \in B^\alpha_{1,\infty}$, then we have that the convolution $\mathfrak{a}*$ can be extended in a unique continuous way from an operator from $\mathcal{S}(\mathbb{R}^d, \mathcal{A})$ into itself, into an operator from $\mathcal{C}^s(\mathbb{R}^d, \mathcal{A})$ into $\mathcal{C}^{s+\alpha}(\mathbb{R}^d, \mathcal{A})$.*

**Proof.** Using the notations of [9], if $\mathfrak{a} \in B^\alpha_{1,\infty}$ then $((-\Delta + 1)^{\alpha/2}(\mathfrak{a})) \in \mathcal{F}L_1$, then the thesis follows from Corollary 6.4 in [9]. □

**Theorem 82.** *Let $s > 0$ and $\ell_1, \ell_2 \in \mathbb{R}$, then we have that the product $\cdot : \mathcal{S}(\mathbb{R}^d, \mathcal{A}) \times \mathcal{S}(\mathbb{R}^d, \mathcal{A}) \to \mathcal{S}(\mathbb{R}^d, \mathcal{A})$ can be (uniquely) extended in a continuous way from $\mathcal{C}^s_{\ell_1}(\mathbb{R}^d, \mathcal{A}) \times \mathcal{C}^s(\mathbb{R}^d, \mathcal{A})$ into $\mathcal{C}^s_{\ell_1 + \ell_2}(\mathbb{R}^d, \mathcal{A})$, furthermore for any $v_1 \in \mathcal{C}^s_{\ell_1}(\mathbb{R}^d, \mathcal{A})$ and $v_2 \in \mathcal{C}^s_{\ell_2}(\mathbb{R}^d, \mathcal{A})$ we have*

$$\|v_1 \cdot v_2\|_{\mathcal{C}^s_{\ell_1 + \ell_2}(\mathbb{R}^d, \mathcal{A})} \lesssim \|v_1\|_{\mathcal{C}^s_{\ell_1}(\mathbb{R}^d, \mathcal{A})} \|v_2\|_{\mathcal{C}^s_{\ell_2}(\mathbb{R}^d, \mathcal{A})},$$

**Proof.** The proof can be found in [10] in the unweighted case. The proof in the Besov spaces with polynomial weights $\rho_\ell$ is similar and can be done using the techniques of Chapter 6 of [175]. □

# Appendix C Convergence of the perturbative series: infinite dimensional case

In this appendix we want to consider the stationary solution to the equation

$$\tilde{\Psi}_{\varepsilon, \lambda^2, t}(x) = \lambda^2 \int_{-\infty}^t e^{(\Delta - m_f^2)(t - \tau)}(F_{\varepsilon, Y}(\tilde{\Psi}_\tau(\cdot))) \mathrm{d}\tau + \tilde{B}^A_{\varepsilon, t} \tag{93}$$

as a function of $\lambda^2$, i.e. we want to study the regularity of the function $\tilde{\lambda} \longmapsto Z_{F_{\varepsilon, Y}, \tilde{B}^A_{\varepsilon, t}}(\tilde{\lambda})$ from a neighborhood $U \subset \mathbb{R}$ of 0 into $C^0(\mathbb{R}, \mathcal{C}^q(\mathbb{R}^2, \mathcal{A} \otimes \mathbb{C}^2))$ (here $0 < q < \gamma$ where $\gamma$ is the constant chosen in Remark 59), where

$$Z_{F_{\varepsilon, Y}, \tilde{B}^A_{\varepsilon, t}}(\tilde{\lambda}) := \tilde{\Psi}_{\varepsilon, \tilde{\lambda}, t}(x) \in C^0(\mathbb{R}, \mathcal{C}^q(\mathbb{R}^2, \mathcal{A} \otimes \mathbb{C}^4)).$$

We will assume that the hypotheses of Theorem 60 are satisfied all along this appendix.

**Theorem 83.** *The function $Z_{F_{\varepsilon, Y}, \tilde{B}^A_t}(\tilde{\lambda})$ is real analytic (see Definition 91 below) in a neighborhood of 0.*

An easy consequence is the following corollary on the infinite volume Schwinger functions

$$\mathcal{S}_n(f_1, ..., f_n | \lambda) := \lim_{h_n \to 1} \frac{\omega\left(\tilde{\psi}(f_1) \cdots \tilde{\psi}(f_k) e^{\mathcal{V}_{\varepsilon, h_n}(\tilde{\psi})}\right)}{\omega\left(e^{\mathcal{V}_{\varepsilon, h_n}(\tilde{\psi})}\right)} = \omega(\tilde{\Psi}^\mathfrak{s}_{\varepsilon, \lambda^2, t}(f_1) \cdots \tilde{\Psi}^\mathfrak{s}_{\varepsilon, \lambda^2, t}(f_k))$$

where $h_n$ is a sequence of spatial cut-off as defined in Section 5 and $f_1,...,f_n \in C_c^\infty(\mathbb{R}^2) \otimes \mathbb{C}^4/$

**Corollary 84.** *The Schwinger functions $\mathcal{S}_n(f_1,...,f_n|\lambda)$ of the (regularized) Yukawa$_2$ model are real analytic functions of the coupling constant $\lambda$.*

**Proof.** We have that the functional $(Y_1,...,Y_n) \mapsto \omega(Y_1(f_1)\cdots Y_n(f_n))$ (where $(Y_1,...,Y_n \in \mathcal{C}^q(\mathbb{R}^2, \mathcal{A} \otimes \mathbb{C}^4))$ is a multilinear continuous functional (and so real analytic), from $\mathcal{C}^q(\mathbb{R}^2, \mathcal{A} \otimes \mathbb{C}^4)$ into $\mathbb{C}$ (since $\omega$ is continuous and linear from $\mathcal{A}$ into $\mathbb{C}$ and the product in $\mathcal{A}$ is also continuous and multilinear). Since the composition of real analytic functionals is real analytic the thesis follows. $\square$

**Remark 85.** Theorem 83 and Corollary 84 do not depend on the fact that the non-linearity $F_{\varepsilon,h,Y}$ has the specific form of the case of Yukawa interaction (i.e., for example, that is a third degree polynomial or on the specific non local form). In particular Theorem 83 and Corollary 84 extend easily to other (UV regularized) Euclidean fermionic models.

The rest of this appendix contains the proof of Theorem 83. Theorem 83 can be seen as an infinite dimensional generalization of the results in Appendix A, where, between the other things, it is proved that the real analyticity of Theorem 83 holds for finite dimensional systems. In our case we can consider the finite dimensional approximation

$$\tilde{\Psi}^{h_R}_{\varepsilon,\tilde{\lambda},N,R,t} = \tilde{\lambda} \int_{-\infty}^t e^{(\Delta - m_f^2)(t-\tau)} P_N(F_{\varepsilon,h_R,Y}(\tilde{\Psi}^{h_R}_{\varepsilon,R,\tau})) \mathrm{d}\tau + \tilde{B}^A_{\varepsilon,R,t}(v) \tag{94}$$

and the relative function $Z_{P_N(F_{\varepsilon,h_R,Y}),\tilde{B}^A_{\varepsilon,R,t}}(\tilde{\lambda}) := \tilde{\Psi}^{h_R}_{\varepsilon,\tilde{\lambda},N,R,t}$.

**Lemma 86.** *The function $\tilde{\lambda} \mapsto Z_{P_N(F_{\varepsilon,h_R,Y}),\tilde{B}^A_{\varepsilon,R,t}}(\tilde{\lambda})$ is real analytic in a neighborhood of $0$*

**Proof.** The result follows from the proof of convergence of perturbation series given in Appendix A (in particular from inequalities (91) and (92)) for finite dimensional Grassmannian equation and from the fact that the topology on $\mathcal{G}(P_N(C^\infty(\mathbb{T}_R^2, \mathcal{A} \otimes \mathbb{C}^4)))$ does depend on the norm chosen on $P_N(C^\infty(\mathbb{T}_R^2, \mathcal{A} \otimes \mathbb{C}^4))$ since all the norms on $P_N(C^\infty(\mathbb{T}_R^2, \mathcal{A} \otimes \mathbb{C}^4))$ are equivalent, this space being finite dimensional. $\square$

Instead of proving the Theorem 83 by trying to generalize the computation done in Appendix A we present a proof using the approximation method described in Section 5, since it is closer to the SPDE approach followed in the main part of the present paper. For this reason we want to introduce some sort of a-priori estimate characterizing real analytic functions. This is done exploiting the notion of formal power series and majorants method (see, e.g., [176] for an introduction to the subject). We call a formal power series around $0$ taking values in the Banach space $X$ the following object:

$$f = \sum_{n=0}^{+\infty} f^{(n)} \lambda^n, \quad \lambda \in \mathbb{C}, \tag{95}$$

where $f^{(n)} \in X$. We write $f(0) = f^{(0)}$ and also

$$\frac{\mathrm{d}^n f}{\mathrm{d}\lambda^n} = \sum_{k=n}^{+\infty} f^{(n)} \lambda^{n-k} \prod_{i=0}^{n-1} (k-i).$$

When the power series (95) is absolutely convergent for $0 < |\lambda| < \varepsilon$ (for some $\varepsilon > 0$) we identify it with the real analytic function $\lambda \mapsto f(\lambda) := \sum_{n=0}^{+\infty} f^{(n)} \lambda^n \in \mathbb{C}$ defined for $|\lambda| < \varepsilon$.

**Definition 87.** *Let $f_X$ be a formal power series around $0$ taking values in the Banach space $X$ and let $K$ be a formal power series around $0$ taking values in $\mathbb{C}$. We say that $f$ is majored by $K$, and we write $f_X \trianglelefteq K$, if for any $n \in \mathbb{N}$ we have*

$$\left\|\frac{d^n f_X}{d\lambda^n}(0)\right\|_X = n!\|f^{(n)}\|_X \leqslant \frac{d^n K}{d\lambda^n}(0) = n! K^{(n)},$$

*where $K(\lambda) = \sum_{n=0}^{+\infty} K^{(n)} \lambda^n$.*

**Remark 88.** *If $f_\mathbb{C}$ is a formal power series in $\mathbb{C}$ and $K$ is an analytic function such that $f_\mathbb{C} \trianglelefteq K$ then also $f_\mathbb{C}$ is an analytic function. Furthermore if the power series defining $g$ (in a neighborhood of $0$) converges absolutely for $\lambda = R > 0$, then the power series defining $f_\mathbb{C}$ converges absolutely for $\lambda \in \mathbb{C}$ with $|\lambda| \leqslant R$ and furthermore*

$$\sup_{\lambda \in \mathbb{C}, |\lambda| \leqslant R} |f_\mathbb{C}(\lambda)| \leqslant K(R).$$

In general it is possible to define the sum and product of formal power series. Furthermore if $g \colon \mathbb{C} \to \mathbb{C}$ is an entire function and $f_\mathbb{C}$ is a formal power series taking values in $\mathbb{C}$ it is possible to define the composition $g \circ f_\mathbb{C}$ as a formal power series taking values in $0$.

Let $g$ be an analytic function such that $0 \trianglelefteq g$, i.e. $g(\lambda) := \sum_{n=0}^{+\infty} g^{(n)} \lambda^n$ such that $g^{(n)} \geqslant 0$, consider $c \geqslant 0$ and define

$$G_{g,c}(a, \lambda) = \lambda g(a) + c.$$

**Lemma 89.** *For any entire function $g$ such that $0 \trianglelefteq g$ and $c \geqslant 0$ there exists a unique analytic function $K \colon U \to \mathbb{R}$ (where $U$ is a suitable neighborhood of $0$) such that*

$$K(\lambda) = G_{g,c}(K(\lambda), \lambda) \tag{96}$$

*for any $\lambda \in U$. Furthermore if $f_\mathbb{C}$ is a formal power series such that $f_\mathbb{C} \trianglelefteq G_{g,c}(f_\mathbb{C}(\cdot), \cdot)$ then $f_\mathbb{C} \trianglelefteq K$.*

**Proof.** Equation (96) admits a unique real analytic solution $K(\lambda)$ defined in a neighborhood $U$ of $0$ and taking values in a neighborhood $V$ of $c \in \mathbb{R}$. Indeed we have that the point $(c, 0) = (k, \lambda)$ is solution to the equation

$$k - G_{g,c}(k, \lambda) = 0. \tag{97}$$

Furthermore we have that

$$\partial_k (k - G_{g,c}(k, \lambda))|_{(k,\lambda)=(c,0)} = (1 - \lambda g'(k))|_{(k,\lambda)=(c,0)} = 1.$$

Thus by the implicit function theorem for holomorphic functions (see, e.g., [116] Chapter 3, Theorem 3.11), there is are two sets $U, V \subset \mathbb{R}$ (such that $V \times U$ is a neighborhood of $(c, 0)$) and a unique real analytic function $K \colon U \to V$ solution to equation (96). Furthermore this $K$ is the unique solution to equation (96) of the form $K \colon U \to \mathbb{R}$ since $(c, 0)$ is the unique solution to equation (97) of the form $(\cdot, 0) \in \mathbb{R}^2$.

The fact that $K(\lambda) = \sum_{n=0}^{+\infty} K^{(n)} \lambda^n$ is the unique solution to equation (96) in a neighborhood of $0$ implies that

$$K^{(1)} = \frac{d}{d\lambda}(\lambda g(K(\lambda)))\bigg|_{\lambda=0} = g(K(0)) = \sum_{k=0}^{+\infty} g^{(k)} (K^{(0)})^k$$

and also, for any $n \geqslant 2$,

$$\begin{aligned} K^{(n)} &= \frac{1}{n!}\frac{\mathrm{d}^n}{\mathrm{d}\lambda^n}(\lambda g(K(\lambda)))\bigg|_{\lambda=0} \\ &= \frac{n}{n!}\frac{\mathrm{d}^{n-1}}{\mathrm{d}\lambda^{n-1}}g(K(\lambda))\bigg|_{\lambda=0} \\ &= \sum_{k=1}^{+\infty} \sum_{\ell_1+\cdots+\ell_k=n-1} g^{(k)} K^{(\ell_1)}\cdots K^{(\ell_k)}. \end{aligned} \quad (98)$$

Suppose that $f_{\mathbb{C}} \trianglelefteq G_{g,c}(f_{\mathbb{C}}(\cdot),\cdot)$, with $f_{\mathbb{C}}(\lambda) = \sum_{n=0}^{+\infty} f^{(n)}\lambda^n$, then we want to prove that $|f^{(n)}| \leqslant K^{(n)}$ by induction on $n$. For $n=0$ we have

$$|f^{(0)}| \leqslant G_{g,c}(f(0),0) = c = K^{(0)}.$$

For $n=1$ we have

$$|f^{(1)}| \leqslant \frac{\mathrm{d}}{\mathrm{d}\lambda}(G_{g,c}(f_{\mathbb{C}}(\lambda),\lambda)) = g(f_{\mathbb{C}}(0)) = \sum_{k=0}^{+\infty} g^{(k)}(f^{(0)})^k \leqslant \sum_{k=0}^{+\infty} g^{(k)}(K^{(0)})^k = K^{(1)}.$$

Supposing that $|f^{(k)}| \leqslant K^{(k)}$ for any $k < n$, where $n \geqslant 2$, then we will have

$$\begin{aligned} |f^{(n)}| &\leqslant \frac{1}{n!}\frac{\mathrm{d}^n}{\mathrm{d}\lambda^n}G_{g,c}(f_{\mathbb{C}}(\lambda),\lambda)\bigg|_{\lambda=0} \\ &\leqslant \frac{n}{n!}\frac{\mathrm{d}^{n-1}}{\mathrm{d}\lambda^{n-1}}g(f_{\mathbb{C}}(\lambda))\bigg|_{\lambda=0} \\ &\leqslant \sum_{k=1}^{+\infty} \sum_{\ell_1+\cdots+\ell_k=n-1} g^{(k)} f^{(\ell_1)}\cdots f^{(\ell_k)} \\ &\leqslant \sum_{k=1}^{+\infty} \sum_{\ell_1+\cdots+\ell_k=n-1} g^{(k)} K^{(\ell_1)}\cdots K^{(\ell_k)} \\ &\leqslant K^{(n)} \end{aligned}$$

where we used equality (98) and the induction hypothesis. This proves that $f_{\mathbb{C}} \trianglelefteq K$. $\square$

**Remark 90.** From formula (98) it is clear that the map $(g,c) \mapsto K$ is monotone, i.e. if $K_1(\lambda) = G_{g_1,c_1}(K_1(\lambda),\lambda)$, $K_2(\lambda) = G_{g_2,c_2}(K_2(\lambda),\lambda)$ and $g_1 \trianglelefteq g_2$ and $c_1 \leqslant c_2$ then $K_1 \trianglelefteq K_2$.

We want to use Lemma 89 for giving an a-priori bound on the norm of $\frac{\partial^n \tilde{\Psi}_{\varepsilon,\lambda^2,t}}{\partial \lambda^n}$ (with $\tilde{\Psi}^{\mathfrak{s}}_{\varepsilon,\lambda^2,t}$ as in equation (93)) when $\lambda = 0$. This would imply that the $\tilde{\Psi}_{\varepsilon,\lambda^2,t}$ is real analytic in $\lambda$ when $\lambda$ is small enough.

First of all we introduce the notion of real-analyticity for functions taking values in a Banach space.

**Definition 91.** *Let $X$ be a Banach space and let $Z: U \to X$ be an infinitely differentiable function (where $U \subset \mathbb{R}$ is an open set). We say that $Z$ is real analytic in $U$ if for any $\lambda_0 \in U$ we have*

$$Z(\lambda) = \sum_{k=0}^{+\infty} \frac{1}{n!}\frac{\mathrm{d}^n Z}{\mathrm{d}\lambda^n}(\lambda_0)(\lambda - \lambda_0)^n \quad (99)$$

for $|\lambda - \lambda_0| > 0$ *small enough (where the series on the right hand side of (99) is supposed to converge absolutely in X).*

If $Z\colon U \to X$ is differentiable an infinite number of times at 0 we can define the formal power series

$$\tilde{Z}(\lambda) := \sum_{n=0}^{+\infty} Z^{(n)} \lambda^n \qquad (100)$$

where $Z^{(n)} \in X$ is given by

$$Z^{(n)} := \frac{1}{n!} \frac{\mathrm{d}^n Z}{\mathrm{d}\lambda^n}(0) \in X.$$

We can consider also the formal power series $S_Z(\lambda)$ defined as

$$S_Z(\lambda) := \sum_{n=0}^{+\infty} \|Z^{(n)}\| \lambda^n. \qquad (101)$$

It is useful to recall an equivalent characterization of real analytic functions in Banach spaces.

**Remark 92.** We have that $\tilde{Z}$ (as given by (100)) is convergent in a neighborhood of 0 (and so it defines an analytic function in a neighborhood of 0) if and only if $S_Z$ (as given by (101)) is real analytic. Furthermore for any $0 \trianglelefteq g$, we have $\tilde{Z} \trianglelefteq g$ if and only if $S_Z \trianglelefteq g$.

It is important to note that if $\tilde{Z}$ is convergent (and so it defines an analytic function in a neighborhood of 0) this is not enough for proving that $Z$ is also real analytic in a neighborhood of 0. Indeed, in general, if $\tilde{Z}$ is real analytic we could have $\tilde{Z}(\lambda) \neq Z(\lambda)$ for $\lambda \neq 0$.

In order to overcome this problem we introduce the following useful theorems.

**Theorem 93.** *Let $X$ be a Banach space and let $Z\colon U \subset \mathbb{R} \to X$ be a function defined in a neighborhood of 0. Then $Z$ is real-analytic in $U$ if and only if, for any $\ell^* \in X^*$, the function $\lambda \mapsto \ell^*(Z(\lambda))$ is real analytic in $U$.*

**Proof.** The proof can be found in [39] Proposition 1 and Remark 2. □

**Theorem 94.** *Suppose that $Z_n(\lambda)$ is a sequence of functions from an open neighborhood $U \subset \mathbb{R}$ of 0 into a Banach space $X$, converging point-wise $Z^n(\lambda) \to Z(\lambda)$ to some function from $U$ into $X$. Suppose furthermore that each $Z_n(\lambda)$ is real analytic and that there is a real analytic function $0 \trianglelefteq K\colon U \to \mathbb{R}$ such that $Z_n \trianglelefteq K$, then $Z(\lambda)$ is real analytic in a neighborhood of 0 and $Z \trianglelefteq K$.*

**Proof.** Consider $\ell^* \in X^*$ such that $\|\ell^*\|_{X^*} \leqslant 1$ and define $\ell^*_{Z_n}(\lambda) := \ell^*(Z_n(\lambda))$. Since $Z_n$ are real analytic also $\ell^*_{Z_n}$ is real analytic, and also $\ell^*_{Z_n}(\lambda) \to \ell^*_Z(\lambda) := \ell^*(Z(\lambda))$ for any $\lambda \in U$. Furthermore $\ell^*_{Z_n} \trianglelefteq S_{Z_n} \trianglelefteq K$ which implies that, by Remark 88, for any $R$ such that $K(R)$ is well defined, we have

$$\sup_{\lambda \in \mathbb{C}, |\lambda| \leqslant R} |\ell^*_{Z_n}(\lambda)| \leqslant K(R),$$

i.e. $\ell^*_{Z_n}$ are uniformly bounded in a (complex)-neighborhood of 0. This means by Vitali-Porter theorem (see, i.e. [154] Chapter 2, Section 2.4) that also $\ell^*_Z(\lambda)$ is real analytic for $\lambda \in U \cap \{|\lambda| < R\}$. Since $\ell^*_Z$ is real analytic for any $\ell^* \in X^*$ such that $\|\ell^*\|_{X^*} = 1$, the function $\tilde{\ell}(Z(\lambda))$ is real analytic in $U \cap \{|\lambda| < R\}$ which implies, by Theorem 93, that $Z$ is real analytic in $U \cap \{|\lambda| < R\}$ which is a non-empty neighborhood of 0.

Finally we have

$$\left|\frac{d^k S_Z}{d\lambda^k}(0)\right| = \sup_{\ell^* \in X^*, \|\ell^*\|_{X^*}=1} \left|\frac{d^k \ell_Z^*}{d\lambda^k}(0)\right| \leqslant \sup_{\ell^* \in X^*, \|\ell^*\|_{X^*}=1, n \in \mathbb{N}} \left|\frac{d^k \ell_{Z_n}^*}{d\lambda^k}(0)\right| \leqslant K^{(k)}$$

where we used that, from the Vitali-Porter theorem, $\frac{d^k \ell_{Z_n}^*}{d\lambda^k}$ converges to $\frac{d^k \ell_Z^*}{d\lambda^k}$ uniformly. $\square$

**Remark 95.** It is easy to generalize the proof of Theorem 94 replacing the convergence of $Z_n(\lambda)$ to $Z(\lambda) \in X$ in the norm of $X$, by the convergence of $Z_n(\lambda)$ to $Z(\lambda) \in X$ with respect to a (weaker) norm of a Banach space $X' \supset X$ in which $X$ is densely contained.

Let $P: X \to X$ be a polynomial map, i.e. there are some $P_0, P_1, ..., P_n$ such that $P_k: X^k \to X$ are continuous $k$-linear symmetric maps such that

$$P(x) := \sum_{k=0}^n P_k(x, ..., x), \quad x \in X.$$

We also suppose that there is $\mathfrak{P}: \mathbb{R} \to \mathbb{R}$ defined as

$$\mathfrak{P}(\lambda) := \sum_{k=0}^n \|P_k\|_{\mathcal{L}(X^k, X)} \lambda^k, \quad \text{for } \lambda \in \mathbb{C},$$

where $\|\cdot\|_{\mathcal{L}(X^k, X)}$ is the natural operator norm of $\mathcal{L}(X^k, X)$ that is the space of $k$-linear maps between $X^k$ into $X$. Suppose that $Z_{P,C}: U \to X$ (where $U \subset \mathbb{R}$ is a neighborhood of 0) is the unique solution to the equation

$$Z_{P,C}(\lambda) = \lambda P(Z(\lambda)) + C \tag{102}$$

where $\lambda \in U$ and $C \in X$. We are now ready to prove a general a-priori estimate for equations of the form (102).

**Theorem 96.** *Suppose that a solution $Z_{P,C}$ to equation (102) is differentiable an infinite number of time in $\lambda = 0$, then we have $Z_{P,C} \trianglelefteq K_{P,C}$ where $K_{P,C}$ is the unique analytic solution (in a neighborhood of 0) to the equation*

$$K_{P,C}(\lambda) = \lambda \mathfrak{P}(K_{P,C}(\lambda)) + \|C\|_X. \tag{103}$$

**Proof.** If we take the derivatives in the variable $\lambda$ in the point $\lambda = 0$ and write

$$Z_{P,C}^{(n)} := \frac{1}{n!} \frac{d^n Z_{P,C}}{d\lambda^n}(0) \in X,$$

we obtain from the formal series (100) and equation (102)

$$Z_{P,C}^{(0)} = C, \quad Z_{P,C}^{(1)} = P(Z^{(0)}) = P(C)$$

and, for $r \geqslant 2$,

$$Z_{P,C}^{(r)} = \sum_{k=0}^n \sum_{j_1 + \cdots + j_k = r-1} P_k(Z_{P,C}^{(j_1)}, ..., Z_{P,C}^{(j_k)}).$$

This implies that

$$\|Z_{P,C}^{(0)}\|_X = \|C\|_X, \quad \|Z_{P,C}^{(1)}\| \leqslant \mathfrak{P}(\|C\|_X)$$

and, for $r \geqslant 2$,

$$\|Z_{P,C}^{(r)}\| \leqslant \sum_{k=0}^n \sum_{j_1 + \cdots + j_{2k+1} = r-1} \|P_{2k+1}\|_{\mathcal{L}(X^{2k+1}, X)} \|Z_{P,C}^{(j_1)}\|_X \cdots \|Z_{P,C}^{(j_{2k+1})}\|_X. \tag{104}$$

If we set $S_{Z_{P,C}} = \sum_{n=0}^{+\infty} \|Z_{P,C}^{(n)}\| \lambda^n$, then we have

$$\begin{aligned} S_{Z_{P,C}}(\lambda) &= \|C\|_X + \sum_{r=1}^{+\infty} \|Z_{P,C}^{(k)}\| \lambda^r \\ &\trianglelefteq \|C\|_X + \sum_{r=1}^{+\infty} \lambda^r \sum_{k=0}^{n} \sum_{j_1+\cdots+j_k=r-1} \|P_k\|_{\mathcal{L}(X^k,X)} \|Z_{P,C}^{(j_1)}\|_X \cdots \|Z_{P,C}^{(j_k)}\|_X \\ &\trianglelefteq \|C\|_X + \lambda \sum_{k=0}^{n} \sum_{r=0}^{+\infty} \lambda^r \sum_{j_1+\cdots+j_k=r} \|P_k\|_{\mathcal{L}(X^k,X)} \|Z_{P,C}^{(j_1)}\|_X \cdots \|Z_{P,C}^{(j_k)}\|_X \\ &\trianglelefteq \|C\|_X + \lambda \mathfrak{P}(S_{Z_{P,C}}(\lambda)) \end{aligned}$$

where we used inequality (104) and the fact that, by Theorem 93, the function $S_{P,C}$ is real analytic in a neighborhood of 0 with positive derivatives in 0. By Lemma 89 the thesis follows by taking $G_{g,c}(a,\lambda) = \lambda \mathfrak{P}(a) + \|C\|_X$. □

**Proof of Theorem 83.** We want to apply the previous theorems for proving Theorem 83. Here the Banach space $X$ is $X = C^0(\mathbb{R}, \mathcal{C}^q(\mathbb{R}^2, \mathcal{A} \otimes \mathbb{C}^4))$ with the norm given by

$$\|\tilde{\Psi}\|_X = \sup_{t \in \mathbb{R}} \|\tilde{\Psi}\|_{\mathcal{C}^q(\mathbb{R}^2, \mathcal{A} \otimes \mathbb{C}^4)}.$$

The function $Z_{P_N(F_{\varepsilon,h_R,Y}), \tilde{B}_{\varepsilon,R,t}^A}(\tilde{\lambda})$ is real analytic in $\tilde{\lambda}$ in its radius of convergence. Furthermore it is a solution to equation (94), which is of the form (102) with

$$P_{3,h_R,N}(\tilde{\Psi})(t,\cdot) = \int_{-\infty}^{t} e^{(\Delta - m_f^2)(t-\tau)} P_N(F_{\varepsilon,h_R,Y}(\tilde{\Psi}(\tau,\cdot))) \mathrm{d}\tau, \quad (105)$$

$C = \tilde{B}_{\varepsilon,R,t}^A$, and $\mathfrak{P}_{h_R,N}(\tilde{\lambda}) = \|P_{3,h_R,N}\|_{\mathcal{L}(X^3,X)} \tilde{\lambda}^3$. This means that, by Theorem 96, we have

$$Z_{P_N(F_{\varepsilon,h_R,Y}), \tilde{B}_{\varepsilon,R,t}^A} \trianglelefteq K_{P_{3,h_R,N}, \tilde{B}_{\varepsilon,R}^A},$$

where $K_{P_{3,h,N}, \tilde{B}_{\varepsilon,R}^A}$ is the solution to the equation

$$K_{P_{3,h_R,N}, \tilde{B}_{\varepsilon,R}^A}(\tilde{\lambda}) = \tilde{\lambda} \|P_{3,h_R,N}\|_{\mathcal{L}(X^3,X)} \left( K_{P_{3,h_R,N}, \tilde{B}_{\varepsilon,R_X}^A}(\tilde{\lambda}) \right)^3 + \|\tilde{B}_{\varepsilon,R}^A\|_X.$$

Since $\|P_{3,h_R,N}\|_{\mathcal{L}(X^3,X)} \leqslant \|P_{3,h_R}\|_{\mathcal{L}(X^3,X)}$ (i.e. the polynomial $P_{3,h_R,N}$ without the finite dimensional projection), by Remark 90, we have $Z_{P_N(F_{\varepsilon,h_R,Y}), \tilde{B}_{\varepsilon,R,t}^A} \trianglelefteq K_{P_{3,h_R}, \tilde{B}_{\varepsilon,R}^A}$. Since $Z_{P_N(F_{\varepsilon,h_R,Y}), \tilde{B}_{\varepsilon,R,t}^A}$ converges point-wise to $Z_{F_{\varepsilon,h_R,Y}, \tilde{B}_{\varepsilon,R,t}^A}$ (i.e. the solution to equation (94) without the finite dimensional projection in the non-linearity) we have that, by Theorem 94, $Z_{F_{\varepsilon,h_R,Y}, \tilde{B}_{\varepsilon,R,t}^A}$ is real analytic and $Z_{F_{\varepsilon,h_R,Y}, \tilde{B}_{\varepsilon,R,t}^A} \trianglelefteq K_{P_{3,h_R}, \tilde{B}_{\varepsilon,R}^A}$.

Using again Remark 90 and the fact that $\|\tilde{B}_{\varepsilon,R,t}^A\| \leqslant \|\tilde{B}_{\varepsilon,t}^A\|_X$ (where, we recall, that $\tilde{B}_{\varepsilon,R,t}^A$ is the periodized version of $\tilde{B}_{\varepsilon,t}^A$), we get $Z_{F_{\varepsilon,h_R,Y}, \tilde{B}_{\varepsilon,R,t}^A} \trianglelefteq K_{P_{3,h_R}, \tilde{B}_{\varepsilon,R}^A} \trianglelefteq K_{P_{3,h}, \tilde{B}_{\varepsilon}^A}$. Furthermore, since $Z_{F_{\varepsilon,h_R,Y}, \tilde{B}_{\varepsilon,R,t}^A}$ converges to $Z_{F_{\varepsilon,h,Y}, \tilde{B}_{\varepsilon,t}^A}$ (where $Z_{F_{\varepsilon,h,Y}, \tilde{B}_{\varepsilon,t}^A}$ is the solution to equation (94) where $\tilde{B}_{\varepsilon,R,t}^A$ is replaced by the non-periodized version $\tilde{B}_{\varepsilon,t}^A$ of $\tilde{B}_{\varepsilon,R,t}^A$) in $X' = C^0(\mathbb{R}, \mathcal{C}_\ell^q(\mathbb{R}^2, \mathcal{A} \otimes \mathbb{C}^4)) \supset X$ (where we consider the weighted space $\mathcal{C}_\ell^q$ instead of the unweighted one), by Theorem 94 and Remark 95, we obtain that $Z_{F_{\varepsilon,h,Y}, \tilde{B}_{\varepsilon,t}^A}$ is real analytic with respect to $\tilde{\lambda}$ and also $Z_{F_{\varepsilon,h,Y}, \tilde{B}_{\varepsilon,t}^A} \trianglelefteq K_{P_{3,h}, \tilde{B}_{\varepsilon}^A}$. The thesis follows in a similar way, using the fact that, by Remark 90 and the fact that $\|P_{3,h}\|_{\mathcal{L}(X^3,X)} \leqslant \|P_3\|_{\mathcal{L}(X^3,X)}$ (where $P_3$ is the nonlinearity defined in equation (105) where all the cut-offs have been removed), we have $K_{P_3,h,\tilde{B}_\varepsilon^A} \trianglelefteq K_{P_3,\tilde{B}_\varepsilon^A}$, and the convergence $Z_{F_{\varepsilon,h,Y},\tilde{B}_{\varepsilon,t}^A}$ to $Z_{F_{\varepsilon,Y},\tilde{B}_{\varepsilon,t}^A}$ (i.e. the solution to equation (93)). $\square$

# Bibliography

**[1]** L. Accardi. An outline of quantum probability. Preprint 10.13140/rg.2.1.3078.3844, 2015.

**[2]** L. Accardi, A. Frigerio, and J. T. Lewis. Quantum stochastic processes. *Kyoto University. Research Institute for Mathematical Sciences. Publications*, 18(1):97–133, 1982. 10.2977/prims/1195184017.

**[3]** S. Albeverio, F. C. De Vecchi, and M. Gubinelli. Elliptic stochastic quantization. *Ann. Probab.*, 48(4):1693–1741, 2020. 10.1214/19-AOP1404.

**[4]** S. Albeverio, F. C. De Vecchi, and M. Gubinelli. The elliptic stochastic quantization of some two dimensional Euclidean QFTs. *Annales de l'Institut Henri Poincaré, Probabilités et Statistiques*, 57(4), Nov. 2021. 10.1214/20-AIHP1145.

**[5]** S. Albeverio and R. Høegh-Krohn. The Wightman axioms and the mass gap for strong interactions of exponential type in two-dimensional space-time. *Journal of Functional Analysis*, 16(1):39–82, May 1974. 10.1016/0022-1236(74)90070-6.

**[6]** S. Albeverio and S. Kusuoka. The invariant measure and the flow associated to the $\phi_3^4$-quantum field model. *Annali della Scuola Normale di Pisa - Classe di Scienze*, 2018. 10.2422/2036-2145.201809_008.

**[7]** S. Albeverio and S. Kusuoka. Construction of a non-Gaussian and rotation-invariant $\Phi^4$-measure and associated flow on $\mathbb{R}^3$ through stochastic quantization. *arXiv:2102.08040 [math-ph]*, Apr. 2021. 10.1007/BF01942330.

**[8]** D. Alpay, I. L. Paiva, and D. C. Struppa. Distribution spaces and a new construction of stochastic processes associated with the Grassmann algebra. *Journal of Mathematical Physics*, 60(1):013508, 21, 2019. 10.1063/1.5052010.

**[9]** H. Amann. Operator-valued Fourier multipliers, vector-valued Besov spaces, and applications. *Math. Nachr.*, 186:5–56, 1997. 10.1002/mana.3211860102.

**[10]** H. Amann. *Linear and quasilinear parabolic problems. Vol. II*, volume 106 of *Monographs in Mathematics*. Birkhäuser/Springer, Cham, 2019. 10.1007/978-3-030-11763-4.

**[11]** G. F. D. Angelis, G. Jona-Lasinio, and V. Sidoravicius. Berezin integrals and Poisson processes. *Journal of Physics A: Mathematical and General*, 31(1):289–308, Jan. 1998. 10.1088/0305-4470/31/1/026.

**[12]** D. Applebaum and R. L. Hudson. Fermion diffusions. *J. Math. Phys.*, 25(4):858–861, 1984. 10.1063/1.526236.

**[13]** D. B. Applebaum and R. L. Hudson. Fermion Itô's formula and stochastic evolutions. *Comm. Math. Phys.*, 96(4):473–496, 1984. 10.1063/1.526236.

**[14]** H. Araki. *Mathematical theory of quantum fields*, volume 101 Carow-Watamura of *International Series of Monographs on Physics*. Oxford University Press, New York, 1999. Translated from the 1993 Japanese original by Ursula 10.1088/0305-4470/28/2/004.

**[15]** J. C. Baez, I. E. Segal, and Z.-F. Zhou. *Introduction to algebraic and constructive quantum field theory*. Princeton Series in Physics. Princeton University Press, Princeton, NJ, 1992. 10.1515/9781400862504.

**[16]** D. Bahns, S. Doplicher, G. Morsella, and G. Piacitelli. Quantum spacetime and algebraic quantum field theory. In *Advances in algebraic quantum field theory*, Math. Phys. Stud., pages 289–329. Springer, Cham, 2015. 10.1007/978-3-319-21353-8_7.

**[17]** N. Barashkov and F. C. De Vecchi. Elliptic stochastic quantization of Sinh-Gordon QFT. *arXiv:2108.12664 [math-ph]*, Aug. 2021. 10.1016/0022-1236(74)90070-6.

**[18]** C. Barnett, R. F. Streater, and I. F. Wilde. The Itô-Clifford integral. *J. Functional Analysis*, 48(2):172–212, 1982. 10.1016/0022-1236(82)90066-0.

**[19]** C. Barnett, R. F. Streater, and I. F. Wilde. The Itô-Clifford integral. II. Stochastic differential equations. *J. London Math. Soc. (2)*, 27(2):373–384, 1983. 10.1112/jlms/s2-27.2.373.

**[20]** C. Barnett, R. F. Streater, and I. F. Wilde. The Itô-Clifford integral. III. The Markov property of solutions to stochastic differential equations. *Comm. Math. Phys.*, 89(1):13–17, 1983. 10.2977/prims/1195184017.

**[21]** G. A. Battle and L. Rosen. On the infinite volume limit of the strongly coupled Yukawa$_2$ model. *Journal of Mathematical Physics*, 22(4):770–776, Apr. 1981. 10.1063/1.524982.

[22] H. Baumgärtel. *Operator algebraic methods in quantum field theory*. Akademie Verlag, Berlin, 1995. A series of lectures 10.1007/978-3-319-21353-8_7.

[23] H. Baumgärtel and M. Wollenberg. *Causal nets of operator algebras*, volume 80 of *Mathematische Lehrbücher und Monographien, II. Abteilung: Mathematische Monographien [Mathematical Textbooks and Monographs, Part II: Mathematical Monographs]*. Akademie-Verlag, Berlin, 1992. Mathematical aspects of algebraic quantum field theory 10.1007/978-3-319-21353-8_7.

[24] V. P. Belavkin. Quantum stochastic calculus and quantum nonlinear filtering. *J. Multivariate Anal.*, 42(2):171–201, 1992. 10.1016/0047-259X(92)90042-E.

[25] G. Benfatto, P. Falco, and V. Mastropietro. Functional Integral Construction of the Massive Thirring model: Verification of Axioms and Massless Limit. *Communications in Mathematical Physics*, 273(1):67–118, July 2007. 10.1007/s00220-007-0254-y.

[26] G. Benfatto and G. Gallavotti. *Renormalization Group*. Princeton University Press, 1995. 10.1515/9780691221694.

[27] Y. M. Berezansky and Y. G. Kondratiev. *Spectral methods in infinite-dimensional analysis. Vol. 1*, volume 12/1 and D. V. Malyshev and revised by the authors. 10.1007/978-94-011-0509-5 of *Mathematical Physics and Applied Mathematics*. Kluwer Academic Publishers, Dordrecht, 1995. Translated from the 1988 Russian original by P. V. Malyshev 10.1007/978-94-011-0509-5.

[28] F. A. Berezin. *The method of second quantization*. Translated from the Russian by Nobumichi Mugibayashi and Alan Jeffrey. Pure and Applied Physics, Vol. 24. Academic Press, New York-London, 1966. 10.1007/978-88-7642-378-9.

[29] F. A. Berezin. *Introduction to superanalysis*, volume 9 of *Mathematical Physics and Applied Mathematics*. D. Reidel Publishing Co., Dordrecht, 1987. 10.1007/978-94-017-1963-6.

[30] L. Bertini, G. Jona-Lasinio, and C. Parrinello. Stochastic Quantization, Stochastic Calculus and Path Integrals: Selected Topics. *Progress of Theoretical Physics Supplement*, 111:83–113, Jan. 1993. 10.1143/PTPS.111.83.

[31] P. Biane. Calcul stochastique non-commutatif. In P. Bernard, editor, *Lectures on Probability Theory*, volume 1608, pages 1–96. Springer Berlin Heidelberg, Berlin, Heidelberg, 1995. 10.1007/BFb0095746.

[32] N. N. Bogolubov, A. A. Logunov, A. I. Oksak, and I. Todorov. *General Principles of Quantum Field Theory*. Mathematical Physics and Applied Mathematics. Springer Netherlands, 1990. 10.1007/978-94-009-0491-0.

[33] F. Bonetto and V. Mastropietro. Critical indices for the Yukawa2 quantum field theory. *Nuclear Physics B*, 497(1):541–554, July 1997. 10.1016/S0550-3213(97)00246-0.

[34] F. F. Bonsall and J. Duncan. *Complete normed algebras*. Springer-Verlag, New York-Heidelberg, 1973. Ergebnisse der Mathematik und ihrer Grenzgebiete, Band 80 10.1007/BF01059040.

[35] L. M. Borasi. *Probabilistic and differential geometric methods for relativistic and Euclidean Dirac and radiation fields*. PhD thesis, University of Bonn, Bonn, Germany, July 2019. 10.1007/978-3-0348-7917-0.

[36] M. Bożejko and R. Speicher. Interpolations between bosonic and fermionic relations given by generalized brownian motions. *Mathematische Zeitschrift*, 222(1):135–160, May 1996. 10.1007/BF02621861.

[37] O. Bratteli and D. W. Robinson. *Operator algebras and quantum statistical mechanics. 1*. Texts and Monographs in Physics. Springer-Verlag, New York, second edition, 1987. 10.1007/978-3-662-02520-8.

[38] O. Bratteli and D. W. Robinson. *Operator algebras and quantum statistical mechanics. 2*. Texts and Monographs in Physics. Springer-Verlag, Berlin, second edition, 1997. 10.1007/978-3-662-03444-6.

[39] F. E. Browder. Analyticity and partial differential equations. I. *Amer. J. Math.*, 84:666–710, 1962. 10.2307/2372872.

[40] E. R. Caianiello. *Combinatorics and Renormalization in Quantum Field Theory*, volume 38. Benjamin, 1973. 10.2977/prims/1195193913.

[41] E. Carlen and P. Krée. $L^p$ estimates on iterated stochastic integrals. *Ann. Probab.*, 19(1):354–368, 1991. 10.1016/0047-259X(92)90042-E.

[42] R. Catellier and K. Chouk. Paracontrolled distributions and the 3-dimensional stochastic quantization equation. *The Annals of Probability*, 46(5):2621–2679, 2018. 10.1214/17-AOP1235.

[43] F. Cipriani. Noncommutative potential theory: a survey. *J. Geom. Phys.*, 105:25–59, 2016. 10.1016/j.geomphys.2016.03.016.

[44] A. Connes. *Noncommutative geometry*. Academic Press, San Diego, 1994. 10.1007/BF01059040.

[45] A. Cooper and L. Rosen. The Weakly Coupled Yukawa$_2$ Field Theory: Cluster Expansion and Wightman Axioms. *Transactions of the American Mathematical Society*, 234(1):1–88, 1977. 10.2307/1997994.

[46] Y. Dabrowski. A free stochastic partial differential equation. *Annales de l'Institut Henri Poincaré, Probabilités et Statistiques*, 50(4), Nov. 2014. 10.1214/13-AIHP548.

[47] P. Damgaard and K. Tsokos. Stochastic quantization with fermions. *Nuclear Physics B*, 235(1):75–92, May 1984. 10.1016/0550-3213(84)90149-4.

[48] P. H. Damgaard and H. Hüffel. *Stochastic Quantization*. World Scientific, 1988.

[49] W. De Roeck and M. Salmhofer. Persistence of Exponential Decay and Spectral Gaps for Interacting Fermions. *Commun. Math. Phys.*, 365(2):773–796, Jan. 2019. 10.1007/s00220-018-3211-z.

[50] F. C. De Vecchi and M. Gubinelli. A note on supersymmetry and stochastic differential equations. *arXiv:1912.04830 [math-ph]*, Dec. 2019. arXiv: 1912.04830 10.1007/BFb0074481.

[51] B. De Witt. *Supermanifolds*. Cambridge Monographs on Mathematical Physics. Cambridge University Press, Cambridge, second edition, 1992. 10.1017/CBO9780511564000.

[52] S. Dirksen. Noncommutative stochastic integration through decoupling. *Journal of Mathematical Analysis and Applications*, 370(1):200–223, Oct. 2010. 10.1016/j.jmaa.2010.04.062.

[53] M. Disertori and V. Rivasseau. Continuous Constructive Fermionic Renormalization. *Annales Henri Poincaré*, 1(1):1–57, Feb. 2000. 10.1007/PL00000998.

[54] A. N. Efremov. Stochastic Quantization of Massive Fermions. *International Journal of Theoretical Physics*, 58(4):1150–1156, Apr. 2019. 10.1007/s10773-019-04006-w.

[55] G. G. Emch. *Mathematical and conceptual foundations of 20th-century physics*, volume 100 of *North-Holland Mathematics Studies*. North-Holland Publishing Co., Amsterdam, 1984. Notas de Matemá 10.1007/978-3-319-21353-8_7 .

[56] J. Feldman and H. Knörrer. E. Trubowitz, A two dimensional Fermi liquid. *Commun. Math. Phys*, 247:1–47, 49–111, 113–177, 179–194, 2004. 10.1007/s00220-018-3211-z.

[57] J. Feldman, H. Knörrer, and E. Trubowitz. *Fermionic functional integrals and the renormalization group*, volume 16 of *CRM Monograph Series*. American Mathematical Society, Providence, RI, 2002. 10.1007/BFb0078061.

[58] J. Feldman, H. Knörrer, and E. Trubowitz. Single scale analysis of many fermion systems. *Rev. Math. Phys.*, 15(09):949–994, 995–1038, 1039–1120, 1121–1169, Nov. 2003. 10.1142/S0129055X03001801.

[59] J. Feldman, H. Knörrer, and E. Trubowitz. Convergence of perturbation expansions in fermionic models. *Commun. Math. Phys.*, 247(1):195–242, May 2004. 10.1007/s00220-004-1039-1.

[60] J. Feldman, H. Knörrer, and E. Trubowitz. Convergence of Perturbation Expansions in Fermionic Models. Part 2: Overlapping Loops. *Commun. Math. Phys.*, 247(1):243–319, May 2004. 10.1007/s00220-004-1040-8.

[61] J. Feldman, J. Magnen, V. Rivasseau, and R. Sénéor. A renormalizable field theory: The massive Gross-Neveu model in two dimensions. *Commun.Math. Phys.*, 103(1):67–103, Mar. 1986. 10.1007/BF01464282.

[62] J. S. Feldman, H. Knörrer, H. Knorrer, and E. Trubowitz. *Fermionic Functional Integrals and the Renormalization Group*. Number 16. American Mathematical Soc., 2002. 10.1007/s00220-018-3211-z.

[63] K. O. Friedrichs. *Mathematical Aspects of the Quantum Theory of Fields*. Interscience Publishers, 1953. 10.2977/prims/1195193913.

[64] J. Fröhlich and K. Osterwalder. Is There a Euclidean Field Theory for Fermions. *Helv. Phys. Acta*, 47:781, 1975. 10.1016/0022-1236(80)90092-0.

[65] T. Fukai, H. Nakazato, I. Ohba, K. Okano, and Y. Yamanaka. Stochastic Quantization Method of Fermion Fields. *Progress of Theoretical Physics*, 69(5):1600–1616, May 1983. 10.1143/PTP.69.1600.

[66] K. Gawędzki and A. Kupiainen. Gross-Neveu model through convergent perturbation expansions. *Communications in Mathematical Physics*, 102(1):1–30, 1985. 10.1007/s00220-007-0254-y.

[67] I. Gelfand and M. Neumark. On the imbedding of normed rings into the ring of operators in Hilbert space. *Rec. Math. [Mat. Sbornik] N.S.*, 12(54):197–213, 1943. 10.2307/1968823.

[68] J. Glimm and A. Jaffe. The Yukawa$_2$ quantum field theory without cutoffs. *J. Functional Analysis*, 7:323–357, 1971. 10.1016/0022-1236(71)90039-5.

[69] J. Glimm and A. Jaffe. *Quantum physics*. Springer-Verlag, New York, second edition, 1987. A functional integral point of view 10.1007/978-1-4612-4728-9.

[70] M. Gordina. Stochastic differential equations on noncommutative $L^2$. In *Finite and infinite dimensional analysis in honor of Leonard Gross (New Orleans, LA, 2001)*, volume 317 of *Contemp. Math.*, pages 87–98. Amer. Math. Soc., Providence, RI, 2003. 10.1090/conm/317/05521.

[71] L. Gross. Abstract Wiener spaces. In *Proc. Fifth Berkeley Sympos. Math. Statist. and Probability (Berkeley, Calif., 1965/66), Vol. II: Contributions to Probability Theory, Part 1*, pages 31–42. Univ. California Press, Berkeley, Calif., 1967. 10.1143/PTPS.111.43.

[72] L. Gross. A noncommutative extension of the Perron-Frobenius theorem. *Bulletin of the American Mathematical Society*, 77(3):343–348, May 1971. 10.1090/S0002-9904-1971-12686-1.

[73] L. Gross. Existence and uniqueness of physical ground states. *Journal of Functional Analysis*, 10(1):52–109, May 1972. 10.1016/0022-1236(72)90057-2.

[74] L. Gross. Hypercontractivity and logarithmic Sobolev inequalities for the Clifford-Dirichlet form. *Duke Mathematical Journal*, 42(3):383–396, Sept. 1975. 10.1215/S0012-7094-75-04237-4.

[75] L. Gross. On the formula of Mathews and Salam. *Journal of Functional Analysis*, 25(2):162–209, June 1977. 10.1016/0022-1236(77)90039-8.

[76] M. Gubinelli and M. Hofmanová. Global Solutions to Elliptic and Parabolic $\phi^4$ Models in Euclidean Space. *Communications in Mathematical Physics*, 368(3):1201–1266, 2019. 10.1007/s00220-017-2997-4.

[77] M. Gubinelli and M. Hofmanová. A PDE Construction of the Euclidean $\Phi_3^4$ Quantum Field Theory. *Commun. Math. Phys.*, 384(1):1–75, 2021. 10.1007/s00220-021-04022-0.

[78] M. Gubinelli, P. Imkeller, and N. Perkowski. Paracontrolled distributions and singular PDEs. *Forum of Mathematics. Pi*, 3:e6, 75, 2015. 10.1017/fmp.2015.2.

[79] M. Gubinelli, H. Koch, and T. Oh. Paracontrolled approach to the three-dimensional stochastic nonlinear wave equation with quadratic nonlinearity. *arXiv:1811.07808 [math]*, Nov. 2018. Comment: 49 pages 10.1214/19-AOP1404.

[80] M. Gubinelli, H. Koch, and T. Oh. Renormalization of the two-dimensional stochastic nonlinear wave equations. *Transactions of the American Mathematical Society*, page 1, 2018. 10.1090/tran/7452.

[81] F. Guerra. Local algebras in Euclidean quantum field theory. In *Symposia Mathematica, Vol. XX (Convegno sulle Algebre C\* e loro Applicazioni in Fisica Teorica, Convegno sulla Teoria degli Operatori Indice e Teoria K, INDAM, Roma, 1974)*, pages 13–26. 1976. 10.1007/s11005-006-0124-0.

[82] R. Haag. *Local quantum physics*. Texts and Monographs in Physics. Springer-Verlag, Berlin, second edition, 1996. 10.1007/978-3-642-61458-3.

[83] Z. Haba and J. Kupsch. Supersymmetry in Euclidean Quantum Field Theory. *Fortschritte der Physik/Progress of Physics*, 43(1):41–66, 1995. 10.1002/prop.2190430103.

[84] M. Hairer. A theory of regularity structures. *Inventiones mathematicae*, 198(2):269–504, 2014. 10.1007/s00222-014-0505-4.

[85] K. Hepp. *Théorie de la renormalisation*. Cours donné à l'École Polytechnique, Paris. Lecture Notes in Physics, Vol. 2. Springer-Verlag, Berlin-New York, 1969. 10.1016/S0370-2693(96)01251-8.

[86] T. Hida, H.-H. Kuo, J. Potthoff, and W. Streit. *White Noise: An Infinite Dimensional Calculus*. Springer Science & Business Media, June 2013. 10.2307/1997994.

[87] A. Holevo. *Probabilistic and statistical aspects of quantum theory*, volume 1 of *Quaderni/Monographs*. Edizioni della Normale, Pisa, second edition, 2011. 10.1007/978-88-7642-378-9.

[88] E. P. Hsu. *Stochastic Analysis on Manifolds*. Graduate Studies in Mathematics 38. American Mathematical Society, 2002. 10.2977/prims/1195184017.

[89] R. L. Hudson and K. R. Parthasarathy. Unification of fermion and Boson stochastic calculus. *Communications In Mathematical Physics*, 104(3):457–470, Sept. 1986. 10.1007/BF01210951.

[90] I. A. Ignatyuk, V. A. Malyshev, and V. Sidoravičius. Convergence of a Method of the Stochastic Quantization II. *Theory of Probability & Its Applications*, 37(4):599–620, Jan. 1993. 10.1137/1137117.

[91] C. Itzykson and J. B. Zuber. *Quantum field theory*. McGraw-Hill International Book Co., New York, 1980. International Series in Pure and Applied Physics 10.1017/CBO9780511564000.

[92] V. D. Ivashchuk. Infinite-dimensional Grassmann-Banach algebras. *arXiv:math-ph/0009006*, Sept. 2000. 10.1016/j.bulsci.2006.05.007.

[93] A. Jadczyk and K. Pilch. Superspaces and supersymmetries. *Commun.Math. Phys.*, 78(3):373–390, Jan. 1981. 10.1007/BF01942330.

[94] S. Janson. *Gaussian Hilbert Spaces*. Cambridge University Press, June 1997. 10.2977/prims/1195184017.

[95] G. Jona-Lasinio. Stochastic quantization: a new domain for stochastic analysis. In *Proceedings of the 1st World Congress of the Bernoulli Society, Vol. 1 (Tashkent, 1986)*, pages 535–546. VNU Sci. Press, Utrecht, 1987. 10.1007/BF02099877.

[96] G. Jona-Lasinio and P. K. Mitter. On the stochastic quantization of field theory. *Communications in Mathematical Physics (1965-1997)*, 101(3):409–436, 1985. 10.1143/PTPS.111.43.

[97] G. Jona-Lasinio and P. K. Mitter. Large deviation estimates in the stochastic quantization of $\phi_2^4$. *Communications in Mathematical Physics*, 130(1):111–121, May 1990. 10.1007/BF02099877.

[98] R. Jost. *The general theory of quantized fields*, volume 1960 of *Mark Kac, editor. Lectures in Applied Mathematics (Proceedings of the Summer Seminar, Boulder, Colorado*. American Mathematical Society, Providence, R.I., 1965. 10.2307/1997994.

[99] R. V. Kadison and J. R. Ringrose. *Fundamentals of the theory of operator algebras. Vol. I*, volume 15 Reprint of the 1983 original of *Graduate Studies in Mathematics*. American Mathematical Society, Providence, RI, 1997. Elementary theory, 10.1090/gsm/016/01.

[100] R. V. Kadison and J. R. Ringrose. *Fundamentals of the theory of operator algebras. Vol. II*, volume 16 Corrected reprint of the 1986 original. 10.1090/gsm/016/01 of *Graduate Studies in Mathematics*.

American Mathematical Society, Providence, RI, 1997. Advanced theory, 10.1090/gsm/016/01.

[101] Y. Kakudo, Y. Taguchi, A. Tanaka, and K. Yamamoto. Gauge-independent calculation of S-matrix elements in quantum electrodynamics. *Progress of Theoretical Physics (Kyoto)*, 69(4):1225–1235, 1983. 10.1007/978-88-7642-378-9.

[102] A. Y. Khrennikov. Equations on a superspace. *Mathematics of the USSR-Izvestiya*, 36(3):597, 1991. Publisher: IOP Publishing. 10.1070/IM1991v036n03ABEH002036.

[103] A. Klein. Supersymmetry and a two-dimensional reduction in random phenomena. In *Quantum probability and applications, II (Heidelberg, 1984)*, volume 1136 of *Lecture Notes in Math.*, pages 306–317. Springer, Berlin, 1985. 10.1007/BFb0074481.

[104] A. Klein and J. Fernando Perez. Supersymmetry and dimensional reduction: A non-perturbative proof. *Physics Letters B*, 125(6):473–475, June 1983. 10.1016/0370-2693(83)91329-1.

[105] B. Kümmerer. Survey on a theory of noncommutative stationary Markov processes. In *Quantum probability and applications, III (Oberwolfach, 1987)*, volume 1303 of *Lecture Notes in Math.*, pages 154–182. Springer, Berlin, 1988. 10.1007/BFb0078061.

[106] A. Kupiainen. Renormalization Group and Stochastic PDEs. *Annales Henri Poincaré*, 17(3):497–535, 2016. 10.1007/s00023-015-0408-y.

[107] J. Kupsch. Fermionic and supersymmetric stochastic processes. *Journal of Geometry and Physics*, 11(1):507–516, June 1993. 10.1016/0393-0440(93)90074-O.

[108] O. E. Lanford III. *Construction of quantum fields interacting by a cutoff Yukawa coupling*. PhD thesis, Princeton university, jan 1966. 10.1007/BF01654299.

[109] Y. Le Jan. Temps local et superchamp. In *Séminaire de Probabilités, XXI*, volume 1247 of *Lecture Notes in Math.*, pages 176–190. Springer, Berlin, 1987. 10.1007/BFb0077633.

[110] D. Lehmann. A probabilistic approach to Euclidean Dirac fields. *Journal of Mathematical Physics*, 32(8):2158–2166, Aug. 1991. Publisher: American Institute of Physics. 10.1063/1.529189.

[111] D. A. Leites. Introduction to the theory of supermanifolds. *Uspekhi Mat. Nauk*, 35(1(211)):3–57, 255, 1980. 10.1142/9789812708854.

[112] S. Leppard and A. Rogers. A Feynman-Kac formula for anticommuting Brownian motion. *Journal of Physics A: Mathematical and General*, 34(3):555–568, Jan. 2001. 10.1088/0305-4470/34/3/315.

[113] A. Lesniewski. Effective action for the Yukawa$_2$ quantum field theory. *Communications in Mathematical Physics*, 108(3):437–467, Sept. 1987. 10.1007/BF01212319.

[114] J. Magnen and R. Sénéor. The Wightman axioms for the weakly coupled Yukawa model in two dimensions. *Comm. Math. Phys.*, 51(3):297–313, 1976. 10.1016/0003-4916(71)90243-0.

[115] J. Magnen and R. Sénéor. Yukawa quantum field theory in three dimensions ($Y_3$). In *Third International Conference on Collective Phenomena (Moscow, 1978)*, volume 337 of *Ann. New York Acad. Sci.*, pages 13–43. New York Acad. Sci., New York, 1980. 10.1016/0003-4916(71)90243-0.

[116] A. I. Markushevich. *Theory of functions of a complex variable. Vol. II*. Revised English edition translated and edited by Richard A. Silverman. Prentice-Hall, Inc., Englewood Cliffs, N.J., 1965. 10.1007/978-1-4612-0907-2.

[117] V. Mastropietro. Schwinger functions in Thirring and Luttinger models. *Nuov Cim B*, 108(10):1095–1107, Oct. 1993. 10.1007/BF02827305.

[118] V. Mastropietro. *Non-Perturbative Renormalization*. World Scientific Publishing Co Pte Ltd, Hackensack, NJ, 2008. 10.1007/BF01464282.

[119] P. T. Matthews and A. Salam. Propagators of quantized field. *Il Nuovo Cimento (1955-1965)*, 2(1):120–134, July 1955. 10.1007/BF02856011.

[120] O. A. McBryan. Volume dependence of Schwinger function in the Yukawa$_2$ quantum field theory. *Communications in Mathematical Physics*, 45(3):279–294, 1975. Publisher: Springer-Verlag 10.1016/S0370-2693(96)01251-8.

[121] P.-A. Meyer. *Quantum probability for probabilists*, volume 1538 of *Lecture Notes in Mathematics*. Springer-Verlag, Berlin, 1993. 10.1007/978-3-662-21558-6.

[122] A. Moinat and H. Weber. Space-time localisation for the dynamic $\phi_3^4$ model. *arXiv:1811.05764*, Nov. 2018. Comment: 27 pages 10.1007/s00220-021-04022-0.

[123] J.-C. Mourrat and H. Weber. The dynamic $\Phi_3^4$ model comes down from infinity. *Comm. Math. Phys.*, 356(3):673–753, 2017. 10.1007/s00220-017-2997-4.

[124] J.-C. Mourrat and H. Weber. Global well-posedness of the dynamic $\Phi^4$ model in the plane. *The Annals of Probability*, 45(4):2398–2476, 2017. 10.1214/16-AOP1116.

[125] M. A. Naimark. *Normed algebras*. Wolters-Noordhoff Publishing, Groningen, third, wolters-noordhoff series of monographs and textbooks on pure and applied mathematics edition, 1972. Translated from the second Russian edition by Leo F. Boron, 10.1515/9781400862504.

[126] T. Nakano. Quantum Field Theory in Terms of Euclidean Parameters. *Progress of Theoretical Physics*, 21(2):241–259, Feb. 1959. Publisher: Oxford Academic. 10.1143/PTP.21.241.

**[127]** E. Nelson. Construction of quantum fields from Markoff fields. *Journal of Functional Analysis*, 12(1):97–112, Jan. 1973. 10.1016/0022-1236(73)90091-8.

**[128]** E. Nelson. Notes on non-commutative integration. *Journal of Functional Analysis*, 15(2):103–116, Feb. 1974. 10.1016/0022-1236(74)90014-7.

**[129]** E. P. Osipov. Quantum interaction $\phi_4^4$, the construction of quantum field defined as a bilinear form. *Journal of Mathematical Physics*, 41(2):759–786, 2000. 10.1063/1.533162.

**[130]** K. Osterwalder. Euclidean fermi fields. In *Constructive Quantum Field Theory*, Lecture Notes in Physics, pages 326–331. Springer, Berlin, Heidelberg, 1973. 10.1007/BFb0113094.

**[131]** K. Osterwalder and R. Schrader. Euclidean Fermi fields and a Feynman-Kac formula for Boson-Fermions models. *Helvetica Physica Acta*, 46:277–302, 1973. 10.1088/0305-4470/25/22/027.

**[132]** J. Palmer. Euclidean Fermi fields. *Journal of Functional Analysis*, 36(3):287–312, May 1980. 10.1016/0022-1236(80)90092-0.

**[133]** G. Parisi and Y. S. Wu. Perturbation theory without gauge fixing. *Scientia Sinica. Zhongguo Kexue*, 24(4):483–496, 1981. 10.1007/BF01074107.

**[134]** K. R. Parthasarathy. *An Introduction to Quantum Stochastic Calculus*. Springer Science & Business Media, 1992. 10.1090/tran/7452.

**[135]** W. A. d. S. Pedra and M. Salmhofer. Determinant Bounds and the Matsubara UV Problem of Many-Fermion Systems. *Communications in Mathematical Physics*, 282(3):797–818, Sept. 2008. 10.1007/s00220-008-0463-z.

**[136]** V. G. Pestov. General construction of Banach-Grassmann algebras. *Atti accademia nazionale Lincei classe Scienze fisiche matematiche naturali. Rendiconti Lincei. Matematica e applicationi. Serie 9*, 3(3):223–231, 1992. 10.1016/j.bulsci.2006.05.007.

**[137]** V. G. Pestov. Analysis on superspace: an overview. *Bulletin of the Australian Mathematical Society*, 50(1):135–165, Aug. 1994. Publisher: Cambridge University Press. 10.1017/S0004972700009643.

**[138]** A. Prastaro and T. M. Rassias. *Geometry in Partial Differential Equations*. World Scientific, 1994. 10.1007/BF01654299.

**[139]** W. Ramasinghe. Exterior algebra of a Banach space. *Bulletin des Sciences Mathématiques*, 131(3):291–324, Apr. 2007. 10.1016/j.bulsci.2006.05.007.

**[140]** P. . Renouard. Borel analyticity and summability of Schwinger functions of the two-dimensional Yukawa model II Adiabatic limit. *Annales de l'Institut Henri Poincare Section A, Physique Theorique*, 31(3):235–318, 1979. 10.1007/BF02827305.

**[141]** V. Rivasseau. *From Perturbative to Constructive Renormalization*. Princeton University Press, Princeton, N.J, 2 edition edition, May 1991. 10.1088/0305-4470/25/22/027.

**[142]** A. Rogers. Graded manifolds, supermanifolds and infinite-dimensional Grassmann algebras. *Communications in Mathematical Physics*, 105(3):375–384, Sept. 1986. 10.1007/BF01205932.

**[143]** A. Rogers. Stochastic calculus in superspace. I. Supersymmetric Hamiltonians. *Journal of Physics A: Mathematical and General*, 25(2):447–468, Jan. 1992. 10.1088/0305-4470/25/2/024.

**[144]** A. Rogers. Stochastic calculus in superspace. II. Differential forms, supermanifolds and the Atiyah-Singer index theorem. *Journal of Physics A: Mathematical and General*, 25(22):6043–6062, Nov. 1992. 10.1088/0305-4470/25/22/027.

**[145]** A. Rogers. Path integration, anticommuting variables, and supersymmetry. *Journal of Mathematical Physics*, 36(5):2531–2545, May 1995. 10.1063/1.531049.

**[146]** A. Rogers. *Supermanifolds*. World Scientific Publishing Co. Pte. Ltd., Hackensack, NJ, 2007. 10.1142/9789812708854.

**[147]** J. Rosenberg. Noncommutative variations on Laplace's equation. *Analysis & PDE*, 1(1):95–114, Oct. 2008. Publisher: Mathematical Sciences Publishers. 10.2140/apde.2008.1.95.

**[148]** R. A. Ryan. *Introduction to tensor products of Banach spaces*. Springer Monographs in Mathematics. Springer-Verlag London, Ltd., London, 2002. 10.1007/978-1-4471-3903-4.

**[149]** M. Salmhofer. *Renormalization: An Introduction*. Springer, Berlin ; New York, 1st corrected ed. 1999, corr. 2nd printing 2007 edition edition, Sept. 2007. 10.1007/BFb0078061.

**[150]** M. Salmhofer. Clustering of Fermionic Truncated Expectation Values Via Functional Integration. *J Stat Phys*, 134(5):941–952, Mar. 2009. 10.1007/s10955-009-9698-0.

**[151]** M. Salmhofer and C. Wieczerkowski. Positivity and Convergence in Fermionic Quantum Field Theory. *Journal of Statistical Physics*, 99(1):557–586, Apr. 2000. 10.1023/A:1018661110470.

**[152]** M. Salmhofer and C. Wieczerkowski. Construction of the renormalized $GN_{2-\epsilon}$ trajectory. *Mathematical Physics Electronic Journal*, pages 1–19, mar 2002. 10.1142/9789812777874_0007.

**[153]** A. Y. Savin and B. Y. Sternin. Noncommutative elliptic theory. Examples. *Proceedings of the Steklov Institute of Mathematics*, 271(1):193–211, Dec. 2010. 10.1134/S0081543810040152.

**[154]** J. L. Schiff. *Normal families*. Universitext. Springer-Verlag, New York, 1993. 10.1007/978-1-4612-0907-2.

[155] H.-J. Schmeisser and H. Triebel. *Topics in Fourier analysis and function spaces*. A Wiley-Interscience Publication. John Wiley & Sons, Ltd., Chichester, 1987. 10.1007/978-1-4612-0907-2.

[156] K. Schmüdgen. *Unbounded operator algebras and representation theory*, volume 37 of *Operator Theory: Advances and Applications*. Birkhäuser Verlag, Basel, 1990. 10.1007/978-3-0348-7469-4.

[157] J. Schwinger. On the euclidean structure of relativistic field theory. *Proceedings of the National Academy of Sciences of the United States of America*, 44(9):956–965, Sept. 1958. 10.2307/1969729.

[158] J. Schwinger. Euclidean Quantum Electrodynamics. *Physical Review*, 115(3):721–731, Aug. 1959. Publisher: American Physical Society. 10.1103/PhysRev.115.721.

[159] M. Schürman. Quantum q-white noise and a q-central limit theorem. *Communications in Mathematical Physics*, 140(3):589–615, Oct. 1991. 10.1007/BF02099136.

[160] I. E. Segal. A non-commutative extension of abstract integration. *Ann. of Math. (2)*, 57:401–457, 1953. 10.2307/1969729.

[161] I. E. Segal. A Non-Commutative Extension of Abstract Integration. *Annals of Mathematics*, 57(3):401–457, 1953. 10.2307/1969729.

[162] I. E. Segal. *Mathematical problems of relativistic physics*, volume 1960 of *With an appendix by George W. Mackey. Lectures in Applied Mathematics (proceedings of the Summer Seminar, Boulder, Colorado*. American Mathematical Society, Providence, R.I., 1963. 10.1007/978-3-319-21353-8_7.

[163] E. Seiler. Schwinger functions for the Yukawa model in two dimensions with space-time cutoff. *Communications in Mathematical Physics*, 42(2):163–182, 1975. Publisher: Springer-Verlag 10.1007/BF01654299.

[164] V. V. Shcherbakov. Elements of stochastic analysis for the case of Grassmann variables. I. Grassmann stochastic integrals and random processes. *Theoretical and Mathematical Physics*, 96(1):792–800, July 1993. 10.1007/BF01074107.

[165] V. V. Shcherbakov. Elements of stochastic analysis for the case of Grassmann variables. II. Stochastic partial differential equations for Grassmann processes. *Theoretical and Mathematical Physics*, 97(2):1229–1235, Nov. 1993. 10.1007/BF01016868.

[166] V. V. Shcherbakov. Elements of stochastic analysis for the case of Grassmann variables. III. Correlation functions. *Theoretical and Mathematical Physics*, 97(3):1323–1332, Dec. 1993. 10.1007/BF01015761.

[167] B. Simon. $P(\phi)_2$ *Euclidean (Quantum) Field Theory*. Princeton University Press, Princeton, N.J, 1974. 10.1088/0305-4470/25/22/027.

[168] K. B. Sinha and D. Goswami. *Quantum stochastic processes and noncommutative geometry*, volume 169 of *Cambridge Tracts in Mathematics*. Cambridge University Press, Cambridge, 2007. 10.1017/CBO9780511618529.

[169] R. F. Streater. Classical and quantum probability. *Journal of Mathematical Physics*, 41(6):3556–3603, June 2000. 10.1063/1.533322.

[170] R. F. Streater and A. S. Wightman. *PCT, spin and statistics, and all that*. Princeton Landmarks in Physics. Princeton University Press, Princeton, NJ, 2000. 10.2307/1997994.

[171] F. Strocchi. *An Introduction to Non-Perturbative Foundations of Quantum Field Theory*. OUP Oxford, Oxford, Feb. 2013. 10.1007/BF01059040.

[172] S. J. Summers. A perspective on constructive quantum field theory. *arXiv preprint arXiv:1203.3991*, 2012. 10.1007/BFb0078061.

[173] K. Symanzik. Euclidean quantum field theory. In R. Jost, editor, *Local quantum field theory (Varenna lectures)*, page 152–226. Academic Press, New York, 1969. 10.1007/978-94-009-0491-0.

[174] F. Trèves. *Topological vector spaces, distributions and kernels*. Academic Press, New York-London, 1967. 10.1007/BFb0078061.

[175] H. Triebel. *Theory of function spaces. III*, volume 100 of *Monographs in Mathematics*. Birkhäuser Verlag, Basel, 2006. 10.1007/978-3-030-11763-4.

[176] J. van der Hoeven. *Majorants for formal power series*. Technical Report 2003-15. Université de Paris-Sud. Département de Mathématique, Orsay, 2003. 10.1007/978-1-4612-0907-2.

[177] B. van der Waerden. *Algebra: Volume I*. Algebra : based in part on lectures by E. Artin and E. Noether. Springer New York, 2003. 10.1143/PTPS.111.43.

[178] B. van der Waerden. *Algebra: Volume II*. Algebra : based in part on lectures by E. Artin and E. Noether. Springer New York, 2003. 10.1143/PTPS.111.43.

[179] P. van Nieuwenhuizen and A. Waldron. On Euclidean spinors and Wick rotations. *Physics Letters B*, 389(1):29–36, Dec. 1996. 10.1016/S0370-2693(96)01251-8.

[180] J. von Neumann. On rings of operators. III. *Ann. of Math. (2)*, 41:94–161, 1940. 10.2307/1968823.

[181] W. von Waldenfels. Non-commutative algebraic central limit theorems. In H. Heyer, editor, *Probability Measures on Groups VIII*, Lecture Notes in Mathematics, pages 174–202, Berlin, Heidelberg, 1986. Springer. 10.1007/BFb0077184.

**[182]** D. N. Williams. Euclidean Fermi fields with a hermitean Feynman-Kac-Nelson formula. I. *Communications in Mathematical Physics*, 38(1):65–80, Mar. 1974. 10.1007/BF01651549.

**[183]** X. Xiong, Q. Xu, and Z. Yin. *Sobolev, Besov and Triebel-Lizorkin Spaces on Quantum Tori*, volume 252 of *Memoirs of the American Mathematical Society*. American Mathematical Society, Mar. 2018. 10.1090/memo/1203.

**[184]** S.-S. Xue and T.-c. Hsien. Stochastic quantization of fermions on lattice. *Chinese Physics Letters*, 2(10):474–476, Oct. 1985. 10.1088/0256-307X/2/10/012.

**[185]** W. a. Żelazko. *Banach algebras*. Elsevier Publishing Co., Amsterdam-London-New York; PWN–Polish Scientific Publishers, Warsaw, 1973. Translated from the Polish by Marcin E. Kuczma 10.1007/BFb0078061.

**[186]** J. Zinn-Justin. *Quantum field theory and critical phenomena*, volume 85 of *International Series of Monographs on Physics*. The Clarendon Press, Oxford University Press, New York, second edition, 1993. Oxford Science Publications 10.1007/978-88-7642-378-9.

**[187]** Y. M. Zinoviev. Equivalence of Euclidean and Wightman field theories. *Communications in Mathematical Physics*, 174(1):1–27, 1995. Publisher: Springer-Verlag 10.1088/0305-4470/28/2/004.